\documentclass[mnsc,nonblindrev]{informs3a}

\OneAndAHalfSpacedXI 

\usepackage{natbib}
 \bibpunct[, ]{(}{)}{,}{a}{}{,}%
 \def\bibfont{\small}%
\usepackage{booktabs}

\renewcommand{\bar}{\overline}

\newcommand{\R}{ \mathbb{R} }
\renewcommand{\P}{ \mathbb{P}}

\newcommand{\E}{ \mathbb{E} }

\newcommand{\bw}{{\bm w}}
\newcommand{\bs}{{\bm s}}
\newcommand{\by}{{\bm y}}
\newcommand{\bx}{{\bm x}}
\newcommand{\bq}{{\bm q}}
\newcommand{\bp}{{\bm p}}
\newcommand{\bd}{{\bm d}}
\newcommand{\bD}{{\bm D}}
\newcommand{\bh}{{\bm h}}
\newcommand{\bb}{{\bm b}} 
\newcommand{\tby}{{\tilde{\bm y}}}
\newcommand{\tbw}{{\tilde{\bm w}}}
\newcommand{\ty}{{\tilde{y}}}
\newcommand{\tQ}{{\tilde{Q}}}
 
\usepackage{parskip} 

\def\eq#1{Eq.~(#1)}

\usepackage{makecell}
\usepackage{braket}
\usepackage{algorithm}
\usepackage{algorithmicx}
\usepackage[noend]{algpseudocode}
\def\1{(\mathrm{\uppercase\expandafter{\romannumeral1}})}
\def\2{(\mathrm{\uppercase\expandafter{\romannumeral2}})}
\def\O{\mathcal O}
\newcommand{\mh}{\mathcal{H}}
\RequirePackage[colorlinks,citecolor=blue,linkcolor=blue,urlcolor=blue]{hyperref}
\usepackage{enumitem}
\DeclareMathAlphabet{\mathscr}{OT1}{pzc}{m}{it}

\def\O{\mathcal O}

\usepackage{bm}
\usepackage{hyperref}

\TheoremsNumberedThrough     
\ECRepeatTheorems

\EquationsNumberedThrough    

\MANUSCRIPTNO{MS-DSC-23-00920.R2}
\usepackage{xcolor}
\definecolor{DSgray}{cmyk}{0,1,0,0}

\usepackage{mathabx}
\renewcommand{\hat}{\widehat}
\renewcommand{\tilde}{\widetilde}

\renewcommand{\bibfont}{\footnotesize\linespread{1}\selectfont}

\begin{document}



\RUNTITLE{Minibatch-SGD-Based Meta-Policy for Inventory Systems}

\TITLE{A Minibatch-SGD-Based  Learning  Meta-Policy for Inventory Systems with Myopic Optimal Policy}

\ARTICLEAUTHORS{%
 \AUTHOR{Jiameng Lyu\footnotemark[1]}
 \AFF{Yau Mathematical Sciences Center \& Department of Mathematical  Sciences, Tsinghua University, Beijing 100084, China, \EMAIL{lvjm21@mails.tsinghua.edu.cn}}
 \AUTHOR{Jinxing Xie\footnotemark[1]}
\AFF{Department of Mathematical Sciences, Tsinghua University, Beijing 100084, China, \EMAIL{xiejx@tsinghua.edu.cn}}
 \AUTHOR{Shilin Yuan\footnotemark[1]}
\AFF{Department of Mathematical Sciences, Tsinghua University, Beijing 100084, China, \EMAIL{yuansl21@mails.tsinghua.edu.cn}}
 \AUTHOR{Yuan Zhou\footnotemark[1]}
 \AFF{Yau Mathematical Sciences Center \& Department of Mathematical  Sciences, Tsinghua University, Beijing 100084, China, \EMAIL{yuan-zhou@tsinghua.edu.cn}}
 } 
 \renewcommand{\thefootnote}{\fnsymbol{footnote}}
\footnotetext[1]{Author names listed in alphabetical order.}
\renewcommand{\thefootnote}{\arabic{footnote}}

\ABSTRACT{Stochastic gradient descent (SGD) has proven effective in solving many inventory control problems with demand learning. However,  it often faces the pitfall of an infeasible target inventory level that is lower than the current inventory level. Several recent works (e.g., \cite{huh2009nonparametric,shi2016nonparametric}) are successful to resolve this issue in various inventory systems. However, their techniques are rather sophisticated and difficult to be applied to more complicated scenarios such as multi-product and multi-constraint inventory systems. 
 
In this paper, we address the infeasible-target-inventory-level issue from a new technical perspective -- we propose a novel minibatch-SGD-based meta-policy. Our meta-policy is flexible enough to be applied
to a general inventory systems framework covering a wide range of inventory management problems with myopic clairvoyant optimal policy. 
By devising the optimal minibatch scheme,  our meta-policy achieves  a regret bound of $\mathcal{O}(\sqrt{T})$ for the general convex case and $\mathcal{O}(\log T)$  for the strongly convex case. 
To demonstrate the power and flexibility of our meta-policy, we apply it to three important inventory control problems: multi-product and multi-constraint systems, multi-echelon serial systems, and one-warehouse and multi-store systems by carefully designing application-specific subroutines.
{We also conduct extensive numerical experiments to demonstrate that our meta-policy enjoys competitive regret performance, high computational efficiency, and low variances among a wide range of applications.}
}

\maketitle

 \section{Introduction}\label{sec:intro}
 Inventory control, a fundamental problem in operations management, has received tremendous attention from both academia and industry.
Classic research assumes the firm has full information of the demand distribution, and mainly focuses on capturing the structure of optimal solutions and finding them efficiently. In practice, however, the firm manager may not know the demand distribution a priori. Motivated by this practical consideration, some recent research on inventory control aims to jointly learn the demand distribution and optimize the inventory control decisions based on the observed data for different inventory systems. This dynamic setting lies in the intersection of machine learning and traditional operations management, and becomes increasingly important in recent years.

Gradient-based methods have proved very useful in online learning and operations management \citep{burnetas2000adaptive,besbes2015non}. The stochastic gradient descent (SGD) algorithm further enjoys much flexibility --- it only needs a noisy estimation of the cost function's gradient. Such estimation is usually available in many inventory systems and it is tempting to leverage the power and flexibility of the SGD algorithm and apply its principle in online inventory control tasks. However, there is a potential pitfall --- SGD-based algorithms might often instruct an infeasible target inventory level that is lower than the current inventory level. Such a generic issue occurs in almost all dynamic inventory control problems with demand learning. 
Although keeping the current inventory and ordering nothing may seem like a straightforward fix to the issue, the regret analysis is considerably more challenging. Typically, this involves upper bounding the total ordering differences incurred by the fix.
The seminal work \citep{huh2009nonparametric} proposed a novel technical method to analyze such differences based on the queueing theory, and established the optimal regret bound for an SGD-based algorithm in a single-product inventory system. This queueing-theory-based analytic method turned out to be extremely successful --- \cite{shi2016nonparametric} extended the analysis to a multi-product inventory system and many more follow-up works \citep{chen2020optimal,yuan2021marrying,ding2021feature,yang2022non} recently emerged to apply the technique to a series of popular inventory control scenarios.

On the other hand, however, the analytic method by \cite{huh2009nonparametric} has its own limitations. The queueing-theory-based argument is complex in nature and this complexity hinders the further application of the method to even more complicated inventory systems. For example, the analysis in \cite{shi2016nonparametric} only works for a single warehouse constraint, and it is not clear how to adapt the argument to multiple warehouse constraints, which is common in practice. As we will soon demonstrate in the paper, there are many more complex scenarios that call for a simpler and more powerful approach to achieve the provably effective learning goal.

\subsection{Our Contributions}

The main contribution of this paper is a general meta-policy for inventory systems with myopic optimal policy to learn and make optimal inventory decisions. Our meta-policy is simple, {powerful} and flexible so that it can be applied to a wide class of inventory systems with much more complex decision processes compared to existing works.

\noindent \underline{\textbf{The minibatch-SGD-based meta-policy.}} Recall that the main difficulty faced by traditional gradient-based approaches is the possibility of a lower target inventory level compared to the current inventory level, which cannot be realized by a single-period ordering instruction. While existing methods \citep{huh2009nonparametric,shi2016nonparametric} spend much effort to bound the \emph{magnitude} of the discrepancy between the two inventory levels, we aim at controlling the \emph{frequency} of the lower target inventory levels. 

The key observation that motivates our new meta-policy is as follows. Suppose we have a constant target inventory level for a consecutive epoch of periods. During this epoch, if the low-target-inventory-level issue ever happens, we only have to deal with this trouble at the very beginning of the epoch. Indeed, once the current inventory level gets no higher than the constant target level, it will stay equal to or below the target level throughout the rest of the epoch. In other words, the low-target-inventory-level issue might only be triggered when we \emph{switch} the target level.

According to the above observation, we denote a consecutive epoch of constant target inventory level by a \emph{minibatch}, and describe the high-level idea of our minibatch-SGD-based meta-policy as follows. We would like to limit the number of minibatches in the learning process so that we can control the frequency of lower target inventory levels. 
We need to bound the regret incurred in two scenarios: one where the low-target-inventory-level issue occurs and the other where it does not.
If we could achieve these two bounds, we derive a concrete regret upper bound. In light of this, we need to work on the following technical problems:

\begin{enumerate}[nolistsep]
\item \emph{What is the optimal minibatch scheme?} In particular, we need to carefully decide the total number of minibatches and the length of each minibatch. Depending on the convexity types of the inventory cost function (general convex or strongly convex), we will design different minibatch schemes.

\item  {\it How should we address the low-target-inventory-level issue?} This happens with a low frequency thanks to the minibatch scheme. However, we still need a \emph{transition solver} to control the total number of periods when such an issue lasts. The transition solver may be problem specific, but can be directly plugged into our flexible meta-policy as long as it meets the simple interface requirements.

\item 
{\it How should we design the gradient estimator? }
We need to design a problem-specific estimator for the gradient of the inventory cost function, and prove it is ``well-behaved'' (e.g., unbiased and has bounded variance). This is also needed by most existing SGD-based inventory algorithms. However, since our meta-policy can be applied to more complex scenarios, we need additional techniques in this step, as will be illustrated in our applications.
\end{enumerate}

{Besides the above-mentioned problem-specific technical problems, we need the following conditions which generally holds for a wide range of inventory management problems:
\begin{enumerate}[nolistsep]
    \item \textit{Myopic policy optimality.} The inventory control problem has a myopic optimal policy.
    \item \textit{Smooth, convex (strongly convex) cost function.} The associated inventory cost function should be smooth and convex (strongly convex). 
    \item \textit{Convex compact domain.} We (only) require the domain (constraint set) of the target inventory level to be convex and compact in the analysis of our meta-policy. 
\end{enumerate}}
In Sections~\ref{sec: meta policy} and~\ref{sec:regret of meta policy}, we present the details of our meta-policy. Under the above sufficient conditions, we prove that, given the desired transition solver and gradient estimator, our meta-policy achieves a regret bound of $\mathcal{O}(\sqrt{T})$ for the general convex case and $\mathcal{O}(\log T)$ for the strongly convex case.
We then apply our meta-policy to several applications to illustrate its power and flexibility.
In the first application, we work on multi-product inventory systems with multiple constraints and substantially extend the results of a series of existing works, {which shows the power of our meta-policy}. {We then consider even more complex scenarios to show its flexibility with the help of additional techniques} -- in the second application, the decision made at each time period is multi-dimensional and involves correlations among the different dimensions; in the third application, the decision involves a second-stage linear program and we need additional techniques to prove the smoothness property of the cost function and build the desired gradient estimator. Below is a more detailed overview.

\begin{table}[t]
\caption{Comparison between our results for Application I and the related works.\\ $\mathcal{R}_{\mathrm{cvx}}$:  regret under the general convex case; $\mathcal{R}_{\mathrm{str}}$:  regret under the strongly convex case. }
\label{tab:app1}
\renewcommand{\arraystretch}{1.3}
\begin{center}
\scalebox{0.8}{%
\begin{tabular}{l|l|l|l|l|l}
\hline
           & \makecell[l]{Number of \\ Product(s)} & \makecell[l]{Number of\\ Constraint(s)} & \makecell[l]{Assumption on Product-level\\ Demand Independence} & $\mathcal{R}_{\mathrm{cvx}}$ & $\mathcal{R}_{\mathrm{str}}$ \\ \hline
\cite{huh2009nonparametric}      & Single         & Single            & \textbf{---}                  & $\mathcal{O}(\sqrt{T})$ & $\mathcal{O}(\log{T})$  \\ \hline
\cite{shi2016nonparametric}       & \textbf{Multiple}       & Single            & Required                & $\mathcal{O}(\sqrt{T})$ & $\mathcal{O}(\log{T})$  \\ \hline
Our results & \textbf{Multiple}      & \textbf{Multiple}          & \textbf{Not Required}            & $\mathcal{O}(\sqrt{T})$ & $\mathcal{O}(\log{T})$  \\ \hline
\end{tabular}%
}
\end{center}
\end{table}
\noindent \underline{\textbf{Application I: multi-product and multi-constraint  inventory system.}}  We first apply our meta-policy to multi-product and multi-constraint inventory systems, where there have been a series of prominent related works. Compared to the state-of-the-art algorithms in the literature, our results admit multiple constraints on the inventory capacities and may deal with demand correlations among products (i.e., no need to assume the product-level demand independence; details presented in Table~\ref{tab:app1}). Both of these advancements are useful in practice and not simple to achieve via existing techniques. Please refer to Section~\ref{app: multi-product} for details.

\noindent \underline{\textbf{Application II: multi-echelon serial inventory system.}} We then consider a multi-echelon  inventory system with serially arranged stages, where each stage holds its own inventory and during the selling, the inventory at an upper stage can be transported to a lower stage with additional fares. This may model the serial production and delivery processes from the factory to upstream and downstream warehouses, and finally to the store. The transportation and coordination among multiple stages raise technical challenges to the existing methods -- the inventory process analyzed in the popular SGD-based inventory algorithms (e.g., \cite{huh2009nonparametric, shi2016nonparametric}) grows from single-dimension to \emph{multi-dimension}, and it is not clear how to deal with this extra complexity by queuing-theory-based analytic methods. In contrast, our minibatch SGD algorithm circumvents this issue and greatly simplifies the analysis. Please refer to Section~\ref{app: multi-echelon} and Section~\ref{ecapp:multi-echelon} in the supplementary materials for details.

\noindent \underline{\textbf{Application III: one-warehouse and multi-store inventory system.}}  The simplicity of our minibatch-SGD-based meta-policy enables us to handle even more complex inventory systems. We demonstrate this by considering a one-warehouse and multi-store system where the decision at each time period is \emph{two-stage}: first, the order-up-to levels at the warehouse and all stores are decided, and then, after observing the demand, the delivery amounts from the warehouse to the stores are decided. The optimal delivery decisions can be formulated by a \emph{linear program}, and its inventory process is clearly quite difficult to analyze. Thanks to our minibatch-SGD-based meta-policy, to achieve optimal regret, we ``only'' need to construct a gradient estimator for the solution of the two-stage linear program and prove the smoothness of the problem. Indeed, achieving these two goals already calls for much technical effort --- we build the gradient estimator from the dual program and establish the smoothness property via a novel sample-based analysis technique. We believe our techniques can be useful for learning other inventory control problems involving two-stage optimal decisions. Especially, as far as we know, we are the first to propose such a sample-based analysis technique for smoothness property, and this technique itself has an independent interest. Please refer to Section~\ref{app: owms} for details.

\subsection{Related Works}\label{sec:relatedwork}
In this section, we introduce the following streams of literature related to the three applications of our paper and other related topics such as nonparametric inventory control in other settings and minibatch SGD methods.
We discuss how our paper is appropriately placed into contemporary literature by giving comparisons with closely-related existing works.

\noindent \textbf{Multi-product inventory system with constraint(s).}
For classic research assuming that the firm has full knowledge of the demand distribution on multi-product inventory systems with constraints, there are two streams of literature. Some works consider a single warehouse-capacity constraint \citep{ignall1969optimality,erlebacher2000optimal}, while others consider multiple resource constraints \citep{lau1995multi, decroix1998optimal, downs2001managing} (see \cite{turken2012multi} for a detailed review).

Our first application is closely related to \cite{huh2009nonparametric,shi2016nonparametric}. We have discussed the technical differences between our paper and \cite{huh2009nonparametric,shi2016nonparametric} in the introduction section above. In Section~\ref{app: multi-product}, we will present a more detailed comparison of the settings between our paper and \cite{shi2016nonparametric}.

\noindent  \textbf{Multi-echelon serial inventory system.}
Multi-Echelon serial inventory systems involve multiple serially arranged stages of inventory storage and distribution. 
 As summarized in the book \citet{snyder2019fundamentals}, there are two common interpretations of the stages in a multi-echelon inventory system. First, stages could represent locations in a serial supply chain network, and in this case, links among the stages represent physical shipments of products. Second, stages could represent processes that the finished product must undergo during manufacturing, assembly, and distribution. In this case, links among the stages represent transitions between steps in the process.  
There have been many classic works on these systems \citep{chen1994lower, shang2003newsvendor}. For a detailed discussion of classic research on multi-echelon serial inventory systems, please refer to the book \citet{snyder2019fundamentals}.

As for the research on non-parametric methods for multi-echelon inventory problems,
\cite{zhang2021sampling} applied sample average approximation (SAA) policy to multi-echelon inventory systems and obtained an upper bound of sample size that guarantees the theoretical performance of their policy. 
{
\cite{yang2022non} designed an algorithm based on the online gradient algorithm proposed by \cite{huh2014online} for a two-echelon inventory system with instantaneous replenishment. In Section \ref{sec:Discussions on Multi-echelon inventory Systems} in the supplementary materials, we carefully compared  \cite{yang2022non} with our Application II.} \cite{zhang2022no} combined SAA and SGD to develop learning algorithms for two-echelon serial systems under both centralized and decentralized cases.

\noindent  \textbf{One-warehouse and multi-store inventory system.}
One-Warehouse multi-store inventory systems, first studied by \cite{clark1960optimal}, have also received a great deal of attention in the literature. In these systems, the central warehouse is replenished from external suppliers and stores receive periodic replenishment from the central warehouse. Some works consider the case where the central warehouse can only be replenished at the beginning of the selling season \citep{chao2021adaptive, miao2022asymptotically}. Some consider the problem of dynamic lot-sizing of the warehouse and distribution of inventory to stores
\citep{roundy198598, chen2001coordination}.

\cite{bekci2021inventory} considered the demand learning of a one-warehouse and multi-store system, in which the central warehouse receives an initial replenishment and distributes its inventory to multiple stores. The goal of their algorithm is to learn how to distribute inventory to the stores. We consider a slightly different system, in which the inventory of the warehouse is replenished periodically, and the goal of the learning algorithm is to decide the optimal order-up-to level of all installations. 

\noindent  \textbf{Other related works on dynamic inventory control.} 
    In recent years, there has been a growing body of research on inventory system management in cases where the demand distribution is unknown to the decision-maker. The decision-maker should learn the demand distribution on the fly, 
    and the demand distribution is either assumed to be in a specific parametric form (see the review in \cite{huh2009nonparametric} for a detailed discussion of the parametric demand setting) or nonparametric. Some works assume that the demand data can be fully observed, in which case the SAA method is commonly used to learn the distribution based on the observed samples of the demand \citep{levi2015data,cheung2019sampling,keskin2021nonstationary,lin2022data},
    while others consider the case where only sales data (i.e., censored demand) can be observed, which poses additional challenges in balancing the exploration and exploitation trade-off. For example, 
    \cite{huh2011adaptive} proposed a data-driven algorithm for a repeated newsvendor problem based on the Kaplan-Meier estimator; 
    \cite{besbes2013implications} characterized the implications of demand censoring on performance; \cite{yuan2021marrying} considered inventory systems with fixed ordering costs; \cite{zhang2018perishable} considered periodic-review perishable inventory systems with fixed lifetime;
    \cite{gao2022efficient} proposed an efficient learning framework for multi-product inventory systems with customer choices;
    \cite{chen2022learning} studied the inventory systems with  uncertain supplies. 
    Interested readers can refer to the review \citep{Gao2022} for a detailed discussion of inventory control with censored demand.

{\noindent  \textbf{Related works on cycle update policies.}
Several cycle-update learning algorithms that operate in cycles and use the same target variable in each cycle have been proposed for different settings of inventory management problems in the literature (e.g., \cite{huh2009adaptive,zhang2018perishable,chen2020optimal,zhang2020closing}). 
{Their algorithm details and analysis are also quite different from ours mainly from the following three aspects. 
First, compared with our minibatch-SGD-based meta-policy, their methods are designed to resolve different technical issues.  
Second, the cycle/minibatch-length patterns are different. Third, the update of the target variables is based on different gradient estimators.
In Section~\ref{sec:discussion on difference} in the supplementary materials, we carefully discuss these differences.}
}

\noindent  \textbf{Related works on minibatch SGD.}  Minibatch SGD, a variant of SGD, updates parameters based on the average gradient of the loss function with respect to the parameters for each minibatch. 
In the literature, the study on minibatch SGD mainly focuses on the convergence rate analysis of the minibatch SGD \citep{bubeck2015convex,ghadimi2016mini,bottou2018optimization,jofre2019variance} or the advantages of minibatch SGD  handling distributed
learning and making more efficient use of the available computational resources \citep{dekel2012optimal,li2014efficient}. 
{In our algorithm, we adopt the minibatch SGD as the update method in our inventory control meta-policy. Unlike the previous works on minibatch SGD, in the context of inventory control problems,  we need to study the cumulative regret of minibatch SGD with the limited number of minibatches, where this limit is determined by the optimal regret of different inventory management settings and may be a function of $T$ (i.e., $\O(\sqrt{T})$ in the general convex setting and $\O(\log T)$ in the strongly convex setting, see Section~\ref{sec:regret of meta policy} for details). 
}

\subsection{Notations}
The vectors throughout this paper are all column vectors.
We use $[n]$ to denote the set $\{1,2,\dots,n\}$  for any $n\in\mathbb N^+$.
For vectors $\bx, \by \in \R^n$, we use $\bx\leq \by$ ($\bx \geq \by$ respectively) to denote $x_i\leq y_i$ ($x_i \geq y_i$ respectively) for all $i\in[n]$. We use  $\mathbb{I}[\bx\geq\by]$ to denote the indicator vector, where $(\mathbb{I}[\bx\geq\by])_i=0$ when $x_i<y_i$ and $(\mathbb{I}[\bx\geq\by])_i=1$ when $x_i\geq y_i$. We use $\bx\odot\by$ to denote the element-wise product of $\bx$ and $\by$.
For $\bx\in\R^n$, $\|\bx\|_2$ is defined as $(\sum_{i=1}^n(x_i)^2)^{1/2}$.
We use $\mathbf{I}_{n}$ to denote the identity matrix of order $n$, and we use
$\mathrm{Diag}\{a_1,a_2,\dots,a_n\}$ to refer to the diagonal matrix with the  diagonal vector $(a_1,a_2,\dots,a_n)^{\top}$.
We use $X\sim Y$ to denote two random variables $X$ and $Y$ have the same distribution.
We use $A_{i,j}$ to denote the $i$-th row and $j$-th column element of matrix $A$.
We use $\bm A_{j}$ to denote the $j$-th column vector of matrix $A$.
We use the big-$\O$ notation $f(T)=\O(g(T))$  to denote that $\limsup_{T \to \infty}f(T)/g(T)< \infty$.
\subsection{Organization}
  
The remainder of this paper is organized as follows. In Section~\ref{sec:a genaral framework}, we propose a general inventory system framework, which includes lots of inventory systems with myopic optimal policy. In Section~\ref{sec: meta policy}, we present our minibatch-SGD-based meta-policy and discuss the high-level ideas of the policy design. 
In Section~\ref{sec:regret of meta policy}, we present the main theorems that upper bound the regret of our meta-policy in convex and strongly convex cases.
In Sections~\ref{app: multi-product},~\ref{app: multi-echelon}, and~\ref{app: owms}, we apply our minibatch-SGD-based meta-policy to three inventory control problems.
To demonstrate the empirical performance of our policy, we conduct several numerical experiments of the first application and present the results in Section~\ref{sec:numerical experiments}. In the end, we give a summary of our paper in Section~\ref{sec:conclusion}. The proofs of most technical lemmas and the additional experimental results are included in the supplementary materials.

\section{A General  Framework for Inventory Systems}\label{sec:a genaral framework}
We consider the following general framework that covers a wide range of inventory management problems with myopic optimal policy.

In each period of a $T$-period selling season, the demand $\bd_t$ is sampled from $\bD_t$, which is independent and stationary over time with common distribution $\bD.$\footnote{{We assume $\P[\bD>0]>0$.}} A firm needs to decide its order-up-to level $\by_t$ to minimize its cost. We assume the firm does not know the distribution of the demand a priori but adjusts its decision adaptively.

Specifically, in each period $t \in [T]$, we assume the following sequence of events.
\begin{enumerate}[nolistsep]
    \item At the beginning of period $t$, the firm observes the on-hand inventory $\bx_t$, where $\bx_t$ could be the inventory of products at different installations. The initial inventory is given as $\bx_1$. 
    \item The firm decides its order-up-to level $\by_t \geq \bx_t$ under the constraint $\by_t \in \Gamma$, where $\Gamma$ is a compact convex constraint set including practical constraints such as the space constraint of the warehouses (please refer to Section~\ref{app: multi-product} for the detailed discussion of the constraint set).  We assume the replenishment of inventory is instantaneous.\footnote{In this paper, we adopt the commonly used setting in the literature that the lead time of inventory systems is zero \citep{huh2009nonparametric,shi2016nonparametric,miao2022asymptotically,yang2022non,zhang2022no}. Studies such as \cite{miao2022asymptotically,lei2022joint}, which focus on multi-echelon inventory systems, point out that the zero lead time assumption is reasonable for many applications of these systems, particularly when the stations are located within a few hours or overnight driving distance.}
    
    \item The demand $\bD_t$ is realized as $\bd_t$. {The inventory transition function is denoted as $\mathfrak{f}(\by_t,\bd_t)$}. At the end of the period, the leftover inventory (on-hand inventory at the beginning of the next period) becomes {$\bx_{t+1}=\mathfrak{f}(\by_t,\bd_t)\leq \by_t$} according to the dynamics of specific inventory systems.  
\end{enumerate}
The expected cost of each period $Q(\by)$, whose formulation relies on specific applications, is decided by the order-up-to inventory level $\by$.
The goal of the firm is to design a policy $(\pi_t)_{t \geq 1}$ 
with $\pi_t:\mh_t\mapsto \by_t\in \{\by_t\in\Gamma\mid\by_t\geq\bx_t\}$ (where 
$\mh_t$ is the historical data at the beginning of period $t$)  and minimize the expected total inventory cost $\E\left[\sum_{t=1}^TQ(\by_t)\right]$.

Now we turn to the discussion of the performance measure of a given policy. It is easy to note that when the distribution is known a priori, the optimal policy is myopic,\footnote{{Note that the total inventory cost is $\sum_{t=1}^T\E\left[Q(\by_t)\right] {\geq} T\min_{\by\in\Gamma}Q(\by) := T \cdot Q(\by^*)$. Besides, due to $\bx_{t+1}=\mathfrak{f}(\by_t,\bd_t)\leq \by_t$, we can apply the myopic policy $\pi^* := (\by^*,\dots,\by^*)$ to get the optimal inventory cost $T\cdot Q(\by^*)$.}} which set the order-up-to level as a solution to the following problem:
\begin{equation*}
    \by^* = \argmin_{\by\in\Gamma} Q(\by),
\end{equation*}
in each period of the $T$-period selling season.
Thus we could define the regret of a given policy as 
\begin{equation}\label{eq:regretdef}
    \mathcal R(T) = \E \left[\sum_{t=1}^T (Q(\by_t)-Q(\by^*))\right],
\end{equation}
where $(\by_t)_{t \geq 1}$ is the sequence of order-up-to level given by the policy $(\pi_t)_{t \geq 1}$.

\section{A Minibatch-SGD-Based Meta-Policy}\label{sec: meta policy}

In this section, we propose a minibatch-SGD-based meta-policy (Algorithm~\ref{alg1}) for the learning-while-doing inventory control problem under the general inventory system framework.

Our meta-policy is motivated by the following key observation: implementing the SGD method in our inventory system framework is difficult because it switches the target inventory level $\Theta(T)$ times over $T$ periods, causing frequent low-target-inventory-level issues. However, when the current inventory level does not exceed the target level, it will remain at or below the target level if we do not \emph{switch} the target level over the following periods. 

This motivates us to adopt the minibatch SGD as our update method, which has fewer switches than the SGD. Specifically, we maintain a constant target inventory level during an epoch of  consecutive periods and we denote the epoch by a \emph{minibatch}. If we only have a few number of minibatches, the low-target-inventory-level issues  happen less frequently. Nevertheless, when such an issue does arise between two minibatches, we require a \emph{transition solver} to control the regret incurred during the issue. 
We will also need to design a \emph{gradient estimator} for the inventory cost function.  These two subroutines may be problem specific, but can be directly plugged into our meta-policy as long as it meets the simple interface requirements, which are defined as follows.

\noindent \underline{\bf Application specific subroutines.}
 Since our meta-policy is a gradient-based algorithm, the gradient estimator of the expected gradient needs to be designed carefully for different applications.  Specifically, a {\sc Well-behaved Gradient Estimator} satisfying these properties is defined as follows.
  
 \begin{definition}[Well-behaved Gradient Estimator]\label{def:gradient estimator} We call a gradient estimator  $\hat \nabla Q(\by)$  a well behaved gradient estimator, if it satisfies:
\begin{enumerate}[nolistsep]
        \item The estimator $\hat\nabla Q(\by)$ is unbiased, i.e., $\E\left[\hat\nabla Q(\by)\right] =   \nabla Q(\by) $.
        \item There exists some constant $\sigma >0$ such that the variance of the estimator is bounded by $\sigma^2$, i.e.,  $\E\left[\|\hat\nabla Q(\by) -\nabla Q(\by)\|_2^2\right] \leq \sigma^2.$
\end{enumerate}
\end{definition}

In each period $t$, if the target inventory level $\bw_t\in\Gamma$ is larger than the on-hand inventory $\bx_t$ element-wisely, i.e., $w_{i,t} \geq x_{i,t}$ for all $i$, we call $\bw_t$ a  
\emph{feasible} target. 
Otherwise, we call $\bw_t$ a 
\emph{infeasible} target.
 In the case where the  target inventory level $\bw_t$ is infeasible, we will keep calling the procedure  {\sc Transition Solver}, which is defined as follows, until the on-hand inventory drop below the target inventory level.
\begin{definition}[Transition Solver]\label{def:transition solver}
  A transition solver $\zeta^{\mathrm T}(\bx_t,\bw)$ takes $\bx_t$ and $\bw$ as inputs and returns the order-up-to level $\by_t\in \{\by_t\in\Gamma\mid\by_t\geq\bx_t\}$.
  Given the initial on-hand inventory $\bx_0$ and a fixed infeasible target inventory level $\bw$ with respect to $\bx_0$,  produce a sequence of on-hand inventory level $(\bx_t)_{t\geq1}$ by iterations \[\by_t=\zeta^{\mathrm T}(\bx_t,\bw)\text{ \quad and \quad } \bx_{t+1}=\mathfrak{f}(\by_t,\bd_t).\] 
  The transition solver makes sure that there exists a universal positive constant $M$ only depending on the problem instance that for all possible $\bx_0$ and  $\bw\in\Gamma$, we have $$\E[s(\bx_0,\bw)]\leq M,$$ where $s(\bx_0,\bw)$ is the first time when $\bw$ become feasible, i.e., 
  $s(\bm{x}_{0},\bw):=\min\{t\mid\bw \geq \bx_{t}\}.$ 
\end{definition}
The design of the {\sc Transition Solver} is specific to different inventory systems. Moreover, it is worth noting that for each application, there may exist many transition solvers that satisfy the above definition (for example, in some cases, simply setting the order quantity to $0$ in the waiting period can form a trivial transition solver). Any of these solvers can be implemented in our meta-policy to obtain a theoretical guarantee for regret bound.
In this paper, we aim to design transition solvers that not only satisfy the requirements of the above definition, but are also practical and heuristic in reality. In particular, we should strive to prevent a situation where we stop replenishing all products (or installations) just because we want to stop the replenishment of one product.

\begin{algorithm}[t!]
\caption{The Minibatch-SGD-Based Meta-Policy}
\label{alg1}
\begin{algorithmic}[1]
\State \textbf{Initialization: }Input a minibatch size sequence $\{n_1, n_2,\dots,n_{\tau_{\mathrm{max}}}\}$ where $\tau_{\mathrm{max}} = \min\{k:\sum_{\tau = 1}^k n_{\tau} \geq T\}$ and the stepsize $\eta$. Set a working period counter $l=0$, epoch step $\tau=1$, estimated gradient sequence $\mathcal G = \emptyset$ for the updating of minibatch SGD. Let the initial target level be the initial inventory level, i.e., $\bw_1 = \bx_1$.
\For{$t = 1,2,\dots,T$}

\If{$\bx_t \leq \bw_t$} {\color{blue}\Comment{Step 1: Decide the actually implemented order-up-to level.}}
\State Set order-up-to level $\by_t=\bw_t$ and let $l = l+1$.
\State Use the observed data to compute the {\sc Well-behaved Gradient Estimator} $\hat\nabla Q(\by_t)$. {\color{blue}\Comment{Step 2: Check whether to compute the gradient estimator.}} 

\State Add $\hat\nabla Q(\by_t)$ to the estimated gradient sequence $\mathcal G$.
\EndIf
\State \textbf{else } Call the {\sc Transition Solver} $\zeta^{\mathrm T}(\bx_t,\bw_t)$ to get the order-up-to level $\by_t$.
    \If{ $l = n_\tau$} {\color{blue}\Comment{Step 3: Update the target inventory level in a low switching manner.}}
    \State Update $\bw_{t+1}$ by minibatch SGD: $\displaystyle{
    \bw_{t+1}=\Pi_{\Gamma}[\bw_t- \frac{\eta}{n_\tau} \sum_{\hat g\in\mathcal G}\hat g ]}$. 
    \State Set $l=0$, $\tau=\tau+1$ and $\mathcal G = \emptyset$. 
    \EndIf
    \State \textbf{else}  Let $\bw_{t+1}=\bw_t.$

\State The on-hand inventory level updated as $\bx_{t+1} = \mathfrak f(\bx_t,\bd_t)$, where $\bd_t$ is the realized demand.

\EndFor
\State \textbf{end for}
\end{algorithmic}
\end{algorithm}

\noindent \underline{\bf Key parameters and steps of the meta-policy.}
With the {\sc Well-behaved Gradient Estimator} $\hat \nabla Q(\by)$  and the {\sc Transition Solver $\zeta^{\mathrm T}(\bx,\bw)$} in hand, the details of our meta-policy is specified in Algorithm~\ref{alg1}. As the input, our meta-policy needs the minibatch size sequence $\{n_1, n_2, \dots\}$ and the stepsize (a.k.a., learning rate) $\eta$. These are the key design parameters and will be carefully decided in various settings.

The main algorithm could be captured in the following three steps.

\noindent \underline{\it Step 1: Decide the actually implemented order-up-to level.}
In each period $t$, if the target inventory level $\bw_t$  updated by minibatch SGD is feasible, we simply use $\bw_t$ as our actually implemented inventory level $\by_t$.  If target inventory level $\bw_t$ is infeasible w.r.t.~the on-hand inventory level $\bx_t$, we call the {\sc Transition Solver} $\zeta^{\mathrm T}(\bx_t,\bw_t)$ to get the actually implemented inventory level $\by_t$.

\noindent \underline{\it Step 2: Check whether to compute the gradient estimator.}
We call the period $t$ a \emph{working period} if the target inventory level is feasible, and a \emph{waiting period} if the target inventory level is infeasible.
 In each working period $t$, we use observed data to compute the gradient estimator at $\bw_t$ and add it to the estimated gradient sequence $\mathcal{G}$. In waiting periods, we do not compute gradient estimators.  

\noindent \underline{\it Step 3: Update the target inventory level in a low switching manner.}
In addition, we maintain a \emph{working period counter} $l$, a batch number $\tau$, and a batch size sequence $\{ n_1,n_2,\dots,n_{\tau_{\mathrm{max}}} \}$.
 At the end of each period $t$ in batch $\tau$, we check whether gradient samples are sufficient for updating, that is $l = n_\tau$. If $l = n_{\tau}$, we update target inventory level $\bw_{t+1}$ by minibatch SGD, and reset the counter $l = 0$ and gradient information $\mathcal{G}=\emptyset$.

\section{Regret Analysis of the Meta-Policy}\label{sec:regret of meta policy}
{The regret of the meta-policy could be decomposed into two parts. 
The first part is the regret of waiting periods, and its magnitude is, at most, of the same order as the number of waiting periods. As elaborated in Section~\ref{sec:proof-of-thm-inv-general}, the number of waiting periods is of the order of the number of minibatches due to the definition of the {\sc Transition Solver}.
Consequently, to meet the lower bounds for the considered inventory management problems (see Proposition~\ref{prop:cvxlowerbound} and Proposition~\ref{prop:strcvxlowerbound}  in the supplementary materials for formal descriptions of the lower bounds), the order of minibatch number needs to be bounded (limited) by the order of the lower bound, i.e., $\O(\sqrt{T})$ in the general convex scenario and $\O(\log T)$ in the strongly convex scenario.

The second part is the regret of the working periods, which is highly related to the cumulative regret analysis of minibatch SGD with the limited number of minibatches. 
To showcase the optimality of our meta-policy,  when the expected cost function $Q(\by)$  is smooth\footnote{We say $Q(\by)$ is smooth ($\beta$-smooth), if there exists $\beta>0$ such that $\|\nabla Q(\by_1)-\nabla Q(\by_2)\|_2\leq \beta\|\by_1-\by_2\|_2$ for any $\by_1,\by_2\in\Gamma$.} and convex, we establish a cumulative regret of $\O(\sqrt{T})$ for the minibatch SGD algorithm employing $\O(\sqrt{T})$ minibatches; when $Q(\by)$ is smooth and strongly-convex\footnote{We say $Q(\by)$ is strongly convex ($\alpha$-strongly convex), if there exists $\alpha>0$ such that $ Q(\by_1)\geq Q(\by_2) + \nabla Q(\by_1)^\top (\by_1-\by_2)+ (\alpha/2) \|\by_1-\by_2\|_2$ for any $\by_1,\by_2\in\Gamma$.}, we establish a cumulative regret of $\O(\log T)$  for the minibatch SGD algorithm using $\O(\log T)$ minibatches in  Section~\ref{sec: minibatch sgd and its regret}. }

Then, in Sections~\ref{sec:proof of metapolicy}~and~\ref{sec:proof-of-thm-inv-general}, we will combine our analysis for minibatch SGD algorithm with our meta-policy, and prove our main theorems for the dynamic inventory problems.

\subsection{Minibatch SGD and its Regret Analysis}\label{sec: minibatch sgd and its regret}

Consider a general stochastic optimization problem $$\min_{\tbw\in\Gamma}F(\tbw)=  \E[f(\tbw;\xi)],$$
where $\Gamma\subset\mathbb{R}^n$ is a nonempty bounded closed convex set, and $\xi$ is a random seed.
Minibatch SGD update the variable $\tbw_t$ through the following way: 
\begin{equation}
\label{eq:minibatch sgd}
    \tbw_{\tau+1}=\Pi_{\Gamma}[\tbw_{\tau}-\eta g(\tbw;\bm{\xi}_{\tau})],
\end{equation}
 where  $\bm\xi = (\xi_{i,\tau})_{i=1}^{n_\tau}$, $\xi_{i,\tau}\sim\xi$,
$
g\left(\tbw_{\tau}; \bm{\xi}_{\tau}\right)=
\frac{1}{n_\tau} \sum_{i=1}^{n_\tau} \nabla f\left(\tbw_\tau ; \xi_{i,\tau}\right)
$, and $\nabla f(\tbw;\xi)$ is a unbiased gradient estimator of $\nabla F(\tbw)$.

Compared to the SGD algorithm, the minibatch SGD algorithm uses a batch of samples to compute the unbiased stochastic gradient estimator in each update. While the traditional research of the minibatch SGD algorithm either focuses on its convergence rate analysis or the advantages of handling distributed optimization (please refer to Section~\ref{sec:relatedwork} for a detailed review), in this section, we focus on the cumulative regret analysis minibatch SGD with the limited number of minibatches. In other words, we would like to achieve the unconditional optimal regret (i.e., the regret without the minibatch number constraint) using {the limited number of minibatches, i.e., $\O(\sqrt{T})$ in the general convex scenario and $\O(\log T)$ in the strongly convex scenario.}
Toward this goal, we need to carefully design the sequence of minibatch sizes and the stepsize (learning rate) and conduct the regret analysis under both convex and strongly convex cases.

At a higher level, we first combine the standard analysis of stochastic smooth convex optimization with the variance reduction property of minibatch SGD (i.e., $\E[\|\nabla F(\tbw_{\tau})-g(\tbw_\tau;\bm{\xi}_\tau)\|_2^2]=\E[\|\nabla F(\tbw_{\tau})-\nabla f(\tbw_\tau ; \xi)\|_2^2] /n_\tau$)) to derive the descent lemmas (please refer to Eq.~\eqref{eq:descent lemma cvx} and \eqref{eq:le:str:regret} in the supplementary materials) in the convex and strongly convex cases. Then we design the batch size and the stepsize carefully according to the limited number of minibatches and descent lemmas in two cases. Finally, we derive the cumulative regret by putting things together in the convex case and by an inductive method in the strongly convex case. We defer the proofs to Section~\ref{sec:proofs of metapolicy} in the supplementary materials.

For convenience, we define a $T$-period sequence $\bw_1, \bw_2,\dots,\bw_T $, where $ \bw_t=\tbw_1$ for $1\leq t< n_1$ and $\bw_t=\tbw_{\tau}$ for $\sum_{i=1}^{{\tau}-1}n_{i}<t\leq \sum_{i=1}^{{\tau}}n_{i}$. The expected cumulative regret is defined as 
\begin{equation*}
    \E\left[\sum_{t=1}^T (F(\bw_{t})-F(\bw^*))\right], 
\end{equation*}
where $\bw^* = \argmin_{\bw\in\Gamma} F(\bw).$
Next, we discuss the regret analysis of minibatch SGD for the convex case and strongly convex case individually.

{In the following lemma,  we propose two batch schemes, fixed-time batch scheme and any-time batch scheme, that establish the optimal cumulative regret of minibatch SGD with $\O(\sqrt{T})$ minibatches for stochastic convex and smooth constrained optimization problems.  The design of the fixed-time batch scheme requires the knowledge of the time horizon $T$ in advance, while the any-time batch scheme enjoys the any-time property -- the horizon $T$, does not need to be fixed upfront, and thus is more flexible and adaptive.}
\begin{lemma}\label{le:cvx}
        Suppose that $F(\bw)=\E[f(\bw;\xi)]$ is  convex and $\beta$-smooth 
        satisfying $\max_{\bw\in\Gamma}\|\nabla F(\bw)\|_2\leq G$ and the bounded variance property $\E[\|\nabla f(\bw;\xi)-\nabla F(\bw)\|_2^2]\leq \sigma^2$, and  $\Gamma$ be a bounded convex set satisfying $\max_{\bw_1,\bw_2\in\Gamma}\|\bw_1-\bw_2\|_2\leq R$. Run the minibatch stochastic gradient algorithm \eqref{eq:minibatch sgd} with stepsize $\eta<1/\beta$ and batch size $n_{\tau}$ specified in the following two cases for all $1\leq \tau\leq \tau_{\mathrm{max}}$, where   $\tau_{\mathrm{max}} = \min\{k:\sum_{\tau= 1}^k n_{\tau} \geq T\}$.

\begin{enumerate}[nolistsep]
    \item {\bf Fixed-time batch scheme.} When $n_{\tau}= \lceil\sqrt{T}\rceil$ for all $1\leq \tau\leq \tau_{\mathrm{max}}$, we have 
    \begin{equation*}
 \E\left[\sum_{t=1}^T (F(\bw_{t})-F(\bw^*))\right]\leq (\sqrt{T}+1)\left(\frac{R^2}{2\eta}+\frac{\eta\sigma^2}{2(1-\eta\beta)}+  GR\right).
    \end{equation*}
\item { {\bf Any-time batch scheme.} When $n_{\tau}= K\tau$ where $K$ is a positive integer, we have 
    \begin{equation*}
 \E\left[\sum_{t=1}^T (F(\bw_{t})-F(\bw^*))\right]\leq \left(1+\sqrt{\frac{2T}{K}+1}\right) \left(\frac{KR^2}{\eta}+ \frac{\eta\sigma^2}{(1-\eta\beta)}\right)+KGR.
    \end{equation*}}
\end{enumerate}

\end{lemma}
\begin{remark}\label{remark:cvx lower bound minibatch}
  It is well known that the optimal convergence rate is $\O(1/\sqrt{T})$ for convex smooth stochastic optimization \citep{nemirovskij1983problem}, which implies the $\O(\sqrt{T})$ regret established by the above lemma is optimal. 
\end{remark}

The following lemma establishes the cumulative regret of minibatch SGD with $\O(\log T)$ minibatches for stochastic strongly convex and smooth constrained optimization problems.
\begin{lemma}\label{le:str}
   Suppose that $F(\bw)=\E[f(\bw;\xi)]$ is $\alpha$-strongly convex and $\beta$-smooth 
   satisfying $\max_{\bw\in\Gamma}\|\nabla F(\bw)\|_2\leq G$ and the bounded variance property $\E[\|\nabla f(\bw;\xi)-\nabla F(\bw)\|_2^2]\leq \sigma^2$, and  $\Gamma$ be a bounded convex set satisfying $\max_{\bw_1,\bw_2\in\Gamma}\|\bw_1-\bw_2\|_2\leq R$. Running the minibatch stochastic gradient algorithm \eqref{eq:minibatch sgd} with stepsize $\eta\leq \min\{\alpha/2,1/\alpha,1/2\beta\}$ and batch size  $n_{\tau}=\lceil\varsigma^{\tau-1}\rceil$ for $1\leq \tau\leq \tau_{\mathrm{max}}$, where   $\tau_{\mathrm{max}} = \min\{k:\sum_{\tau= 1}^k n_{\tau} \geq T\}$,  $\varsigma=1/{\gamma}$ and $\gamma=1-\eta\alpha+2\eta^2\in(0,1)$, we have,  
   \begin{align*}
    \E\left[\sum_{t=1}^T (F(\bw_{t})-F(\bw^*))\right]\leq
     GR +\left(\frac{\kappa\varsigma^2}{\eta}+\frac{\varsigma^2\eta^2 \sigma^2}{(1-\eta \beta)}\right)\left(\frac{\ln{((\varsigma-1)T+1)}}{\ln{\varsigma}}+2\right),
\end{align*}
where $\kappa = \max\{R^2,\sigma^2\}$.
\end{lemma}
\begin{remark}\label{remark:str lower bound minibatch}
{\cite{hazan2014beyond} established the $\Omega(\ln T)$ regret lower bound for stochastic strongly-convex optimization (see Theorem 18 therein),  which shows the $\O(\log{T})$ regret established by the above lemma is optimal.} From  Lemma~\ref{le:str}, it is also easy to note that the minibatch SGD algorithm achieves the optimal regret with the exponentially increasing batch size $n_{\tau}=\lceil\varsigma^{\tau-1}\rceil$ for $1\leq \tau\leq \tau_{\mathrm{max}}$, using only $\O(\log{T})$ minibatches.
\end{remark}

{  \cite{bubeck2015convex} (Section 6.2 therein) conducted a convergence rate analysis of the minibatch SGD for smooth and convex functions. In the following, we discuss its connections to and differences from our lemmas:}
\begin{enumerate}[nolistsep]
    \item {In terms of results, \cite{bubeck2015convex} focused on the batch scheme where the batch sizes are the same, and analyzed the convergence rate given the batch size as a parameter. In our work, the number of minibatches is subject to a constraint due to the inventory problems, and we provide various batch schemes (of both constant batch sizes and non-constant batch sizes that enjoy the any-time property) to achieve the optimal regret. }

    \item {In terms of the proofs, both \cite{bubeck2015convex} and ours adopt the standard techniques (using variance reduction to derive the descent lemma). By setting the constant batch size to be $\lceil\sqrt{T}\rceil$ and a straightforward derivation, \cite{bubeck2015convex} may achieve a similar result as the fixed-time batch scheme (Item 1) in our Lemma~\ref{le:cvx} with a different stepsize. 
However, our any-time batch scheme in the convex scenario (Item 2 of Lemma~\ref{le:cvx}) cannot be directly deduced from \cite{bubeck2015convex}. 
Besides, our proof for the strongly convex scenario (Lemma~\ref{le:str}) diverges from \cite{bubeck2015convex}, primarily due to the use of a different type of descent lemma (Eq.~\eqref{eq:le:str:regret}).}
\end{enumerate}

\subsection{The Main Theorems} \label{sec:proof of metapolicy}
We present the main theorems for our meta-policy.  For notational simplicity, let $R$ be the diameter of $\Gamma$ (i.e., $R := \max_{\by,\by' \in \Gamma} \Vert \by-\by' \Vert_2$);
let $\bar{Q}$ be the maximal expected inventory cost (i.e.,
$\bar{Q} := \max_{y \in \Gamma} Q(\by)$); let $G$ be the maximal norm of the gradient (i.e., $G:=\max_{y \in \Gamma} \nabla Q(\by)$); let $\kappa := \max\{R^2,\sigma^2\}$.
Since $Q(\by)$ is smooth and $\Gamma$ is compact, it holds that all of $R$, $\bar{Q}$ and $G$ are finite.

Our first main theorem (proof deferred to Section~\ref{sec:proof-of-thm-inv-general}) shows that the regret of our meta-policy is $\O(\sqrt{T})$ when the expected inventory cost $Q(\by)$ is convex and smooth, matching the lower bound established in Proposition~\ref{prop:cvxlowerbound}.
\begin{theorem}
\label{thm:inv,general}
Suppose that the expected inventory cost $Q(\by)$ is convex and $\beta$-smooth,  and the meta-policy is executed with the {\sc Well-behaved Gradient estimator} and the {\sc Transition Solver} satisfying Definitions~\ref{def:gradient estimator} and~\ref{def:transition solver}. Run the meta-policy with stepsize $\eta<1/\beta$ and batch size $n_{\tau}$ specified in the following two cases for all $1\leq \tau\leq \tau_{\mathrm{max}}$, where $\tau_{\mathrm{max}} = \min\{k:\sum_{\tau= 1}^k n_{\tau} \geq T\}$,
\begin{enumerate}[nolistsep]
    \item {\bf Fixed-time batch scheme.} when $n_{\tau}= \lceil\sqrt{T}\rceil$ for all $1\leq \tau\leq \tau_{\mathrm{max}}$, we have that
\[
\mathcal R_{\mathrm{cvx}}(T)  \leq 
(\sqrt{T}+1)\left(\bar{Q}M+\frac{R^2}{2\eta}+\frac{\eta\sigma^2}{2(1-\eta\beta)}+  GR\right)
.
\]
 \item {{\bf Any-time batch scheme.} When $n_{\tau}= K\tau$ where $K$ is a positive integer, we have that}
$${\mathcal R_{\mathrm{cvx}}(T)  \leq \left(1+\sqrt{\frac{2T}{K}+1}\right) \left(\bar{Q}M+\frac{KR^2}{\eta}+ \frac{\eta\sigma^2}{(1-\eta\beta)}\right)+KGR
.}$$
\end{enumerate}
\end{theorem}
When the expected inventory cost $Q(\by)$ is strongly convex and smooth, our second main theorem (proof deferred to Section~\ref{sec:proof-of-thm-inv-str}) shows that the regret of our meta-policy can be improved to $O(\log T)$, matching the lower bound established in Proposition~\ref{prop:strcvxlowerbound} for this case.

\begin{theorem}
\label{thm:inv,str}
Suppose that the expected inventory cost $Q(\by)$ is  $\alpha$-strongly convex and $\beta$-smooth,  and the meta-policy is executed with the {\sc Well-behaved Gradient estimator} and the {\sc Transition Solver} satisfying Definitions~\ref{def:gradient estimator} and~\ref{def:transition solver}. Running the meta-policy with stepsize $\eta\leq \min\{\alpha/2,1/\alpha,1/2\beta\}$ and batch size  $n_{\tau}=\lceil\varsigma^{\tau-1}\rceil$ for $1\leq \tau\leq \tau_{\mathrm{max}}$, where  $\tau_{\mathrm{max}} = \min\{k:\sum_{\tau= 1}^k n_{\tau} \geq T\}$,  $\varsigma= 1/{\gamma}$ and $\gamma=1-\eta\alpha+2\eta^2\in(0,1)$, we have that
\[
\mathcal R_{\mathrm{str}}(T)  \leq 
 \bar{Q}M\frac{\ln T}{\ln \varsigma}+GR +\left(\frac{\kappa\varsigma^2}{\eta}+\frac{\varsigma^2\eta^2 \sigma^2}{(1-\eta \beta)}\right)\left(\frac{\ln{((\varsigma-1)T+1)}}{\ln{\varsigma}}+2\right)
.
\]
\end{theorem}

\section{Application I: Multi-Product and Multi-Constraint Inventory System}\label{app: multi-product}
In this section, we apply our meta-policy to a multi-product inventory system with multiple constraints and censored demand.
In each period, given the on-hand inventory and the historical sale quantities (censored demand), the firm needs to decide the order quantities of $n$ products under $m$ constraints to minimize its total cost. 

\subsection{Problem Formulation}
For each time $t\in[T]$ and product $i\in[n]$, the demand of $i$-th product at time $t$ is denoted by random variable $D_{i,t}\sim D_i$.
For convenience, we denote the random vector $(D_{1,t},D_{2,t},\dots,D_{n,t})^{\top}$ as $\bD_t$, its realization $(d_{1,t},d_{2,t},\dots,d_{n,t})^{\top}$ as $\bd_t$, and  $(D_{1},D_{2},\dots,D_{n})^{\top}$ as $\bD$.

\noindent \underline{\bf Sequence of events.}
In each period $t \in [T]$, we assume the following sequence of events:
\begin{enumerate}[nolistsep]
    \item At the beginning of period $t$, the firm observes the on-hand inventory $\bx_t = (x_{1,t},x_{2,t},\dots,x_{n,t})^{\top}$. The initial on-hand inventory is given as $\bx_1$.
    \item The firm decides its order-up-to inventory level  $\by_t \geq \bx_t$ and order quantity $\bq_t = \by_t-\bx_t$. 
    We assume instantaneous replenishment. Therefore after-delivery inventory is $\by_t$, which is restricted by the following inequalities
\[
\by_t \in \Gamma = \{\by \in \mathbb{R}^n_+ \mid A \by \leq \bm{\rho} \},
\]
where $A \in \mathbb{R}^{m\times n}_+, \bm{\rho} = (\rho_1,\rho_2,\dots,\rho_m)^{\top} \in \mathbb{R}^m_+$.\footnote{In our theoretical analysis, we only require the constraint set $\Gamma$ to be a compact convex set.} Moreover, we assume that the initial on-hand inventory $\bx_1 \in \Gamma$. 
    \item The demand $\bD_t\sim\bD$ is realized, denoted by $\bd_t$, which is satisfied to the maximum extent using after-delivery inventory.  Unsatisfied demand is lost, and the firm only observes the sales data (censored demand),  $\bs_t=\min (\bd_t,\by_t)$. The leftover inventory level
    $$\bx_{t+1}=\mathfrak f(\by_t,\bd_t):=\max\{\bm{0},\by_t-\bd_t\}.$$
    \end{enumerate}
   The inventory cost\footnote{It is shown by \cite{huh2009nonparametric} that the ordering cost can be absorbed into the lost-sales cost term. Therefore, we simply assume the ordering cost is zero.} of the period $t$ is  
\[
 {\bh}^\top (\by_t - \bd_t)^+ +  {\bb} ^\top (\bd_t-\by_t)^+,
\]
where $ {\bh}=( {h}_1,h_2,\dots, {h}_n)^{\top}$ and $ {\bb}=( {b}_1,b_2,\dots, {b}_n)^{\top}$ are the per unit holding and lost-sales penalty cost.
The expected inventory cost $Q(\by)$ is  
\[
Q(\by) = \bh^\top \E[(\by - \bD)^+] + \bb ^\top \E[(\bD-\by)^+].
\]
The goal of the firm is to design a policy
to minimize the following expected total inventory cost:
\begin{equation}\nonumber
\E\left[\sum_{t=1}^TQ(\by_t)\right] = \sum_{t=1}^T[\bh^\top \E[(\by_t - \bD)^+] + \bb ^\top \E[(\bD-\by_t)^+] ].
\end{equation}

\noindent \underline{\bf Regularity assumptions about demand.}
\begin{assumption}
\label{ass1:generalass}
 In this application, we make the following 
basic assumptions on demand.
\begin{enumerate}[nolistsep]
    \item For each product $i\in [n]$, $D_{i,t}$ is i.i.d.~across time period $t \in [T]$ and we assume $D_{i,t}\sim D_i$. 
    \item For each product $i \in [n]$, $D_i$ has PDF $f_{D_i}(\cdot)$ and $f_{D_i}(x) \leq \beta_0$ for some $\beta_0 >0$.  
\end{enumerate}
\end{assumption}
\begin{remark}
    The above assumption is quite benign and similar assumptions are used in other inventory-related learning literature (see, e.g., \cite{huh2009nonparametric,shi2016nonparametric,zhang2020closing,yuan2021marrying}).
Assumption 1.2 that density functions $f_{D_i}(x), i \in [n]$ of all demand distributions are upper bounded by $\beta_0$ is a mild technique condition ensuring that cost function $Q(\by)$ is smooth. 
\end{remark}

The following assumption requires that the density $f_{D_i}(x), i \in [n]$ are lower bounded by some positive constant, which ensures strong convexity of cost function $Q(\by)$. As we have shown in Section~\ref{sec:regret of meta policy}, with the strong convexity of cost function $Q(\by)$, the performance of our meta-policy can be improved.
\begin{assumption}
\label{ass2:strongly}
  For each product $i \in [n]$, the PDF of $D_i$  is lower bounded by a positive constant. That is $f_{D_i}(x) \geq \alpha_0$ for some $\alpha_0 >0$.
\end{assumption}
\begin{remark}
    Similar assumptions appear in many inventory-related learning studies to ensure strong convexity and improve the performance of algorithms (see, e.g., \cite{huh2009nonparametric,besbes2013implications,shi2016nonparametric,zhang2018perishable}).
\end{remark}

The above-described multi-product inventory system setting is similar to (and inspired by) the setting considered in \cite{shi2016nonparametric}.  However, our setting is more general and practical as explained in the following two aspects.

\noindent \underline{\bf Multiple constraints.} In contrast to the single warehouse-capacity constraint considered in \cite{shi2016nonparametric}, we consider a more general multi-constraint setting. Multi-product inventory systems with multiple constraints have been investigated by many researchers in classic research. However, this paper is the first work considering the multiple constraints in the demand learning setting.
There are lots of interpretations of the constraints set $\Gamma = \{\by_t \in \mathbb{R}^n_+ \mid A \by_t \leq \bm{\rho} \}$ in classic works 
\citep{lau1995multi,decroix1998optimal,downs2001managing}.
Generally, $A_{ij}$ represents resource requirement of resource $i$ for keeping one unit product $j$, and each constraint $(A \by_t)_i\leq \rho_i$ in $\Gamma$ represents the resource constraint of resource $i$.

\noindent \underline{\bf Removing the product-level demand independence assumption.}
Another major difference is that we do not require that the demands of different products are independent, which is assumed in \cite{shi2016nonparametric}. 
Demands of products could be highly correlated in lots of cases. If two products have similar functions, their demands are usually negatively correlated, e.g., the products of the same type but different brands. If two products have complementary functions, their demands are usually positively correlated, e.g., the products that are often used together.  Therefore, our algorithm can be applied to a wider scope of multi-product inventory systems by allowing dependence among products.

To the best of our knowledge, it is not immediately clear how to adapt the method in \cite{shi2016nonparametric} to handle this more general setting. For example, the proof of the Lemma 5 in \cite{shi2016nonparametric} heavily relies on the single warehouse-capacity constraint structure and is hard to work for the multi-constraint version. Meanwhile, the independent assumption among products is essential for the analysis of Lemmas 6, 7 and 8 in \cite{shi2016nonparametric}. Without such assumption, the inter-arrival distribution $\tilde{D}_t$ defined therein may be zero, and the expected busy period cannot be bounded, hence the existing technique fails. {Please refer to Section \ref{subsec:shicong} for detailed discussions of the roles of these lemmas and the reasons why the method in \cite{shi2016nonparametric} fails in handling multiple constraints and demand correlations among products.}

\subsection{Design of the {\sc Well-behaved Gradient Estimator} and the {\sc Transition Solver}}\label{app1:design}
As mentioned before, we denote the expected one-period cost when the after-delivery inventory level is $\by$ as $$Q(\by) = \bh^\top \E[(\by - \bD)^+] + \bb ^\top \E[(\bD-\by)^+].$$
And the  gradient of $Q(\by)$ at $\by$ is given by $$\nabla Q(\by) = \E \left[ \bh\odot\mathbb{I}[\by>\bD]-\bb\odot\mathbb{I}[\by\leq\bD] \right].$$
Assuming $\bd$ is a realization of $\bD$, we could design the following gradient estimator of $Q(\by)$
    \begin{equation}\label{app1:gradient estimator}
    \hat\nabla Q(\by) = \bh\odot\mathbb{I}[\by>\bd]-\bb\odot\mathbb{I}[\by\leq\bd].
    \end{equation}
    Note that although we can only observe sales data, it is sufficient to compute gradient estimator $\hat\nabla Q(\by)$.
    The properties of $Q(\by)$ and $\hat\nabla Q(\by)$ are summarized in the following lemma.
\begin{lemma}\label{le:app1properties}
Properties of $Q(\by)$ and $\hat\nabla Q(\by)$:
\begin{enumerate}[nolistsep]
    \item Under 
Assumption \ref{ass1:generalass}, $Q(\by)$ is convex and $(\max_{i\in[n]}(h_i+b_i)\beta_0)$-smooth.
    \item Under Assumptions \ref{ass1:generalass} and \ref{ass2:strongly}, $Q(\by)$ is $(\min_{i\in[n]}(h_i+b_i)\alpha_0)$-strongly convex.
    \item The norm of gradient is bound as $\|\nabla Q(\by)\|_2\leq \sigma_0$, where $\sigma_0 := (\|\bh\|_2^2+\|\bb\|_2^2)^\frac12.$
    \item The estimator $\hat\nabla Q(\by)$ is unbiased, i.e.,$\E\left[\hat\nabla Q(\by)\right] =   \nabla Q(\by) $.
    \item It holds that $\E\left[\|\hat\nabla Q(\by) -\nabla Q(\by)\|_2^2\right] \leq \E\left[\|\hat\nabla Q(\by)\|_2^2\right] \leq \|\bh\|_2^2+\|\bb\|_2^2=\sigma_0^2.$
\end{enumerate}
\end{lemma}
The last two items of Lemma~\ref{le:app1properties} show that the gradient estimator $\hat\nabla Q(\by)$ (\ref{app1:gradient estimator}) satisfies Definition~\ref{def:gradient estimator}, hence is a {\sc Well-behaved Gradient Estimator}.
We now turn to the design of the {\sc Transition Solver}.\footnote{We would like to point out that simply setting the order quantity to $0$ is a trivial transition solver satisfying Definition~\ref{def:transition solver}. However, if we adopt this trivial solver, the inventory level of some products may drop to $0$ and thus become out of stock during the waiting periods, which is not practical. Our greedy projection transition solver could avoid this problem, and theoretically, it could handle a more general case that there are additional minimal inventory level constraints $\by_t\geq\underline{\bm{\rho}}$  in $\Gamma$.}

Our {\sc Transition Solver} $\zeta^{\mathrm T}(\bx_t,\bw)$ takes a on-hand inventory level $\bx_t$ and some infeasible target inventory level $\bw$ as inputs,  uses the following \emph{greedily projection} method to get an actually implemented inventory level $\by_t$: 
\begin{align}\label{app1:opt}
\by_t &= \argmin_{\by\in\R^n} \Vert \by- \bw \Vert_2^2
     \qquad \text{s.t.~}  A \by \leq \bm{\rho} \quad\text{and}\quad \by \geq \bx_t. 
\end{align}

The above optimization problem (\ref{app1:opt}) is a convex quadratic program, and hence can be solved efficiently. Moreover, since $\bx_t \geq 0$ is a feasible solution, the optimal solution to this optimization problem exists. The following lemma shows the solver designed above meets Definition~\ref{def:transition solver}.
\begin{lemma}\label{app1:le:transition solver}
  The solver $\zeta^{\mathrm T}$ defined in \eqref{app1:opt} satisfies Definition~\ref{def:transition solver}  with $M=n + 6n\beta_0 R$, where $R = \max_{\by,\by' \in \Gamma} \Vert \by-\by' \Vert_2$.
\end{lemma}

\subsection{Regret Analysis}\label{app1:regret}
Plugging the parameter constants in Lemmas~\ref{le:app1properties} and \ref{app1:le:transition solver} to the meta-results (Theorems~\ref{thm:inv,general} and \ref{thm:inv,str}), we obtain the following regret upper bounds for  multi-product multi-constraint inventory systems with censored demand in two cases. 

For convenience, we define notations 
$\beta = \max_{i\in[n]}(h_i+b_i)\beta_0$,  $\alpha = \min_{i\in[n]}(h_i+b_i)\alpha_0$, $\bar{Q} = \max_{\by \in \Gamma} Q(\by)$, {$R = \max_{\by,\by' \in \Gamma} \Vert \by-\by' \Vert_2\leq  2\sqrt{\sum_{j=1}^n(\min_{i\in[m]}\rho_j/A_{i,j})^2}$ (see Lemma~\ref{le:Upper bound of R} in the supplementary materials for this inequality),} $\sigma_0 = (\|\bh\|_2^2+\|\bb\|_2^2)^\frac12$  and $M=n + 6n\beta_0 R$ in the following two theorems.
\begin{theorem}\label{app1:thm:inv,general}
Under Assumption \ref{ass1:generalass}, invoking the meta-policy with  {\sc Well-behaved Gradient estimator} \eqref{app1:gradient estimator} and the {\sc Transition Solver} \eqref{app1:opt}, and
 running the meta-policy with stepsize $\eta<1/\beta$ and batch size $n_{\tau}= \lceil\sqrt{T}\rceil$ or $n_{\tau}= K\tau$, where $K$ is a positive integer, we have
$
 \mathcal R_{\mathrm{cvx}}(T) \leq \O(\sqrt{T})
$. Here the $\O(\cdot)$ notation hides the polynomial dependence on $\beta$,  $\bar{Q}$, $\eta$, $R$, $M$, $K$ and $\sigma_0$.
\end{theorem}

\begin{theorem}
\label{app1:thm:inv,str}
Under Assumptions \ref{ass1:generalass} and \ref{ass2:strongly}, invoking the meta-policy with  {\sc Well-behaved Gradient estimator} \eqref{app1:gradient estimator} and the {\sc Transition Solver} \eqref{app1:opt}, and
running the meta-policy with stepsize $\eta\leq \min\{\alpha/2,1/\alpha,1/2\beta\}$ and batch size  $n_{\tau}=\lceil\varsigma^{\tau-1}\rceil$, where  $\varsigma= 1/{\gamma}$ and $\gamma=1-\eta\alpha+2\eta^2\in(0,1)$, we have $\mathcal R_{\mathrm{str}}(T)  \leq \O(\log{T})$. Here the $\O(\cdot)$ notation hides the polynomial dependence on $\beta$, $\alpha$, $\varsigma$, $\bar{Q}$, $R$, $M$, and $\sigma_0$.
\end{theorem}

\section{Application II: Multi-Echelon Serial Inventory System} \label{app: multi-echelon}

We apply our meta-policy to a multi-echelon serial inventory system with censored demand. The problem in Section~\ref{app: multi-product} can be viewed as a single-echelon special case of the one discussed here.\footnote{For simplicity, we only discuss the single-product case. It is easy to extend our method to the multi-product case.}
\begin{figure}
    \centering
    \includegraphics[width =0.8\textwidth]{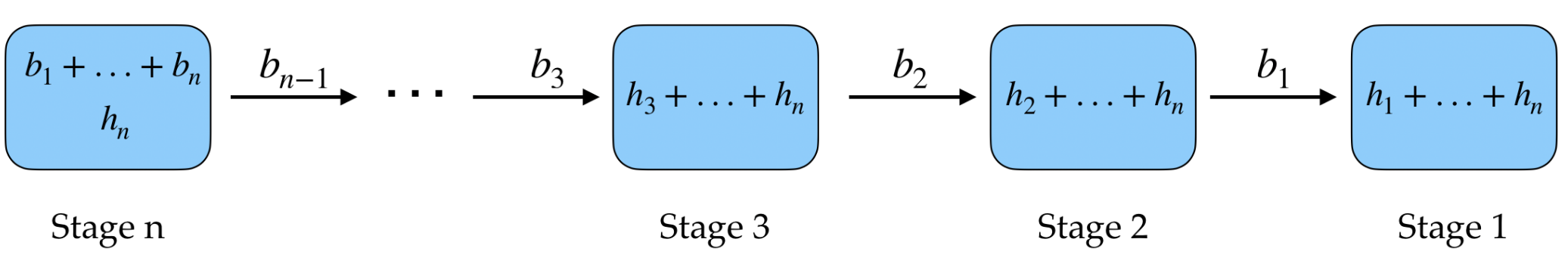}
    \caption{A multi-echelon serial 
inventory system with $n$ stages. Parameters $h_i$ and $b_i$ will be specified in Section~\ref{app2:problem formulation}.}
    \label{fig:multiechelon}
\end{figure}

Consider  a multi-echelon
inventory system with serially arranged stages in Figure~\ref{fig:multiechelon}. At the beginning of each period, given the on-hand inventory of all stages and the historical sale quantities (censored demand), the firm needs to decide the order-up-to levels at all stages. After that, 
Stage $1$ faces stochastic demand and could get replenishment from Stage $2$ if out of stock by paying extra transportation fares. If Stage $2$ is out of stock, it could get replenishment from Stage 3, and so on up the line to Stage $n$. In the same manner as the single echelon case, we assume the firm does not know the distribution of demand as a prior, but observes sales quantity and adjusts its policy in each period. We assume the replenishment is instantaneous and the ordering cost is zero.\footnote{It is reasonable to assume that the ordering cost of each stage plus the total transportation cost for downstream delivery remains constant. By a similar argument with \cite{huh2009nonparametric}, we could absorb the ordering cost into the lost-sales cost term. As a result, we only consider additional transportation fares during the selling and assume zero ordering costs.}

We give two motivating examples of the above multi-echelon inventory system setting as follows.
\begin{enumerate}[nolistsep]
    \item Consider a two-echelon serial inventory system managed by a firm, where the first stage is a store and the second stage is a warehouse. The store faces external demand directly, and if out of stock, the firm will transport the inventory from the warehouse to the store with extra transportation fares to satisfy the demand. Such a system is common and simple in practice.
    \item 
Blood banks, which are important in healthcare, can store both fresh blood and frozen blood. Fresh blood will be used first to satisfy the blood demand, and if the fresh blood is not enough, the bank will transform the frozen blood to fresh blood with the cost of thawing frozen blood and combining them with plasma (see \cite{jennings1973blood,nahmias1976myopic}).

\end{enumerate}

{
Due to space constraints, we present the detailed problem formulations and technical results of our Application II in Section~\ref{ecapp:multi-echelon} in the supplementary materials. We highlight the main ingredients needed in this application as follows.
\begin{enumerate}[nolistsep]
     \item \textit{Formulating the cost function.} 
     According to the above motivating examples, we formulate the inventory cost function carefully according to the inventory management process of the multi-echelon inventory systems. 
    \item \textit{Cost function simplification.}
    The natural formulation of the inventory cost function turns out to be quite complex to analyze. We simplify the original formulation to an equivalent form via case-by-case discussion, variable substitution and transformation.
    \item \textit{Invoking our meta-policy.}
    We establish the desired properties (convexity and smoothness) of the simplified cost function. Finally, we invoke our meta-policy and main theorems to obtain the regret bound for this application.
\end{enumerate}}

\section{Application III: One-Warehouse and Multi-Store Inventory System}\label{app: owms}
In this section, we use our meta-policy to solve a one-warehouse and multi-store inventory control problem, which serves as an illustrative example to show the efficacy of our meta-policy to handle inventory control problems involving two-stage decisions, which are extensively studied in the literature (e.g.  Section 1.1 of \cite{bertsimas2020predictive}, Section 1.3 of \cite{shapiro2021lectures}, Section 7.4 of \cite{snyder2019fundamentals}).  After designing the necessary subroutines (see Section \ref{app3:design}) and proving desired properties of the cost function (see Section~\ref{app3:convexity and smoothness}), we can apply our meta-policy to these inventory systems with two-stage decisions.   The complexity of inventory systems with two-stage decisions makes their learning tasks largely under-explored in literature, and also raises new challenges in our algorithm design and analysis. We believe our methods are of independent interest and may be adopted in other similar inventory control problems.

\begin{figure}
    \centering
    \includegraphics[width =0.4\textwidth]{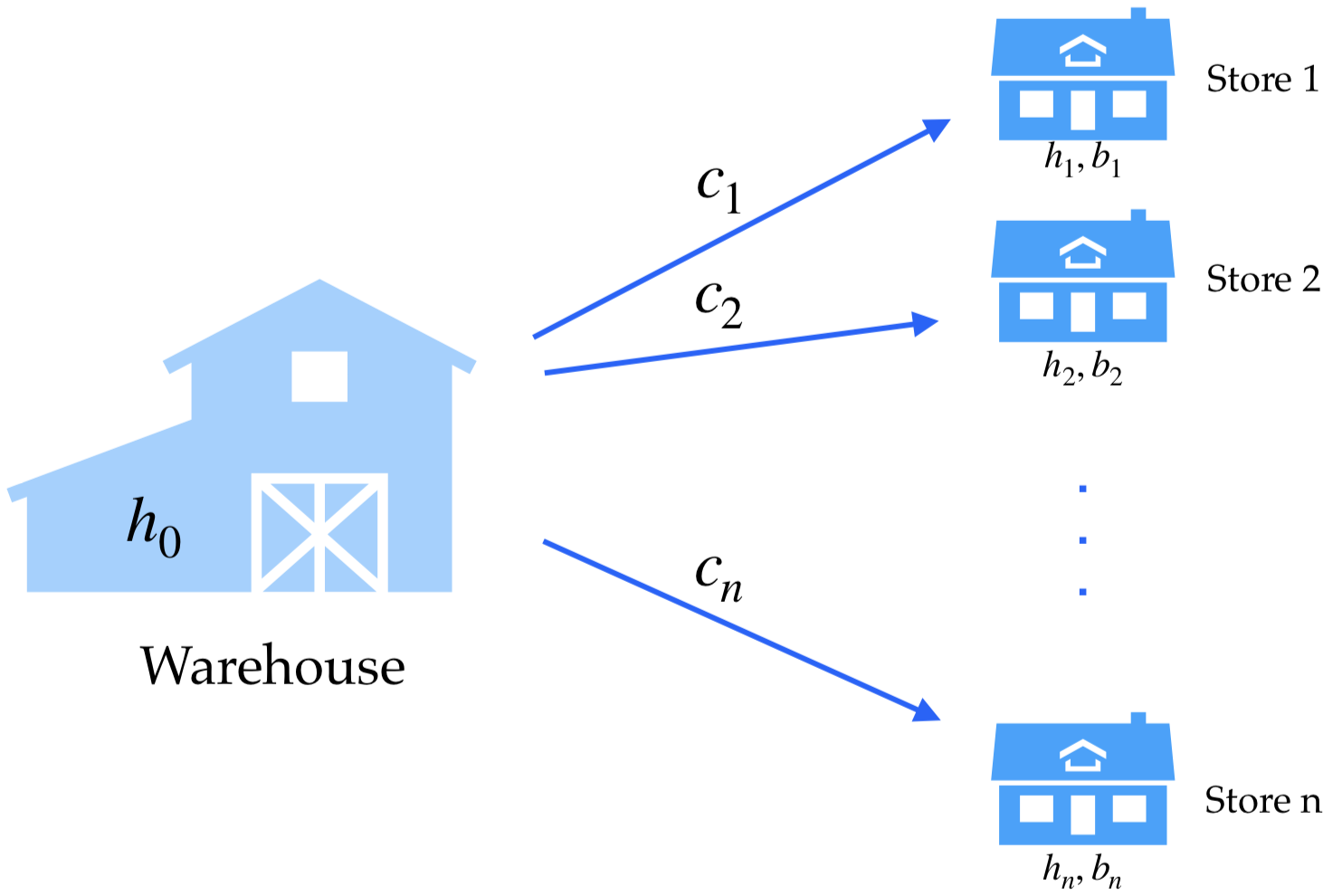}
    \caption{Inventory system with one warehouse \& $n$ stores. Parameters $h_i,$ $c_i$ and $b_i$ are specified in Section~\ref{app3:problem formulation}.}
    \label{fig:owms}
\end{figure}

\subsection{Problem Formulation}\label{app3:problem formulation}
We consider an inventory system with one warehouse and $n$ stores as Figure~\ref{fig:owms}. In the following, we also refer to the warehouse or a store as an \emph{installation}. The warehouse (labeled $0$), holding buffer inventory, can deliver products to other stores with some extra transportation fares. Each store, labeled as $1$ through $n$, independently faces external demands. The warehouse and all stores are managed by a firm to minimize its total cost. At the beginning of each period $t$, the firm decides how many products to order (first-stage decision). After the demand is realized, the firm decides how to deliver products from the warehouse to the stores (second-stage decision).

For each period $t\in[T]$ and store $i\in[n]$, the demand faced by store $i$ at time $t$ is denoted by the random variable $D_{i,t}\sim D_i$. For convenience, we use vector notations similar to the first application.

\noindent \underline{\bf Sequence of events.}
In each period $t \in [T]$, we assume the following sequence of events:
\begin{enumerate}[nolistsep]
    \item At the beginning of period $t$, the firm reviews the on-hand inventory $\bx_t = (x_{0,t},x_{1,t},\dots,x_{n,t})^{\top}$ at all installations. The initial on-hand inventory is given as $\bx_1$.
    \item The firm decides its order-up-to level $\by_t \geq \bx_t$ and order quantity $\bq_t = \by_t-\bx_t$ of the warehouse and all stores, We assume instantaneous replenishment and zero ordering cost (for the same reason as the second application). After the replenishment, inventory is $\by_t$, which is restricted as $\by_t \in \Gamma$,
   where $ \Gamma = \{\by\in \mathbb{R}^{n+1}_+ \mid \by\leq \bm{\rho} \}$ and  $\bm{\rho} = (\rho_0,\rho_1,\dots,\rho_{n})^{\top}\in\R_+^{n+1}$ is the capacity constraints of the installations.
    \item The demand $\bD_t\sim \bD$ is realized, denoted by $\bd_t$. Let  $z_{i,t} \geq 0$ be the delivery quantity from the warehouse to the $i$-th store. Then the firm decides the optimal delivery quantity $\bm{z}_t=(z_{1,t},z_{2,t},\dots,z_{n,t})^{\top}$, which is the optimal solution of the optimization problem $C(\by_t, \bd_t)$ where $C(\cdot, \cdot)$ is defined as below.
        \begin{equation}\label{eq:delivery}
\begin{aligned} 
C(\by,\bm{d}) := &\min_{\bm z} \bm{c}^\top \bm{z} + h_0(y_{0}-\sum_{i=1}^n z_{i})^+ + \sum_{i=1}^n \left[ h_i(z_{i}+y_{i}-d_{i})^+ + b_i(d_{i}-z_{i}-y_{i})^+\right] \\
&\begin{array}{r@{\quad}r@{}l@{\quad}l}
s.t. 
     &\sum_{i=1}^n z_{i} &\leq y_{0}, \\
     &z_i &\geq 0,\forall i \in [n] .\\
\end{array}
\end{aligned}
\end{equation}
Here, $\bm z=(z_1,z_2,\dots,z_n)^{\top}$, $\bh = (h_0,h_1,\dots,h_n)^{\top}$ is the unit holding cost at the warehouse and all stores, $\bb = (b_1,b_2,\dots,b_n)^{\top}$ is the unit lost-sales cost at all stores, and $\bm{c} = (c_1,c_2,\dots,c_n)^{\top}$ is the unit transportation fares from the warehouse to all stores.
     After the delivery, demand is satisfied to the maximum extent using after-delivery inventory $z_{i,t}+y_{i,t}$ at each store $i \in [n]$. Unsatisfied demand is lost.
 The leftover inventory (on-hand inventory of the next period $t+1$) is given by 
    $$
    \bx_{t+1} =\mathfrak f(\by_t,\bd_t) = \left( (y_{0,t}-\sum_{i=1}^n z_{i,t}),(y_{1,t}+z_{1,t}-d_{1,t})^+,\dots,(y_{n,t}+z_{n,t}-d_{n,t})^+\right)^{\top}.
    $$
\end{enumerate}
 We use $Q(\by)=\E[C(\by,\bm{D})]$ to denote the expected cost of the firm when the order-up-to level is $\by$. 
The goal of the firm is to design a policy
to minimize the following expected total inventory cost :
\begin{equation}\nonumber
\E\left[\sum_{t=1}^TQ(\by_t)\right] = \sum_{t=1}^T\E[C(\by_t,\bm{D}_t)].
\end{equation}
\noindent \underline{\bf Regularity assumptions about demand.}
Like Assumption~\ref{ass1:generalass} and Assumption~\ref{ass1:generalass:app2} made for the first two applications, we make the following  assumption for the one-warehouse and multi-store systems. 
\begin{assumption}
\label{ass1:generalass:app3}
In this application, we make the following basic assumptions on demand.
\begin{enumerate}[nolistsep]
    \item $\bm{D}_{t}$ is i.i.d.~across time period $t \in [T]$ with common distribution $\bD$. 
    \item $\bm D$ has PDF $f_{\bm D}(\cdot)$ and $f_{\bm D}(\bm x) \leq \beta_0$ for some $\beta_0 >0.$
\end{enumerate}
\end{assumption}

\subsection{Design of the {\sc Well-behaved Gradient Estimator} and the {\sc Transition Solver}}\label{app3:design}

While the gradient estimator can be obtained directly in the previous two applications, in this application the design of the gradient estimator is trickier due to the lack of a closed-form expression for $Q(\by)$. We will achieve this goal by leveraging the dual formulation of Problem (\ref{eq:delivery}).

\noindent \underline{\bf Design of the gradient estimator.} 
By introducing intermediate variables $z_0^{(1)}$, $z_i^{(1)}$ and $z_i^{(2)}$ for all $i \in [n]$, we could transform Problem (\ref{eq:delivery}) to the following linear program.
\begin{equation}\nonumber
\begin{aligned} 
&\min \bm{c}^\top \bm{z} + \sum_{i=0}^n h_i z_i^{(1)} + \sum_{i=1}^n b_iz_i^{(2)}\\
&\begin{array}{r@{\quad}r@{}l@{\quad}l}
s.t. 
     &\sum_{i=1}^n z_i &\leq y_0,\\
     &y_0-\sum_{i=1}^n z_i &\leq z_0^{(1)}, \\
    &z_i + y_i-d_i &\leq z_i^{(1)},i \in [n], \\
     &d_i - z_i-y_i &\leq z_i^{(2)},i \in [n], \\
     &z_i,z_0^{(1)},z_i^{(1)},z_i^{(2)}  &\geq 0,i \in [n].\\
\end{array}
\end{aligned}
\end{equation}
We further transform the inequality constraints into equality constraints by introducing the intermediate slack variables (in the forms of $z_0^{(3)}$, $z_0^{(4)}$, $z_i^{(4)}$, $z_i^{(5)}$). We then collect all decision variables ($z_i$ and intermediate variables) to form the vector $z'$. By correspondingly defining vectors $\bm{c}'$, $\bd'$ and matrices $H$, $W$, we may transform the above linear program to the following standard form. (This step is standard and the details are deferred to Section~\ref{subsec: proofoflemmaapp3} in the supplementary materials.) 
\begin{equation}\label{eq:app3standard}
\begin{aligned} 
C(\by,\bm{d}) := &\min (\bm{c}')^\top \bm{z}'\\
&\begin{array}{r@{\quad}r@{}l@{\quad}l}
s.t. 
     & H\by  + W\bm{z}' &= \bd' \\
     &\bm{z}' &\geq \bm{0}.\\
\end{array}
\end{aligned}
\end{equation}

Given a sample $\bd$, we need to compute $\nabla_{\by} C(\by,\bd)$ to get a gradient estimator of $Q(\by) = \E[C(\by,\bD)]$. We first derive the dual problem of \eqref{eq:app3standard} as
\begin{equation}\label{problem:dual}
    \max_{\bm{\pi}} \bm{\pi}^\top (\bd'- H\by ) 
    \quad\text{s.t. }  W^\top \bm{\pi} \leq \bm{c}',
\end{equation} 
which have the same value as the primal problem \eqref{eq:app3standard} due to the strong duality. 
We now invoke Danskin's theorem which states that, if $C(\by,\bd)$ is defined by a maximization problem, its gradient can be computed under any optimal solution $\bm{\pi}^*(\by,\bd) = \argmax_{\bm{\pi} \in \{W^\top \bm{\pi} \leq \bm{c}'\}} \bm{\pi}^{\top} (\bd'- H\by ).$  As a consequence, we can choose the gradient of $C(\by,\bd)$ as our gradient estimator,
\begin{equation}\label{estimator:app3}
   \hat{\nabla} Q(\by) := \nabla_{\bm y} C(\by,\bd) = -H^\top \bm{\pi}^*(\by,\bd).
\end{equation}
 The properties of $\hat{\nabla}Q(\by)$ are summarized in the following lemma and the proof is referred to Section~\ref{sec:app-app3-properties-estimator}.
\begin{lemma}
\label{le:app3-properties-estimator}
Under Assumption \ref{ass1:generalass:app3}, there exists $\sigma_0 := [(h_0 + \sum_{i=1}^n b_i)^2+ \sum_{i=1}^n(h_i+b_i)^2]^\frac 12$ such that: 
\begin{enumerate}[nolistsep]
    \item The norm of gradient is bounded by $\sigma_0$, i.e., $\|\nabla Q(\by)\|_2\leq \sigma_0$.
    \item The estimator $\hat\nabla Q(\by)$ is unbiased, i.e., $\E\left[\hat\nabla Q(\by)\right] =   \nabla Q(\by) $.
    \item It holds that $\E\left[\|\hat\nabla Q(\by) -\nabla Q(\by)\|_2^2\right] \leq \E\left[\|\hat\nabla Q(\by)\|_2^2\right] \leq  \sigma_0^2.$
\end{enumerate}
\end{lemma}
Lemma~\ref{le:app3-properties-estimator} show that the gradient estimator $\hat\nabla Q(\by)$ is a {\sc Well-behaved Gradient Estimator} satisfying Definition~\ref{def:gradient estimator}.

\noindent \underline{\bf Design of the {\sc Transition Solver}.} 
In the following, we design a  {\sc Transition Solver}   which not only satisfies the requirements in Definition~\ref{def:transition solver} but also balances inventory cost intuitively. 

Our {\sc Transition Solver} $\zeta^{\mathrm T}(\bx_t,\bw)$ takes an on-hand inventory level $\bx_t$ and some infeasible target inventory level $\bw$ as inputs, and output $\by_t$ as the following two cases. 
\begin{equation}\label{app3transitionsolver}
\begin{aligned}
\begin{array}{ll}
    &~~y_{i,t} = \max \{x_{i,t},w_{i} \},i \in \{0,1,\dots,n\} \text{~~~~~~~~~~~if } x_{0,t} \leq w_0;
    \\ \ \ \\ \ \
    &\left \{
\begin{array}{ll}
     y_{i,t} &= \max \{x_{i,t},w_{i} \}, i \in \{0,1,\dots,n\}\backslash i^* \\
     y_{i^*,t} &= x_{i^*,t}
\end{array}
\right.  \text{~~~if } x_{0,t} > w_0, 
\end{array}
\end{aligned}
\end{equation}
where $i^* := \argmin_{i \in [n]} \{b_i-c_i \}.$
If the target inventory level of the warehouse is feasible, we set the order-up-to level of all installations greedily. On the other hand, if the target inventory level of the warehouse is infeasible, besides setting the order-up-to level in other installations in a greedy manner
we also halt the replenishment of store $i^*$, which is the store that benefits most from deliveries of the warehouse. This helps reduce the inventory in the warehouse since stock-outs at store $i^*$ would call for deliveries in warehouse inventory.

The following lemma (proved in Section~\ref{sec:app-lemma-transition-solver}) shows the above solver meets Definition~\ref{def:transition solver}.
\begin{lemma}\label{app3:le:transition solver}
  The solver $\zeta^{\mathrm T}$ defined in \eqref{app3transitionsolver} satisfies Definition~\ref{def:transition solver}  with $M = n+6\beta_0 \sum_{i=0}^n \rho_i$.
\end{lemma}

\subsection{Convexity and Smoothness of $Q(\by)$}\label{app3:convexity and smoothness}
We have the following convexity and smoothness properties of $Q(\by)$.

\begin{lemma}
\label{le:app3-properties-qy}
Under Assumption \ref{ass1:generalass:app3}, we have that $Q(\by)$ is convex and $\beta$-smooth, where $\beta = (n+4)^2\beta_0 \left(h_0 + \sum_{i=1}^n (h_i + b_i + c_i)\right) = \O(n^3)$.
\end{lemma}

We briefly discuss the novel technical ingredients in the proof of Lemma~\ref{le:app3-properties-qy}. The convexity proof is quite simple -- since $Q(\by)$ is the expected value of $C(\by, \bd)$, we only need to show the convexity of $C(\by, \bd)$. For the latter goal, we turn to analyze the dual of $C(\by, \bd)$, which is the maximum of a set of linear functions $\bm{\bm{\pi}}^\top(\bd'- H\by )$ and therefore convex.
Establishing the smoothness of $Q(\by)$, however, presents greater challenges. At a higher level, we adopt a sample-based approach, which is also used in \cite{wang1985distribution}. However, the general approach of \cite{wang1985distribution} could only achieve a smoothness bound of $2^{\O(n)}$. In contrast, we exploit the delivery dynamics of the one-warehouse multi-store system and come up with a clever decomposition of the cost function $C(\by, \bd)$. Each term in the decomposition turns out easier to deal with, and this finally helps us prove the $\O(n^3)$ smoothness bound. {We present a more detailed discuss about technical novelties in our smoothness proof in Section~\ref{sec:app3smoothnesshighlevel} in the supplementary materials.}
The formal proof of Lemma~\ref{le:app3-properties-qy} is included in Section~\ref{proofofsmoothnessofapp3}.

\subsection{Regret Analysis}\label{app3:regret}
Plugging the parameter constants in Lemmas~\ref{le:app3-properties-estimator}, \ref{app3:le:transition solver} and  \ref{le:app3-properties-qy} into the meta-result in Theorem~\ref{thm:inv,general}, we obtain the following regret upper bounds for one-warehouse and multi-store inventory systems. 

For convenience, we define notations 
 $\beta = (n+4)^2\beta_0 \left(h_0 + \sum_{i=1}^n (h_i + b_i + c_i)\right)$, $\bar{Q} = \max_{\by \in \Gamma} Q(\by)$, $R = \max_{\by,\by' \in \Gamma} \Vert \by-\by' \Vert_2$, $\sigma_0 = [(h_0 + \sum_{i=1}^n b_i)^2+ \sum_{i=1}^n(h_i+b_i)^2]^\frac 12$, and $M = n+6\beta_0\sum_{i=0}^n \rho_i$, 
in the following theorem.
\begin{theorem}\label{app3:thm:inv,general}
Under Assumption \ref{ass1:generalass:app3}, invoking the meta-policy with  {\sc Well-behaved Gradient estimator} \eqref{estimator:app3} and the {\sc Transition Solver} \eqref{app3transitionsolver}, and
 running the meta-policy with stepsize $\eta<1/\beta$ and 
  batch size $n_{\tau}= \lceil\sqrt{T}\rceil$ or $n_{\tau}= K\tau$, where $K$ is a positive integer, we have
$
 \mathcal R_{\mathrm{cvx}}(T) \leq \O(\sqrt{T})
$. Here the $\O(\cdot)$ notation hides the polynomial dependence on $\beta$,  $\bar{Q}$, $\eta$, $R$, $M$, $K$ and $\sigma_0$.
\end{theorem}

{\section{Numerical Experiments}\label{sec:numerical experiments}
We conduct extensive numerical experiments to demonstrate that our meta-policy enjoys many advantages such as low relative average regret, low variance of regret and distance from the optimal solution, and high computational efficiency. 

The theoretical analysis of our meta-policy is based on the continuous assumption of the demand distribution. In some scenarios,\footnote{{In an example of \cite{den2020discontinuous},  discontinuity of the demand function is mainly because there is a substantial mass of customers earning a minimum wage.
Similarly, in inventory control problems, demand distribution usually has a positive mass at $0$. However, this case can be addressed by the theoretical analysis of our meta-policy via a straightforward adaptation.}} the demand might have point masses, or even be discrete. To evaluate the numerical robustness of our meta-policy, we also conduct numerical experiments under discrete distributions, such as Poisson distribution and geometric distribution.

\noindent\underline{\textbf{Experimental settings.}} Due to space constraints, we now briefly introduce the settings of each experiment, and defer the detailed descriptions and the numerical results to Section~\ref{detailedexp}.
\begin{enumerate}[nolistsep] 
    \item {\it Single product newsvendor problem.} We consider the newsvendor problem and compare our meta-policy with the other four classical algorithms under four different demand distributions. The numerical results in Section~\ref{subsec:num-nvp} show that our meta-policy is competitive among many classic algorithms.
    \item {\it Application I: Multi-product and multi-constraint inventory system.}
We conduct numerical experiments for both the multi-product single constraint system management problem studied in~\cite{shi2016nonparametric} and the multi-product multi-constraint system management problem described in Section~\ref{app: multi-product}. We compare both the regret and running time of SGD and our meta-policy. Please refer to Section~\ref{subsec:num-mpmc} in the supplementary materials.

\item \textit{Application II: multi-echelon serial inventory system.}
We conduct numerical experiments for the multi-echelon serial inventory system management problem described in Section~\ref{app: multi-echelon}. Please refer to Section~\ref{subsec:num-me}.

\item \textit{Application III: one-warehouse and multi-store inventory system.}
We conduct numerical experiments for one-warehouse and multi-store inventory system described in Section~\ref{app: owms}. Please refer to Section~\ref{subsec:num-owms}.

\item \textit{Inventory system with positive lead times.}
For inventory systems with positive lead times, we consider learning the optimal base-stock policy. We compute
the sample path derivative as a gradient estimator of the long-run average function and conduct
numerical experiments on several demand distributions and system instances. Please refer to Section~\ref{subsec:num-leadtime}. 

\item We also consider the two-echelon problem described in~\cite{zhang2022no}. The authors proposed an algorithm that optimizes the decision of the retailer by SAA and optimizes the decision of the supplier by SGD. We replace the SGD subroutine by minibatch-SGD with the batch size scheme derived in our paper, and compare the new algorithm with the one proposed by~\cite{zhang2022no}. Numerical results in Section~\ref{subsec:num-te} show that the method based on our minibatch-SGD has a better numerical performance. We also establish regret analysis for the minibatch-SGD-based method.
\end{enumerate} 

\noindent\underline{\textbf{Numerical strengths of our meta-policy.}} We summarize the numerical strengths of our meta-policy demonstrated in the above numerical experiments as follows.
\begin{enumerate}[nolistsep]
    \item \textit{Competitive regret performance.} 
    Our meta-policy consistently demonstrates competitive regret performance among the wide-range of selected classic and baseline algorithms, which is consistent across our extensive experiment settings.
    \item \textit{High computational efficiency.}  Our meta-policy enjoys high computational efficiency due to fewer projection operations (infrequent decision updates). Through the running time test (please refer to Tables~\ref{tab:run1} and~\ref{tab:run2} in the supplementary materials), we observe that our method achieves about $15\sim 150$ times speed-up for different values of $T$ compared with the SGD-based algorithms. 
   \item \textit{Low variances.}  
      In the numerical experiment, we observe that the relative average regret and the distance between the optimal solution and the one learned by our meta-policy have smaller variances than the SGD-based algorithms in several experimental instances (see Figure~\ref{fig:nvp-var} in the supplementary materials).  
    \item \textit{Flexible applications.} We have demonstrated the flexibility of our meta-policy through three applications via theoretical analysis. In the numerical experiments, we apply our meta-policy to more applications such as inventory systems with positive lead times and the two-echelon systems considered in~\cite{zhang2022no}. The performance of our meta-policy in these two new applications further highlights the flexibility of our meta-policy framework.
\end{enumerate}}

\section{Conclusion and Future Directions}\label{sec:conclusion}
We propose a minibatch-SGD-based meta-policy to achieve joint learning and optimization for inventory systems with myopic optimal policies. Our meta-policy provides a new method to deal with the low-target-inventory-level issue encountered by gradient-based methods, which turns out to be more flexible than existing techniques in literature, especially when applied to more complicated inventory systems. We illustrate the power and flexibility of our meta-policy by applying it to a wide range of inventory control problems to derive new and optimal theoretical guarantees.

SGD-based algorithms and queueing-theory-based arguments have emerged to be the predominant methods for inventory control problems with demand learning. We believe our minibatch-SGD-based algorithm, which is more flexible, could also be applied to solve many other problems in inventory management beyond the scope of this paper. The following are three interesting future directions.

{
\noindent\underline{\bf Inventory system with positive lead times.}
For single-echelon lost-sales inventory systems with positive lead times, \cite{huh2009asymptotic} proved that the base-stock policies are asymptotically optimal. Thus, existing results mainly focus on learning optimal base-stock policies for single-echelon systems \citep{huh2009adaptive,zhang2020closing,agrawal2019learning}.

It is also hopeful to use our meta-policy (with possible modifications) to learn the optimal base-stock policies for the above single-echelon setting. To achieve this, we need to establish several key steps, which are listed as follows. First, we need to construct a gradient estimator for the long-run cost function to the base-stock parameter. Second, we need to prove the estimator has desired properties such as unbiasedness (in the long run) and bounded variance. Third, we need to prove the long-run cost function has desired properties such as convexity and smoothness. 
 To further illustrate the potential of our meta-policy in learning the optimal base-stock policy,  we compute the sample path derivative as a gradient estimator of the long-run average function and conduct numerical experiments on several demand distributions and system instances in Section \ref{subsec:num-leadtime} in the supplementary materials. All numerical examples consistently demonstrate a good convergence performance of our meta-policy. 
 Since the long-run cost functions depend on the steady-state inventory level under a specific base-stock policy and have no closed-form expressions, tackling the second and third steps would call for heavy machinery that is tailored for inventory systems with positive lead times (such as analyzing a complex and problem-specific Markov chain). 
By addressing these technical challenges, we would be able to establish the theoretical guarantee of our meta-policy for single-echelon inventory systems with deterministic lead times. Moreover, the simplicity of our meta-policy suggests its applicability to the more complex setting of stochastic lead times which is less studied in literature. We leave these questions as an interesting future work.
}

{
\noindent\underline{\bf Feature-based inventory management problems.} Another important direction is to consider feature-based problems, where demand depends on (usually high-dimensional) features and decision-makers can observe the feature information to make better inventory decisions. \cite{ding2021feature} considered feature-based inventory control problems under censored demand and linear assumptions. They used a gradient method to update the estimation of linear model parameters and queueing-theory-based methods to establish regret analysis.  To better exploit the high-dimensional data, many studies consider the clustering of features. \cite{aouad2023market} proposed tree-based methods and \cite{keskin2023data} design spectral methods for problems in different areas. 
It would be interesting to combine our minibatch meta-policy with the above-mentioned techniques to address more feature-based inventory management problems. Since feature vectors usually change over time, one technical challenge for future research is to design appropriate batch schemes in our meta-policy to deal with the changing features in a single batch.

\noindent\underline{\bf Inventory system with nonstationary demand.}  
Our paper and most classic works assume the demand distribution is stationary over time, i.e., for all $t$, $D_t$ follows the same distribution. 
However, 
nonstationarity captures broader practical scenarios, and therefore it has attracted more attention in the field of inventory management in recent years. 
In the literature, there are mainly two types of nonstationarity settings,
1) each change-point in parameters results in a minimum positive shift \citep{keskin2022data};
2) there is a variation budget on the total variation \citep{keskin2021nonstationary}. It would be quite interesting to extend our minibatch meta-policy to handle the nonstationary demand, yet it would lead to a more challenging problem due to the changing optimal order-up to level. 
}








\bibliographystyle{informs2014}
\bibliography{main.bib}

\ECSwitch


\begin{center}
 		\Large{\bf Supplementary materials to}\\ ``A Minibatch-SGD-Based Learning Meta-Policy for Inventory Systems with Myopic Optimal Policy'' 
 	\end{center}
 	
 	\vspace{10pt}

\section{Formulation and Theoretical Results for Application II: Multi-Echelon Serial Inventory System} \label{ecapp:multi-echelon}

{In Section~\ref{app: multi-echelon}, we briefly introduced our second application. In this section, we present the detailed problem formulations and technical results.}
\subsection{Problem Formulation}\label{app2:problem formulation}
\noindent \underline{\bf Sequence of events.}
In each period $t \in [T]$, we assume the following sequence of events:
\begin{enumerate}[nolistsep]
    \item At the beginning of period $t$, the firm observes the on-hand inventory of $n$ stages $\bx_t = (x_{1,t},x_{2,t},\dots,x_{n,t})^{\top}$. The initial on-hand inventory at $n$ stages is given as $\bx_1$.
    \item The firm decides the order-up-to level of $n$ stages  $\by_t \geq \bx_t$. The order-up-to level is restricted by the following inequalities
\[
\by_t \in \Gamma = \{\by\in \mathbb{R}^{n}_+ \mid \by \leq \bm{\rho}\},
\]
where $\bm \rho = (\rho_1,\rho_2,\dots,\rho_n)^{\top} \in \mathbb{R}^n_+$ is the capacity constraint of each stage.
    \item The demand $D_t\sim D$ is realized in Stage $1$, denoted by $d_t$, which is firstly satisfied to the maximum extent using the on-hand inventory of Stage $1$. If Stage $1$ is out of stock, it could get replenishment from Stage $2$ by paying the per-unit transportation fares $b_1$. If Stage $2$ is out of stock, it could get replenishment from Stage $3$ by paying the per-unit transportation fare $b_2$, and so on. Demand that can not be satisfied by all stages is lost, and the firm only observes the sales data (censored demand),  $s_t=\min (d_t,y_{1,t}+\dots+y_{n,t})$. The leftover inventory level 
    $$\bx_{t+1}=\mathfrak f(\by_t,d_t):=((y_{1,t}-d_t)^+,(y_{2,t}-(d_t-y_{1,t})^+)^+,\dots,(y_{n,t}-(d_t-y_{1,t}-\dots-y_{n-1,t})^+)^+)^{\top}.$$
      \end{enumerate} 
   At the end of period $t$, the firm suffer inventory cost \begin{align}\nonumber
        &\underbrace{(h_1+h_2+\dots+h_n)(y_{1,t}-d_t)^+ }_{\text{holding cost in Stage $1$}}+
        \underbrace{b_1\min\{y_{2,t}+y_{3,t}+\dots + y_{n,t},(d_t-y_{1,t})^+\} }_{\text{transportation cost from Stage $2$ to Stage $1$}}\\&+\nonumber
     \underbrace{(h_2+h_3+\dots+h_n)(y_{2,t}-(d_t-y_{1,t})^+)^+}_{\text{holding cost in Stage $2$}}+
                \underbrace{b_2\min\{y_{3,t}+y_{4,t}+\dots + y_{n,t},(d_t-y_{1,t}-y_{2,t})^+\} }_{\text{transportation cost from Stage $3$ to Stage $2$}}
                \\&+\dots+\label{eq:complex formulation}
     \underbrace{h_n(y_{n,t}-(d_t-y_{1,t}-\dots-y_{n-1,t})^+)^+}_{\text{holding cost in Stage $n$}}+
       \underbrace{(b_1+b_2+\dots+b_n)(d_t-y_{1,t}-\dots-y_{n,t})^+}_{\text{lost sales cost}},
    \end{align}
 where $b_i$ and $h_i$ are positive constants for all $i\in [n]$,  $\sum_{k=i}^n h_k$ is the per unit holding cost in Stage $i$ for $i\in[n]$,\footnote{In practice, the holding cost of the higher stage is greater than the holding cost of the lower stage, and thus for convenience, we denote the per unit holding cost of stages in this way.} $b_i$ is the per unit transportation cost from Stage $i+1$ to Stage $i$ for $i\in[n-1]$, and $(b_1+\dots+b_n)$ is the per unit lost sales cost.\footnote{Assume the lost sales cost is greater than the total transportation cost from Stage $n-1$ to Stage 1, we use $(b_1+\dots+b_n)$ to denote the lost sales cost.}  

Although the definition of the inventory cost (Eq.~\eqref{eq:complex formulation}) is clear and practical, it seems that the above cost formulation is too complex to analyze. Fortunately,  we find that it could be simplified to an equivalent formulation (Eq.~\eqref{eq:simple formulation}) by the following derivation.

For any $i\in[n-1]$, by discussing three cases $d_t\leq y_{1,t}+y_{2,t}+\dots+y_{i-1,t}$,  $d_t\geq y_{1,t}+y_{2,t}+\dots+y_{i,t}$ and   $y_{1,t}+y_{2,t}+\dots+y_{i-1,t}< d_t< y_{1,t}+y_{2,t}+\dots+y_{i,t}$ individually, we could obtain \begin{align*}
    (y_{1,t}+\dots+y_{i,t}-d_t)^++(y_{i+1,t}- (y_{1,t}+y_{2,t}+\dots+y_{i,t}-d_t)^+)^+  = (y_{1,t}+y_{2,t}+\dots+y_{i+1,t}-d_t)^+.
\end{align*}
By the above equation the sum of all the holding cost terms multiplied with the $h_i$ equals to
$$h_i(y_{1,t}+\dots+y_{i,t}-d_t)^+.$$
Therefore, collecting all holding cost terms gives the total holding cost
$$\sum_{i=1}^n[h_i(y_{1,t}+\dots+y_{i,t}-d_t)^+].$$
Similarly, we could also get the sum of the total transportation cost and the lost sale cost as 
\[\sum_{i=1}^n[b_i(d_t-(y_{1,t}+\dots+y_{i,t}))^+].\]
Therefore, the inventory cost \eqref{eq:complex formulation} is equivalent to
    \begin{align}\label{eq:simple formulation}
\sum_{i=1}^n[h_i(y_{1,t}+\dots+y_{i,t}-d_t)^++b_i(d_t-(y_{1,t}+\dots+y_{i,t}))^+].
    \end{align}
The goal of the firm is to design a policy
to minimize the following expected total inventory cost  :
\begin{equation}\nonumber
\E\left[\sum_{t=1}^TQ(\by_t)\right] = \sum_{t=1}^T\sum_{i=1}^n\E[h_i(y_{1,t}+y_{2,t}+\dots+y_{i,t}-D)^++b_i(D-(y_{1,t}+y_{2,t}+\dots+y_{i,t}))^+].
\end{equation}

\noindent \underline{\bf Regularity assumptions about demand.}
Almost the same as the assumptions made for the first application in Section~\ref{app: multi-product}, we make the following two assumptions for the multi-echelon serial inventory system. Please refer to the remark and discussion of the assumptions in Section~\ref{app: multi-product} for a detailed discussion of Assumptions~\ref{ass1:generalass:app2} and \ref{ass2:strongly:app2}.
\begin{assumption}
\label{ass1:generalass:app2}
In this application, we make the following basic assumptions on demand.
\begin{enumerate}[nolistsep]
    \item $D_{t}$ is i.i.d.~across time period $t \in [T]$ with common distribution $D$.
    \item $D$ is a continuous random variable with  PDF $f_{D}(\cdot)$ and $f_{D}(x) \leq \beta_0$ for some $\beta_0 >0$
\end{enumerate}
\end{assumption}
The following assumption is proposed to improve the regret performance of our meta-policy.
\begin{assumption}
\label{ass2:strongly:app2}
  The PDF of $D$  is lower bounded by a positive constant. That is $f_{D}(x) \geq \alpha_0$ for some $\alpha_0 >0$.
\end{assumption}

\subsection{Design of the {\sc Well-behaved Gradient Estimator} and the {\sc Transition Solver}}\label{app2:design}
As mentioned before, we denote the expected one-period  cost when the after-delivery inventory level is $\by$ as $$Q(\by) = \sum_{i=1}^n\E[h_i(y_{1}+\dots+y_{i}-D)^++b_i(D-(y_{1}+\dots+y_{i}))^+].$$
For $i\in[n]$, we define $\tilde y_i= \sum_{k=1}^i y_k$  and $\tilde y_{i,t}= \sum_{k=1}^i y_{k,t}$. Let $\tilde \by=(\tilde y_1,\tilde y_2,\dots,\tilde y_n)^{\top}$ and  $\tilde \by_t=(\tilde y_{1,t},\tilde y_{2,t},\dots,\tilde y_{n,t})^{\top}$. We denote the linear one-to-one correspondence between $\by$ and $\tilde \by$ by a $n\times n$ matrix $$B = 
  \begin{bmatrix}
  1 & 0 & 0& 0 &\cdots & 0\\
  1 & 1 & 0& 0&\cdots & 0  \\
  1 & 1 & 1& 0&\cdots & 0  \\
   \vdots &\vdots &\vdots& \vdots &\vdots & \vdots  \\
    1 & 1& 1& 1 &\cdots & 1  
  \end{bmatrix},$$  and thus we have $\tilde\by = B\by$. With $\bh=(h_1, h_2,\dots,h_n)^{\top}$ and $\bb = (b_1, b_2,\dots,b_n)^{\top}$, we could define
\begin{align*}
    \tQ(\tby) &:= \sum_{i=1}^n\E[h_i(\ty_{i}-D)^++b_i(D-\ty_{i})^+]\\&= \bh^\top \E[(\tby - D\cdot\bm{1})^+] + \bb ^\top \E[(D\cdot\bm{1}-\tby)^+],
\end{align*}
Therefore, minimizing $\E\left[\sum_{t=1}^TQ(\by_t)\right]$ over $\Gamma$ is equivalent to minimize
$\E\left[\sum_{t=1}^T\tQ(\tby_t)\right]$ over $\tilde\Gamma = \{B\by\mid \by\in\Gamma\}$. We only need to conduct minibatch SGD update method on  $\tQ(\tby_t)$ to obtain $\tby_{t+1}$,  and get the target level through $\by_{t+1} =  B^{-1} \tby_{t+1}$.

Since the  gradient of $\tQ(\tby)$ at $\tby$ is given by
$$\nabla \tQ(\tby) = \E \left[ \bh\odot\mathbb{I}[\tby>D\cdot\bm{1}]-\bb\odot\mathbb{I}[\tby\leq D\cdot\bm{1}] \right],$$ assuming $d$ is a realization of $D$, it is easy to 
 design the following gradient estimator of $\tQ(\tby)$
 \begin{equation}\label{gradientestimatorapp2}
     \hat\nabla \tQ(\tby) =\bh\odot\mathbb{I}[\tby>d\cdot\bm{1}]-\bb\odot\mathbb{I}[\tby\leq d\cdot\bm{1}].
 \end{equation}

 By the same proof as Lemma~\ref{le:app1properties}, we have the following lemma.
 \begin{lemma}\label{le:app2properties}
Properties of $\tQ(\tby)$ and $\hat\nabla \tQ(\tby)$:
\begin{enumerate}[nolistsep]
    \item Under 
Assumption \ref{ass1:generalass:app2}, $\tQ(\tby)$ is convex and $(\max_{i\in[n]}(h_i+b_i)\beta_0)$-smooth.
    \item Under Assumptions \ref{ass1:generalass:app2} and \ref{ass2:strongly:app2}, $\tQ(\tby)$ is $(\min_{i\in[n]}(h_i+b_i)\alpha_0)$-strongly convex and $(\max_{i\in[n]}(h_i+b_i)\beta_0)$-smooth.
    \item The norm of gradient is bound as $\|\nabla \tQ(\tby)\|_2\leq \sigma_0$, where $\sigma_0:=(\|\bh\|_2^2+\|\bb\|_2^2)^\frac 12$.
    \item The estimator $\hat\nabla \tQ(\tby)$ is unbiased, i.e.,$\E\left[\hat\nabla \tQ(\tby)\right] =   \nabla \tQ(\tby) $.
    \item It holds that $\E\left[\|\hat\nabla \tQ(\tby) -\nabla \tQ(\tby)\|_2^2\right] \leq \E\left[\|\hat\nabla \tQ(\tby)\|_2^2\right] \leq \|\bh\|_2^2+\|\bb\|_2^2=\sigma_0^2.$
\end{enumerate}
\end{lemma}

The last two items of Lemma~\ref{le:app2properties} show that the gradient estimator $\hat\nabla \tQ(\tby)$ is a {\sc Well-behaved Gradient Estimator} satisfying Definition~\ref{def:gradient estimator}.
We now turn to the design of the {\sc Transition Solver}. 

Our {\sc Transition Solver} $\zeta^{\mathrm T}(\bx_t,\bw)$ takes a on-hand inventory level $\bx_t$ and some infeasible target inventory level $\bw$ as inputs,  and output an actually implemented inventory level $\by_t$ as follows: 
\begin{equation}\label{app2:opt}
\begin{aligned}
    &\left \{
\begin{array}{ll}
     y_{i,t} &= w_{i}, ~~i \in [n] \text{ and }i>i^*; \\ \ \ \\
     y_{i^*,t} &= x_{i^*,t}, i \in [n] \text{ and }i\leq i^*,
\end{array}
\right.   
\end{aligned}
\end{equation}
where $i^* := \max\{i \in [n]\mid w_{i}<x_{i,t}\}.$ Specifically, we set the order-up-to level of stages before $i^*$ as the target inventory level, and set the order quantity of stages behind $i^*$ to $0$, where $i^*$ is the highest stage where the target inventory level is infeasible.

The following lemma shows the solver designed above meets Definition~\ref{def:transition solver}.

\begin{lemma}\label{app2:le:transition solver}
  The solver $\zeta^{\mathrm T}$ defined in \eqref{app2:opt} satisfies Definition~\ref{def:transition solver}  with $M = 1 + 6\beta_0 \sum_{i=1}^n \rho_i$.
\end{lemma}
\proof{Proof of Lemma~\ref{app2:le:transition solver}.}
Consider the worst case that the target inventory level of stage $n$ is infeasible. In this case, all stages stop replenishing. Noting that each stage has at most $\rho_i$ inventory capacity, by Lemma \ref{lem:hittingtime}, we can bound $\E[s(\bx_0,\bw)]$ as follows 
$$\E [s(\bx_0,\bw)] \leq 1 + 6\beta_0 \sum_{i=1}^n \rho_i,$$
where we complete the proof of this lemma.
\Halmos
\endproof

\subsection{Regret Analysis}\label{app2:regret}
Plugging the parameter constants in Lemmas~\ref{le:app2properties} and \ref{app2:le:transition solver} to the meta-results (Theorems~\ref{thm:inv,general} and \ref{thm:inv,str}), we obtain the following regret upper bounds for multi-echelon serial inventory systems with censored
demand  in two cases. 

For convenience, we define notations 
$\beta = \max_{i\in[n]}(h_i+b_i)\beta_0$,  $\alpha = \min_{i\in[n]}(h_i+b_i)\alpha_0$, $\bar{Q} = \max_{\by \in \Gamma} Q(\by)$, $R = \max_{\by,\by' \in \tilde\Gamma} \Vert \by-\by' \Vert_2$,  $\sigma_0=(\|\bh\|_2^2+\|\bb\|_2^2)^\frac 12$ and $M=n + 6n\beta_0 R$ in the following two theorems.
\begin{theorem}\label{app2:thm:inv,general}
Under Assumption \ref{ass1:generalass:app2}, invoking the meta-policy with  {\sc Well-behaved Gradient estimator} \eqref{gradientestimatorapp2} and the {\sc Transition Solver} \eqref{app2:opt}, and
 running the meta-policy with stepsize $\eta<1/\beta$ and 
  batch size $n_{\tau}= \lceil\sqrt{T}\rceil$ or $n_{\tau}= K\tau$, where $K$ is a positive integer, we have
$
 \mathcal R_{\mathrm{cvx}}(T) \leq \O(\sqrt{T})
$. Here the $\O(\cdot)$ notation hides the polynomial dependence on $\beta$,  $\bar{Q}$, $\eta$, $R$, $M$, $K$ and $\sigma_0$.
\end{theorem}

\begin{theorem}
\label{app2:thm:inv,str}
Under Assumptions \ref{ass1:generalass:app2} and \ref{ass2:strongly:app2}, invoking the meta-policy with  {\sc Well-behaved Gradient estimator} \eqref{gradientestimatorapp2} and the {\sc Transition Solver} \eqref{app2:opt}, and
running the meta-policy with stepsize $\eta\leq \min\{\alpha/2,1/\alpha,1/2\beta\}$ and batch size  $n_{\tau}=\lceil\varsigma^{\tau-1}\rceil$, where  $\varsigma= 1/{\gamma}$ and $\gamma=1-\eta\alpha+2\eta^2\in(0,1)$, we have $\mathcal R_{\mathrm{str}}(T)  \leq \O(\log{T})$. Here the $\O(\cdot)$ notation hides the polynomial dependence on $\beta$, $\alpha$, $\varsigma$, $\bar{Q}$, $R$, $M$, and $\sigma_0$.
\end{theorem}
{\subsection{Discussions on Multi-echelon inventory Systems}\label{sec:Discussions on Multi-echelon inventory Systems}
\noindent\underline{\bf Differences from \cite{yang2022non}.}
\cite{yang2022non} considered a periodic-review multi-echelon inventory problem with zero lead time. In each period, the downstream stage can get replenishment from the upstream stage, and backorder is assumed if the inventory in the upstream stage is not enough. The algorithm they designed is based on the online gradient algorithm proposed by \cite{huh2014online}, which considered an online optimization problem satisfying a generalized notion of convexity, called sequentially convexity. To complete the regret analysis, they introduced two auxiliary systems combined with queueing-theory-based
methods.

Our multi-echelon setting is different from their works in the following perspectives. First, 
in contrast to the setting that only the upmost stage could get replenishment from an external supplier with ample stock in \cite{yang2022non}, we allow all stages to get replenishment from an external supplier at the beginning of each period. 
Second, we consider emergency transportation with additional fares from upstream stages to downstream stages, and unsatisfied demand is lost. 
Besides, the cost function of our problem is different from the one in \cite{yang2022non} and is proved to be convex after a variable transformation.}

\noindent\underline{\bf Difficulties of considering positive lead times.}
 For multi-echelon lost-sales inventory systems with positive lead times, the optimal policy is quite complicated and computationally intractable even if the demand distribution is known \citep{snyder2019fundamentals}.
Additionally, to the best of our knowledge, no simple heuristics (such as base-stock or constant-order policies) are proved to be asymptotically optimal so as to qualify a reasonable benchmark for learning algorithms in such inventory systems. Even if we use the base-stock policy as a benchmark, the cost function is not convex in the base-stock vector \citep{huh2010base}. In all, learning multi-echelon lost-sales inventory systems with positive lead times presents multiple challenges due to its complex and not-well-understood solution structure. Making progress towards this direction might call for new techniques that are quite different from this paper.

\section{Discussion on Differences from Existing Cycle-update Policies.}\label{sec:discussion on difference}

Several cycle-update learning algorithms that operate in cycles and use the same target variable in each cycle have been proposed for different settings of inventory management problems in the literature. However, the ideas and analysis of their methods are quite different from ours. 
In the following, we carefully discuss the differences from a methodology perspective to highlight the novelty of our inventory control meta-policy.

The work by \cite{huh2009adaptive} studies lost-sales inventory systems with positive lead times, and proposes the first cycle-update policy in the inventory control literature that achieves a sublinear (but suboptimal) regret $\mathcal O(T^{2/3})$. 
The key difference between \cite{huh2009adaptive} and our meta-policy is about the different technical challenges faced by the two works that call for the cycle/minibatch update policies. In \cite{huh2009adaptive},  since the gradient of their long-run cost function depends on the steady-state inventory level under a specific base-stock policy, using gradient information from each period directly to construct the gradient estimator introduces substantial bias. Therefore, the authors keep the target variable, i.e., the base-stock level unchanged throughout a predetermined sufficiently long cycle, and utilize the gradient information from the single last period of each cycle to derive a gradient estimator with minimal bias and perform policy updates.
In contrast, our meta-policy is conceived to reduce the frequency of target infeasibility through a low-switching policy update scheme. Moreover, our meta-policy utilizes gradient information from all periods to reduce the variance of our unbiased gradient estimators.

We are also aware that similar cycle-update designs were adopted by \cite{zhang2018perishable}, \cite{chen2020optimal}, and \cite{zhang2020closing} in their learning algorithms across different settings.
Specifically, 
\cite{zhang2018perishable} pioneered the development of the first nonparametric learning algorithm for perishable inventory systems. \cite{chen2020optimal} explored backlogging inventory systems with random capacity and unknown demand. \cite{zhang2020closing} introduced a novel simulated cycle-update policy, successfully establishing an optimal $\mathcal O(\sqrt{T})$ regret for the inventory management problem considered in \cite{huh2009adaptive}.

Our meta-policy is mainly technically different from the above three works from the following two aspects.

First, the cycle/minibatch-length patterns are different. In the above three works, new cycles are activated by specific triggering events dependent on sample paths. As a result, the cycles have random lengths and starting periods, which also depend on the algorithms' execution trajectory.\footnote{{The triggering events or random cycles are problem-specific and serve various purposes. For instance, \cite{zhang2018perishable} defined the triggering event when stockout occurs to ensure the target base-stock level can always be attained, while \cite{chen2020optimal} defined the random cycle primarily to deal with the non-convexity of the cost function.}}In contrast, the epoch lengths in our meta-policy are fixed and pre-determined. Besides, the cycle length in these works is a constant in expectation, which means that the number of cycles therein is $\mathcal O(T)$. However,  the number of epochs in our meta-policy is designed to be either $\mathcal O(\sqrt{T})$ or $\mathcal O(\log T)$ according to the optimal regret in different settings. Therefore, our meta-policy has a low-switching property, which helps to control the frequency of the lower target inventory levels.

Second, the update of the target variables is based on different cost functions. In the above three works,
the authors initially formulated cycle cost functions by aggregating the costs across all periods in a cycle. Then, they establish convexity and construct gradient estimators for these cycle cost functions, which involve additional (and very important) technical efforts.\footnote{{It is important to note that our discussion focuses solely on the similarity observed in the gradient descent update component of their algorithms. The complexity of their considered inventory systems necessitates the incorporation of additional algorithmic techniques. For instance, \cite{zhang2020closing} introduced a sophisticated and effective design of an auxiliary simulated system to assist their cycle-update policy.}} Their target variable updates employ SGD on the cycle cost functions.
In contrast, in our meta-policy we directly deal with the original cost functions and establish their convexity and gradient estimators. We conduct minibatch SGD directly on the original cost functions to update the target inventory level.

\section{Proofs Omitted in Section~\ref{sec:regret of meta policy}}\label{sec:proofs of metapolicy}
\subsection{Lower Bounds}\label{subsec:lower bound and maximum switching budget}

We first present the lower bound for the case that the expected cost function $Q(\by)$ is smooth and convex in the following proposition.
\begin{proposition}\label{prop:cvxlowerbound}
There exist problem instances of our inventory system framework whose expected cost function $Q(\by)$ is smooth and convex, such that the regret of any learning algorithm is lower bounded by $\Omega(\sqrt{T})$. 
\end{proposition}
\begin{remark}
  The problem instances satisfying the above proposition have been constructed by \cite{zhang2020closing} based on the discrete demand instances of \cite{besbes2013implications}  for the single product newsvendor problem. Since the newsvendor problem is a special case of our framework, it also provides proof for our proposition. 
\end{remark}
The following proposition establishes the lower bound for the case that the expected cost function $Q(\by)$ is smooth and strongly convex.
\begin{proposition}\label{prop:strcvxlowerbound}
There exist problem instances of our inventory system framework whose expected cost function $Q(\by)$ is smooth and strongly convex, such that the regret of any learning algorithm is lower bounded by $\Omega(\log{T})$. 
\end{proposition}
\begin{remark}
 The problem instances satisfying the above proposition have been constructed by \cite{besbes2013implications} for the single product newsvendor problem. Similarly, it also provides proof for our proposition. 
\end{remark}
\subsection{Proof of Lemma~\ref{le:cvx}}
\proof{Proof.}
The $\beta$-smooth property of $F(\bw)$ gives
\begin{align}
\nonumber&F(\tbw_{\tau+1})-F(\tbw_{\tau})
\\&\leq \nabla F(\tbw_{\tau})^{\top}\left(\tbw_{\tau+1}-\tbw_{\tau}\right)+\frac{\beta}{2}\left\|\tbw_{\tau+1}-\tbw_{\tau}\right\|_2^2 \nonumber \\ \nonumber
&=g(\tbw_{\tau};\bm{\xi}_t)^{\top}\left(\tbw_{\tau+1}-\tbw_{\tau}\right)+\left(\nabla F(\tbw_{\tau})-g(\tbw_{\tau};\bm{\xi}_t)\right)^{\top}\left(\tbw_{\tau+1}-\tbw_{\tau}\right)+\frac{\beta}{2}\left\|\tbw_{\tau+1}-\tbw_{\tau}\right\|_2^2 \\ 
&\leq g(\tbw_{\tau};\bm{\xi}_t)^{\top}\left(\tbw_{\tau+1}-\tbw_{\tau}\right)+\left\|\nabla F(\tbw_{\tau})-g(\tbw_{\tau};\bm{\xi}_t)\right\|_2\|\tbw_{\tau+1}-\tbw_{\tau}\|_2+\frac{\beta}{2}\left\|\tbw_{\tau+1}-\tbw_{\tau}\right\|_2^2,
\label{eq:le:cvx:smoothness}
\end{align}
where the second inequality is due to Cauchy-Schwarz inequality.

Recalling the updating rule $\tbw_{\tau+1}=\Pi_\Gamma[\tbw_{\tau}-\eta g(\tbw_{\tau};\bm{\xi}_{\tau})]$, by the property of projection onto a convex set, we have 
\begin{equation}\label{eq:le:cvx:projection}
    (\bw^*-\tbw_{\tau+1})^\top(( \tbw_{\tau}-\eta g(\tbw_{\tau};\bm{\xi}_{\tau}))-\tbw_{\tau+1})\leq 0.
\end{equation}
And it is easy to note that 
\begin{equation}\label{eq:le:cvx:simple}
    \|\tbw_{\tau}-\bw^*\|_2^2 = \|\tbw_{\tau+1}-\bw^*\|_2^2+\|\tbw_{\tau+1}-\tbw_{\tau}\|_2^2 - 2(\bw^*-\tbw_{\tau+1})^\top(\tbw_{\tau}-\tbw_{\tau+1}).
\end{equation}
Combining \eq{\ref{eq:le:cvx:projection}} and \eq{\ref{eq:le:cvx:simple}}, we have 
\begin{align*}
    \eta g(\tbw_{\tau};\bm{\xi}_{\tau})^\top(\tbw_{\tau+1}-\bw^*)&\leq -(\tbw_{\tau+1}-\bw^*)^\top(\tbw_{\tau+1}-\tbw_{\tau})\\&=   \frac{1}{2}  (\|\tbw_{\tau}-\bw^*\|_2^2 -\|\tbw_{\tau+1}-\bw^*\|_2^2-\|\tbw_{\tau+1}-\tbw_{\tau}\|_2^2).
\end{align*}
And thus, we can obtain that
\begin{align*}
    \eta g(\tbw_{\tau};\bm{\xi}_{\tau})^\top(\tbw_{\tau+1}-\tbw_{\tau})&\leq  \eta g(\tbw_{\tau};\bm{\xi}_{\tau})^\top(\bw^*-\tbw_{\tau}) +  \frac{1}{2}  (\|\tbw_{\tau}-\bw^*\|_2^2 -\|\tbw_{\tau+1}-\bw^*\|_2^2-\|\tbw_{\tau+1}-\tbw_{\tau}\|_2^2).
\end{align*}
Combining the above inequality with \eq{\ref{eq:le:cvx:smoothness}}, we could have
\begin{align}
\nonumber&F(\tbw_{\tau+1})-F(\tbw_{\tau})
\\&\nonumber\leq g(\tbw_{\tau};\bm{\xi}_{\tau})^{\top}\left(\bw^*-\tbw_{\tau}\right)+\frac{1}{2\eta}  (\|\tbw_{\tau}-\bw^*\|_2^2 -\|\tbw_{\tau+1}-\bw^*\|_2^2)\\&\nonumber\quad+\left\|\nabla F(\tbw_{\tau})-g(\tbw_{\tau};\bm{\xi}_{\tau})\right\|_2\|\tbw_{\tau+1}-\tbw_{\tau}\|_2+\frac{\eta\beta-1}{2\eta}\left\|\tbw_{\tau+1}-\tbw_{\tau}\right\|_2^2
\\&\leq g(\tbw_{\tau};\bm{\xi}_{\tau})^{\top}\left(\bw^*-\tbw_{\tau}\right)+\frac{1}{2\eta}  (\|\tbw_{\tau}-\bw^*\|_2^2 -\|\tbw_{\tau+1}-\bw^*\|_2^2)+\frac{\eta\left\|\nabla F(\tbw_{\tau})-g(\tbw_{\tau};\bm{\xi}_{\tau})\right\|^2_2}{2(1-\eta\beta)},\label{eq:le:cvx:varianceform}
\end{align}
where the last inequality is because $\eta\beta\leq 1$ and $at^2+bt\leq -b^2/4a$ for all $t\in \R$ if $a<0$.

Since $F(\bw)$ is  convex, we could have
\begin{align}
F(\bw^*)-F(\tbw_{\tau})
\geq \nabla F(\tbw_{\tau})^{\top}\left(\bw^*-\tbw_{\tau}\right), \label{eq:le:cvx:cvx}
\end{align}

Combining \eq{\ref{eq:le:cvx:varianceform}} and \eq{\ref{eq:le:cvx:cvx}} with the fact that $\E[g(\tbw_{\tau};\bm{\xi}_{\tau})]= \nabla F(\tbw_{\tau})$, we could obtain
\begin{align}
   \nonumber &\E[F(\tbw_{\tau+1})]
\\&\leq  \E[F(\tbw_{\tau})+\nabla F(\tbw_{\tau})^{\top}\left(\bw^*-\tbw_{\tau}\right)]+\frac{1}{2\eta} \E[ \|\tbw_{\tau}-\bw^*\|_2^2 -\|\tbw_{\tau+1}-\bw^*\|_2^2]+\frac{\eta\E[\left\|\nabla F(\tbw_{\tau})-g(\tbw_{\tau};\bm{\xi}_{\tau})\right\|^2_2]}{2(1-\eta\beta)} \nonumber 
\\&\leq F(\bw^*)+\frac{1}{2\eta}  \E[\|\tbw_{\tau}-\bw^*\|_2^2 -\|\tbw_{\tau+1}-\bw^*\|_2^2]+\frac{\eta\E[\|\nabla F(\tbw_{\tau})-g(\tbw_{\tau};\bm{\xi}_{\tau})\|^2_2]}{2(1-\eta\beta)}. \label{eq:le:form00}
\end{align}

It is easy to obtain that
\begin{align*}
    &\E[\left\|\nabla F(\tbw_{\tau})-g(\tbw_{\tau};\bm{\xi}_{\tau})\right\|_2^2]\\&= \E\left[\left\|\nabla F(\tbw_{\tau})-\frac{1}{n_{\tau}} \sum_{i=1}^{n_\tau} \nabla f\left(\tbw_{\tau} ; \xi_{i,\tau}\right)\right\|_2^2\right]\\&= \frac{1}{n_{\tau}^2}\sum_{i=1}^{n_{\tau}} \E[\left\|\nabla F(\tbw_{\tau})-\nabla f\left(\tbw_{\tau} ; \xi_{i,\tau}\right)\right\|_2^2]+ \frac{2}{n_{\tau}^2}\sum_{1\leq i<j\leq n_{\tau}} \E[\braket{\nabla F(\tbw_{\tau})-\nabla f(\tbw_{\tau} ; \xi_{i,\tau}),\nabla F(\tbw_{\tau})-\nabla f(\tbw_{\tau} ; \xi_{j, \tau})}].
\end{align*}
Since $\E[\nabla f(\tbw_{\tau} ; \xi_{i,\tau})]= \nabla F(\tbw_{\tau})$ and  $\xi_{i,\tau}$ and $\xi_{j, \tau}$ are independent for $i\neq j $, we have 
\begin{align*}
 \E[\braket{\nabla F(\tbw_{\tau})-\nabla f(\tbw_{\tau} ; \xi_{i,\tau}),\nabla F(\tbw_{\tau})-\nabla f(\tbw_{\tau} ; \xi_{j, \tau})}]= 0.
\end{align*}
Therefore, we have
\begin{align}
   \nonumber \E[\left\|\nabla F(\tbw_{\tau})-g(\tbw_{\tau};\bm{\xi}_{\tau})\right\|_2^2]&= \frac{1}{n_{\tau}^2}\sum_{i=1}^{n_{\tau}} \E[\left\|\nabla F(\tbw_{\tau})-\nabla f\left(\tbw_{\tau} ; \xi_{i,\tau}\right)\right\|_2^2] \\&\leq \frac{\sigma^2}{n_{\tau}}.\label{eq:le:var}
\end{align}
Combing \eq{\ref{eq:le:form00}} and \eq{\ref{eq:le:var}}, we could have the following \textit{descent lemma}
\begin{align}\nonumber
    &\E[F(\tbw_{\tau+1})-F(\bw^*)]
\\&\leq  \frac{1}{2\eta}\E[ (\|\tbw_{\tau}-\bw^*\|_2^2 -\|\tbw_{\tau+1}-\bw^*\|_2^2)]+\frac{\eta\sigma^2}{2n_{\tau}(1-\eta\beta)}. \label{eq:descent lemma cvx}
\end{align}

\noindent\underline{\it Proof of the fixed-time batch scheme.}
Since $n_{\tau}=\lceil\sqrt{T}\rceil$, from the definition of $\tau_{\mathrm{max}}$ we have $\tau_{\mathrm{max}}\leq \lceil\sqrt{T}\rceil$. Noting the fact that $\max_{\bw_1,\bw_2\in\Gamma}\|\bw_1-\bw_2\|_2\leq R$, we obtain
\begin{align*}
    \sum_{\tau=1}^{\tau_{\mathrm{max}}}\E[ n_{\tau+1}(F(\tbw_{\tau+1})-F(\bw^*))]\leq\frac{\lceil\sqrt{T}\rceil R^2}{2\eta}+\frac{\eta\sigma^2\lceil\sqrt{T}\rceil}{2(1-\eta\beta)}.
\end{align*}
Since $\max_{\bw\in\Gamma}\|\nabla F(\bw)\|_2\leq G$, by the quasi-mean value theorem  we could have 
\begin{align*}
    F(\tbw_1)-F(\bw^*)\leq \max_{\bw\in\Gamma}\|\nabla F(\bw)\|_2\|\tbw_{1}-\bw^*\|_2\leq G R.
\end{align*}
Combining the above two inequalities, we get the conclusion
\begin{align*}
    \E\left[\sum_{t=1}^T (F(\bw_{t})-F(\bw^*))\right]&\leq    \sum_{\tau=1}^{\tau_{\mathrm{max}}}\E[ n_{\tau+1}(F(\tbw_{\tau+1})-F(\bw^*))] +n_1(F(\tbw_{1})-F(\bw^*))\\&\leq \frac{\lceil\sqrt{T}\rceil R^2}{2\eta}+\frac{\eta\sigma^2\lceil\sqrt{T}\rceil}{2(1-\eta\beta)}+ \lceil\sqrt{T}\rceil GR \\&\leq (\sqrt{T}+1)\left(\frac{R^2}{2\eta}+\frac{\eta\sigma^2}{2(1-\eta\beta)}+  GR\right).
\end{align*}

{
\noindent\underline{\it Proof of the any-time batch scheme.}
Since $n_{\tau}=K\tau$, from the definition of $\tau_{\mathrm{max}}$ we have $K(1+\dots+ (\tau_{\mathrm{max}}-1))\leq T,$
and thus it holds that 
\begin{align*}
    \tau_{\mathrm{max}}\leq 1+\sqrt{\frac{2T}{K}+1}.
\end{align*}
By Eq.~\eqref{eq:descent lemma cvx}, for any $1\leq \tau\leq \tau_{\mathrm{max}}$ we obtain
\[ \E[ n_{\tau+1}(F(\tbw_{\tau+1})-F(\bw^*))] = \frac{K(\tau+1)}{2\eta}\E[ (\|\tbw_{\tau}-\bw^*\|_2^2 -\|\tbw_{\tau+1}-\bw^*\|_2^2)]+\frac{\eta\sigma^2}{(1-\eta\beta)}. \]
Taking summation from $\tau=1$ to $\tau=\tau_{\mathrm{max}}$, we have 
\begin{align*}
    \sum_{\tau=1}^{\tau_{\mathrm{max}}}\E[ n_{\tau+1}(F(\tbw_{\tau+1})-F(\bw^*))]&\leq \sum_{\tau=1}^{\tau_{\mathrm{max}}} \E\left[\frac{2K\|\tbw_{\tau}-\bw^*\|_2^2}{2\eta}\right]+ \frac{\eta\sigma^2\tau_{\mathrm{max}} }{(1-\eta\beta)}
 \\&\leq  
  \frac{KR^2\tau_{\mathrm{max}}}{\eta}+ \frac{\eta\sigma^2\tau_{\mathrm{max}} }{(1-\eta\beta)},
\end{align*}
where the second inequality is due to $\max_{\bw_1,\bw_2\in\Gamma}\|\bw_1-\bw_2\|_2\leq R$.

Since $\max_{\bw\in\Gamma}\|\nabla F(\bw)\|_2\leq G$, by the quasi-mean value theorem  we could have 
\begin{align*}
    F(\tbw_1)-F(\bw^*)\leq \max_{\bw\in\Gamma}\|\nabla F(\bw)\|_2\|\tbw_{1}-\bw^*\|_2\leq G R.
\end{align*}
Combining the above two inequalities, we get the conclusion
\begin{align*}
    \E\left[\sum_{t=1}^T (F(\bw_{t})-F(\bw^*))\right]&\leq    \sum_{\tau=1}^{\tau_{\mathrm{max}}}\E[ n_{\tau+1}(F(\tbw_{\tau+1})-F(\bw^*))] +n_1(F(\tbw_{1})-F(\bw^*))
    \\&\leq \frac{KR^2\tau_{\mathrm{max}}}{\eta}+ \frac{\eta\sigma^2\tau_{\mathrm{max}}}{(1-\eta\beta)}+KGR 
    \\&\leq \left(1+\sqrt{\frac{2T}{K}+1}\right) \left(\frac{KR^2}{\eta}+ \frac{\eta\sigma^2}{(1-\eta\beta)}\right)+KGR, 
\end{align*}
where we complete the proof of this lemma.\Halmos
}
\endproof
\subsection{Proof of Lemma~\ref{le:str}}
\proof{Proof.}
The $\beta$-smooth property of $F(\bw)$ gives
\begin{align}
\nonumber&F(\tbw_{\tau+1})-F(\tbw_{\tau})
\\&\leq \nabla F(\tbw_{\tau})^{\top}\left(\tbw_{\tau+1}-\tbw_{\tau}\right)+\frac{\beta}{2}\left\|\tbw_{\tau+1}-\tbw_{\tau}\right\|_2^2 \nonumber \\ \nonumber
&=g(\tbw_{\tau};\bm{\xi}_{\tau})^{\top}\left(\tbw_{\tau+1}-\tbw_{\tau}\right)+\left(\nabla F(\tbw_{\tau})-g(\tbw_{\tau};\bm{\xi}_{\tau})\right)^{\top}\left(\tbw_{\tau+1}-\tbw_{\tau}\right)+\frac{\beta}{2}\left\|\tbw_{\tau+1}-\tbw_{\tau}\right\|_2^2 \\ 
&\leq g(\tbw_{\tau};\bm{\xi}_{\tau})^{\top}\left(\tbw_{\tau+1}-\tbw_{\tau}\right)+\left\|\nabla F(\tbw_{\tau})-g(\tbw_{\tau};\bm{\xi}_{\tau})\right\|_2\|\tbw_{\tau+1}-\tbw_{\tau}\|_2+\frac{\beta}{2}\left\|\tbw_{\tau+1}-\tbw_{\tau}\right\|_2^2 \label{eq:le:str:smoothness}
\end{align}
where the second inequality is due to Cauchy-Schwarz inequality.

Recalling the updating rule $\tbw_{\tau+1}=\Pi_\Gamma[\tbw_{\tau}-\eta g(\tbw_{\tau};\bm{\xi}_{\tau})]$, by the property of projection onto a convex set, we have 
\begin{equation}\label{eq:le:str:projection}
    (\bw^*-\tbw_{\tau+1})^\top(( \tbw_{\tau}-\eta g(\tbw_{\tau};\bm{\xi}_{\tau}))-\tbw_{\tau+1})\leq 0.
\end{equation}
And it is easy to note that 
\begin{equation}\label{eq:le:str:simple}
    \|\tbw_{\tau}-\bw^*\|_2^2 = \|\tbw_{\tau+1}-\bw^*\|_2^2+\|\tbw_{\tau+1}-\tbw_{\tau}\|_2^2 - 2(\bw^*-\tbw_{\tau+1})^\top(\tbw_{\tau}-\tbw_{\tau+1}).
\end{equation}
Combining \eq{\ref{eq:le:str:projection}} and \eq{\ref{eq:le:str:simple}}, we have 
\begin{align*}
    \eta g(\tbw_{\tau};\bm{\xi}_{\tau})^\top(\tbw_{\tau+1}-\bw^*)&\leq -(\tbw_{\tau+1}-\bw^*)^\top(\tbw_{\tau+1}-\tbw_{\tau})\\&=   \frac{1}{2}  (\|\tbw_{\tau}-\bw^*\|_2^2 -\|\tbw_{\tau+1}-\bw^*\|_2^2-\|\tbw_{\tau+1}-\tbw_{\tau}\|_2^2).
\end{align*}
And thus, we can obtain that
\begin{align*}
    \eta g(\tbw_{\tau};\bm{\xi}_{\tau})^\top(\tbw_{\tau+1}-\tbw_{\tau})&\leq  \eta g(\tbw_{\tau};\bm{\xi}_{\tau})^\top(\bw^*-\tbw_{\tau}) +  \frac{1}{2}  (\|\tbw_{\tau}-\bw^*\|_2^2 -\|\tbw_{\tau+1}-\bw^*\|_2^2-\|\tbw_{\tau+1}-\tbw_{\tau}\|_2^2).
\end{align*}
Combining the above inequality with \eq{\ref{eq:le:str:smoothness}}, we could have
\begin{align}
\nonumber&F(\tbw_{\tau+1})-F(\tbw_{\tau})
\\&\nonumber\leq g(\tbw_{\tau};\bm{\xi}_{\tau})^{\top}\left(\bw^*-\tbw_{\tau}\right)+\frac{1}{2\eta}  (\|\tbw_{\tau}-\bw^*\|_2^2 -\|\tbw_{\tau+1}-\bw^*\|_2^2)\\&\nonumber\quad+\left\|\nabla F(\tbw_{\tau})-g(\tbw_{\tau};\bm{\xi}_{\tau})\right\|_2\|\tbw_{\tau+1}-\tbw_{\tau}\|_2+\frac{\eta\beta-1}{2\eta}\left\|\tbw_{\tau+1}-\tbw_{\tau}\right\|_2^2
\\&\leq g(\tbw_{\tau};\bm{\xi}_{\tau})^{\top}\left(\bw^*-\tbw_{\tau}\right)+\frac{1}{2\eta}  (\|\tbw_{\tau}-\bw^*\|_2^2 -\|\tbw_{\tau+1}-\bw^*\|_2^2)+\frac{\eta\left\|\nabla F(\tbw_{\tau})-g(\tbw_{\tau};\bm{\xi}_{\tau})\right\|^2_2}{2(1-\eta\beta)},\label{eq:le:str:varianceform}
\end{align}
where the last inequality is because $\eta\beta\leq 1$ and $at^2+bt\leq -b^2/4a$ for all $t\in \R$ if $a<0$.

The $\alpha$-strongly convex property of $F(\bw)$ gives
\begin{align}
F(\bw^*)-F(\tbw_{\tau})
\geq \nabla F(\tbw_{\tau})^{\top}\left(\bw^*-\tbw_{\tau}\right)+\frac{\alpha}{2}\left\|\tbw_{\tau+1}-\tbw_{\tau}\right\|_2^2. \label{eq:le:str:strongly}
\end{align}

Combining \eq{\ref{eq:le:str:varianceform}} and \eq{\ref{eq:le:str:strongly}} with the fact that $\E[g(\tbw_{\tau};\bm{\xi}_{\tau})]= \nabla F(\tbw_{\tau})$, we could obtain
\begin{align}
\nonumber
    &\E[F(\tbw_{\tau+1})]
\\&\leq  \E[F(\tbw_{\tau})+\nabla F(\tbw_{\tau})^{\top}\left(\bw^*-\tbw_{\tau}\right)]+\frac{1}{2\eta} \E[ \|\tbw_{\tau}-\bw^*\|_2^2 -\|\tbw_{\tau+1}-\bw^*\|_2^2]+\frac{\eta\E[\left\|\nabla F(\tbw_{\tau})-g(\tbw_{\tau};\bm{\xi}_{\tau})\right\|^2_2]}{2(1-\eta\beta)}
\nonumber
\\&\leq F(\bw^*)-\frac{\alpha}{2}\E[\|\tbw_{\tau}-\bw^*\|_2^2] +\frac{1}{2\eta}  \E[\|\tbw_{\tau}-\bw^*\|_2^2 -\|\tbw_{\tau+1}-\bw^*\|_2^2]+\frac{\eta\E[\|\nabla F(\tbw_{\tau})-g(\tbw_{\tau};\bm{\xi}_{\tau})\|^2_2]}{2(1-\eta\beta)}.\label{eq:le:str:value}
\end{align}
Noting $F(\tbw_{\tau+1})\geq F(\bw^*)$ in the above inequality, we obtain the following \textit{descent lemma}.
\begin{align}
  \nonumber  \E[\|\tbw_{\tau+1}-\bw^*\|_2^2]&\leq (1-\eta\alpha )\E[\|\tbw_{\tau}-\bw^*\|_2^2 ]+\frac{\eta^2\E[\left\|\nabla F(\tbw_{\tau})-g(\tbw_{\tau};\bm{\xi}_{\tau})\right\|^2_2]}{(1-\eta\beta)}\\& \leq(1-\eta\alpha )\E[\|\tbw_{\tau}-\bw^*\|_2^2 ]+\frac{\eta^2\sigma^2}{n_{\tau}(1-\eta\beta)},\label{eq:le:str:regret}
\end{align}
where the second inequality is because
\begin{align}
\nonumber
    \E[\left\|\nabla F(\tbw_{\tau})-g(\tbw_{\tau};\bm{\xi}_{\tau})\right\|_2^2]
    &= \E\left[\left\|\nabla F(\tbw_{\tau})-\frac{1}{n_{\tau}} \sum_{i=1}^{n_{\tau}} \nabla f\left(\tbw_{\tau} ; \xi_{i,\tau}\right)\right\|_2^2\right]\\
\nonumber
    &= \frac{1}{n_{\tau}^2}\sum_{i=1}^{n_{\tau}} \E[\left\|\nabla F(\tbw_{\tau})-\nabla f\left(\tbw_{\tau} ; \xi_{i,\tau}\right)\right\|_2^2] \\&\leq \frac{\sigma^2}{n_{\tau}}.\label{eq:le:str:varreduct}
\end{align}
With $\eta\leq \min\{\alpha/2,1/\alpha,1/(2\beta)\}$ and $n_{\tau}=\lceil\varsigma^{\tau-1}\rceil$ where $\varsigma= 1/{\gamma}$ and $\gamma=1-\eta\alpha+2\eta^2\in(0,1)$, we can prove $\mathbb{E}\left[\left\|\tbw_{\tau}-\bw^*\right\|_2^2\right] 
\leq\kappa \gamma^{\tau-1}$ for all $\tau\geq 1$ by induction, where $\kappa = \max\{R^2,\sigma^2\}$.

For $\tau=1$, by the definition of $\kappa$ it is easy to obtain that $\mathbb{E}\left[\left\|\bw_{1}-\bw^*\right\|_2^2\right] 
\leq\kappa$. Assuming that  $\mathbb{E}\left[\left\|\tbw_{\tau}-\bw^*\right\|_2^2\right] 
\leq\kappa \gamma^{\tau-1}$ for $\tau\geq 1$, by \eq{\ref{eq:le:str:regret}} we have
\begin{align*}
\mathbb{E}\left[\left\|\tbw_{\tau+1}-\bw^*\right\|_2^2\right] 
&\leq\left(1- \eta  \alpha\right) \mathbb{E}\left[\left\|\tbw_{\tau}-\bw^*\right\|_2^2\right]+\frac{\eta ^2 \sigma^2\varsigma^{1-\tau}}{1-\eta\beta} \\& \leq \left(1- \eta  \alpha\right)\kappa \gamma^{\tau-1}+\frac{\eta ^2 \sigma^2\varsigma^{1-\tau}}{1-\eta\beta} 
\\& \leq \kappa \gamma^{\tau-1}\left(1- \eta  \alpha+\frac{\eta ^2 \sigma^2(\varsigma\gamma)^{1-\tau}}{\kappa(1-\eta\beta)}\right)
\\& =\kappa \gamma^{\tau-1}\left(1- \eta  \alpha+\frac{\eta ^2 \sigma^2}{\kappa(1-\eta\beta)}\right)\\& \leq \kappa \gamma^{\tau-1}\left(1- \eta  \alpha+\frac{2\eta ^2 \sigma^2}{\kappa}\right)
\\& \leq \kappa \gamma^{\tau},
\end{align*}
where the first inequality is due to $n_{\tau}=\lceil\varsigma^{\tau-1}\rceil\geq \varsigma^{\tau-1}$,  the second inequality is due to the induction assumption, the fourth inequality is due to $\varsigma\gamma=1$, the fifth inequality is due to $\eta\leq 1/(2\beta)$, and the last inequality is due to $\kappa \geq \sigma^2$.

Combining above inequality and \eq{\ref{eq:le:str:value}}, we know
\begin{align}
\nonumber
    &\E[F(\tbw_{\tau+1})]
\\
\nonumber
&\leq F(\bw^*)-\frac{\alpha}{2}\E[\|\tbw_{\tau}-\bw^*\|_2^2] +\frac{1}{2\eta}  \E[\|\tbw_{\tau}-\bw^*\|_2^2 -\|\tbw_{\tau+1}-\bw^*\|_2^2]+\frac{\eta\E[\|\nabla F(\tbw_{\tau})-g(\tbw_{\tau};\bm{\xi}_{\tau})\|^2_2]}{2(1-\eta\beta)} \\
\nonumber
&\leq F(\bw^*)+\frac{1}{2\eta}  \E[\|\tbw_{\tau}-\bw^*\|_2^2] +\frac{\eta\E[\|\nabla F(\tbw_{\tau})-g(\tbw_{\tau};\bm{\xi}_{\tau})\|^2_2]}{2(1-\eta\beta)} \\
&\leq F(\bw^*) + \frac{\kappa \gamma^{\tau-1}}{2\eta} + \frac{\eta^2\sigma^2}{n_{\tau}(1-\eta\beta)}, \label{eq:le:str:value2}
\end{align}
where the last inequality is by \eq{\ref{eq:le:str:varreduct}}.
Since $\max_{\bw\in\Gamma}\|\nabla F(\bw)\|_2\leq G$, by the quasi-mean value theorem  we could have 
\begin{align}
    F(\tbw_{\tau})-F(\bw^*)\leq \max_{\bw\in\Gamma}\|\nabla F(\bw)\|_2\|\tbw_{\tau}-\bw^*\|_2\leq G \|\tbw_{\tau}-\bw^*\|_2 \leq GR.
\end{align}
From the definition of $\tau_{\mathrm{max}}$, it holds that $\tau_{\mathrm{max}}\leq (\ln{((\varsigma-1)T+1)}/\ln{\varsigma})+2)$. Combining this fact with \eq{\ref{eq:le:str:value2}}, we get the conclusion
\begin{align*}
\E\left[\sum_{t=1}^T (F(\bw_{t})-F(\bw^*))\right]&\leq  \E[\sum_{\tau=1}^{\tau_{\mathrm{max}}-1} n_{\tau+1} [F(\tbw_{\tau+1})-F(\bw^*)]+\E[n_1[F(\tbw_{1})-F(\bw^*)]]
\\&\leq \sum_{\tau=1}^{\tau_{\mathrm{max}-1}} \left(\frac{\kappa n_{\tau+1} \gamma^{\tau-1}}{2\eta} + \frac{n_{\tau+1}\eta^2\sigma^2}{n_{\tau}(1-\eta\beta)}\right) + GR
    \\&\leq \sum_{\tau=1}^{\tau_{\mathrm{max}-1}} \left(\frac{\kappa \varsigma^2 (\varsigma\gamma)^{\tau-1}}{\eta}+\frac{\varsigma^2\eta^2 \sigma^2}{(1-\eta \beta)}\right) + GR \\&\leq GR +\left(\frac{\kappa\varsigma^2}{\eta}+\frac{\varsigma^2\eta^2 \sigma^2}{(1-\eta \beta)}\right)\left(\frac{\ln{((\varsigma-1)T+1)}}{\ln{\varsigma}}+2\right),
\end{align*}
where the third inequality is due to $n_{\tau+1}=\lceil\varsigma^{\tau}\rceil\leq \varsigma^{\tau+1}$ and $n_{\tau}=\lceil\varsigma^{\tau-1}\rceil\geq \varsigma^{\tau-1}$.  \Halmos

\endproof

\subsection{Proof of Theorem~\ref{thm:inv,general}} \label{sec:proof-of-thm-inv-general}
\proof{Proof of Theorem~\ref{thm:inv,general}.}
We first partition the index set $I = \{1,2,\dots,T \}$ into two index sets:
\begin{align*}
I_1 = \{t \mid x_{i,t} > w_{i,t}, \exists i \in [n] \}
\quad \text{and} \quad 
I_2 = \{t \mid x_{i,t} \leq w_{i,t}, \forall i \in [n]  \}.
\end{align*}
As mentioned before, we call the period in $I_1$ a waiting period and the period in $I_2$ a working period. Then we express the regret as follows:
\begin{align*}
\mathcal R_{\mathrm{cvx}}(T) &= \E \left[ \sum_{t=1}^T (Q(\by_t)-Q(\by^*))\right] = \Lambda_1(T) + \Lambda_2(T),
\end{align*}
where
\begin{align*}
\Lambda_1(T) = \mathbb{E} \left[ \sum_{t \in I_1} (Q(\by_t)-Q(\by^*)) \right]
\quad \text{and} \quad 
\Lambda_2(T) = \mathbb{E} \left[\sum_{t \in I_2} (Q(\by_t)-Q(\by^*)) \right].
\end{align*}
Obviously, $\Lambda_1(T)$ is the regret incurred in all waiting periods, which is at most of the same order as the number of the waiting periods $|I_1|$.

By the procedures of Algorithm~\ref{alg1}, when a new target inventory level $\bw$ updated by minibatch SGD is infeasible, the algorithm calls {\sc Transition Solver} and enters into a consecutive epoch of waiting periods. 
By the definition of {\sc Transition Solver} (Definition~\ref{def:transition solver}), the expected waiting periods due to this switch are at most $M$.

Combining this with the fact that the switch time of our meta-policy for the fixed-time setting,  is bounded by $(\sqrt{T}+1)$, it holds that
$$|I_1|\leq M(\sqrt{T}+1).$$

Moreover, since $0 \leq Q(\by) \leq \bar{Q}$, for the fixed-time batch scheme  we have
\[
\Lambda_1(T) = \mathbb{E} \left[\sum_{t \in I_1} (Q(\by_t)-Q(\by^*)) \right]
\leq \bar{Q} \E \left[ I_1 \right] \leq \bar{Q}M(\sqrt{T}+1).
\]

$\Lambda_2(T)$ is regret incurred in all working periods, which is highly related to the cumulative regret analysis of minibatch SGD. Therefore, in the following, we bound $\Lambda_2(T)$ by relating it with the cumulative regret of minibatch SGD established in Lemma~\ref{le:cvx}.

It is easy to note that all the conditions of Lemma \ref{le:cvx} are satisfied in our problem. And we could treat the inventory level $\by_t$ in all periods of $I_2$ as a subsequence of $\bw_t$ generated by the minibatch SGD algorithm when used to minimize $Q(\by)$ over $\Gamma$ in $T$ periods.

With the above analysis, 
by Lemma \ref{le:cvx} we could upper bound $\Lambda_2(T)$ as
\[
\Lambda_2(T) \leq (\sqrt{T}+1)\left(\frac{R^2}{2\eta}+\frac{\eta\sigma^2}{2(1-\eta\beta)}+  GR\right).
\]
Combining the upper bounds of $\Lambda_1(T)$ and $\Lambda_2(T)$, for the fixed-time batch scheme we conclude that
\begin{align*}
    \mathcal R_{\mathrm{cvx}}(T) &= \Lambda_1(T) + \Lambda_2(T) \leq (\sqrt{T}+1)\left(\bar{Q}M+\frac{R^2}{2\eta}+\frac{\eta\sigma^2}{2(1-\eta\beta)}+  GR\right).
\end{align*}
Similarly, for the any-time batch scheme we have 
\begin{align*}
    \mathcal R_{\mathrm{cvx}}(T) &= \Lambda_1(T) + \Lambda_2(T) \\&\leq \left(1+\sqrt{\frac{2T}{K}+1}\right) \bar{Q}M+   \left(1+\sqrt{\frac{2T}{K}+1}\right) \left(\frac{KR^2}{\eta}+ \frac{\eta\sigma^2}{(1-\eta\beta)}\right)+KGR \\&\leq
 \left(1+\sqrt{\frac{2T}{K}+1}\right) \left(\bar{Q}M+\frac{KR^2}{\eta}+ \frac{\eta\sigma^2}{(1-\eta\beta)}\right)+KGR,
\end{align*}
which completes the proof of the theorem.
\Halmos
\endproof
\subsection{Proof of Theorem~\ref{thm:inv,str}}\label{sec:proof-of-thm-inv-str}
\proof{Proof of Theorem~\ref{thm:inv,str}.}
We first partition the index set $I = \{1,2,\dots,T \}$ into two index sets:
\begin{align*}
I_1 = \{t \mid x_{i,t} > w_{i,t}, \exists i \in [n] \}
\quad \text{and} \quad 
I_2 = \{t \mid x_{i,t} \leq w_{i,t}, \forall i \in [n]  \}.
\end{align*}
As mentioned before, we call the period in $I_1$ a waiting period and the period in $I_2$ a working period. Then we express regret as follows:
\begin{align*}
\mathcal R_{\mathrm{str}}(T) &= \E \left[ \sum_{t=1}^T (Q(\by_t)-Q(\by^*))\right]
\\
&= \Lambda_1(T) + \Lambda_2(T),
\end{align*}
where
\begin{align*}
\Lambda_1(T) = \mathbb{E} \left[ \sum_{t \in I_1} (Q(\by_t)-Q(\by^*)) \right]
\quad \text{and} \quad 
\Lambda_2(T) = \mathbb{E} \left[\sum_{t \in I_2} (Q(\by_t)-Q(\by^*)) \right]
\end{align*} 
Obviously, $\Lambda_1(T)$ is regret incurred in all waiting periods, which is of the same order as the number of the waiting periods $|I_1|$.

By the procedures of Algorithm~\ref{alg1}, when a new target inventory level $\bw$ updated by minibatch SGD is infeasible, the algorithm calls {\sc Transition Solver} and enters into a consecutive epoch of waiting periods. 
By the definition of {\sc Transition Solver} (Definition~\ref{def:transition solver}), the expected waiting periods due to this switch are at most $M$.

{Since total number of switches is at most $\tau_{\mathrm{max}}-1$ times, where  $\tau_{\mathrm{max}} = \min\{k:\sum_{\tau = 1}^k n_{\tau} \geq T\}$, it holds that $$n_{\tau_{\mathrm{max}}} = \lceil \varepsilon^{\tau_{\mathrm{max}}-1}\rceil\leq T,$$
and thus
\[|I_1| \leq \tau_{\mathrm{max}}-1\leq  M\frac{\ln T}{\ln \varsigma}. \]
}

Moreover, since $0 \leq Q(\by) \leq \bar{Q}$, we have
\[
\Lambda_1(T) = \mathbb{E} \left[\sum_{t \in I_1} (Q(\by_t)-Q(\by^*)) \right]
\leq \bar{Q} \E \left[ I_1 \right] \leq \bar{Q}M(\frac{\ln T}{\ln \varsigma}).
\]

$\Lambda_2(T)$ is regret incurred in all working periods, which is highly related to the cumulative regret analysis of minibatch SGD. Therefore, in the following, we bound $\Lambda_2(T)$ by relating it  with the cumulative regret  of minibatch SGD established in Lemma~\ref{le:str}.

It is easy to note that all the conditions of Lemma \ref{le:str} are satisfied in our problem. And we could treat the inventory level $\by_t$ in all periods of $I_2$ as a subsequence of $\bw_t$ generated by the minibatch SGD algorithm when used to minimize $Q(\by)$ over $\Gamma$ in $T$ periods,

With the above analysis, 
by Lemma \ref{le:str} we could upper bound $\Lambda_2(T)$ as
\[
\Lambda_2(T)  \leq GR +\left(\frac{\kappa\varsigma^2}{\eta}+\frac{\varsigma^2\eta^2 \sigma^2}{(1-\eta \beta)}\right)\left(\frac{\ln{((\varsigma-1)T+1)}}{\ln{\varsigma}}+2\right).
\]
Combining the upper bounds of $\Lambda_1(T)$ and $\Lambda_2(T)$, we could get the conclusion,
\begin{align*}
    \mathcal R_{\mathrm{str}}(T) &= \Lambda_1(T) + \Lambda_2(T) \\
    &\leq \bar{Q}M\frac{\ln T}{\ln \varsigma}+GR +\left(\frac{\kappa\varsigma^2}{\eta}+\frac{\varsigma^2\eta^2 \sigma^2}{(1-\eta \beta)}\right)\left(\frac{\ln{((\varsigma-1)T+1)}}{\ln{\varsigma}}+2\right),
\end{align*}
which completes the proof of this theorem.
\Halmos
\endproof
\section{Proofs and Discussions Omitted in Section~\ref{app: multi-product}}
\subsection{Detailed Discussions about \cite{shi2016nonparametric}} \label{subsec:shicong}
{In this section, we discuss the content and roles of some lemmas in the regret analysis of \cite{shi2016nonparametric} and explain why their proof techniques fail in the setting of Application I. For readers' convenience, we keep all notations the same as the ones in \cite{shi2016nonparametric}. 

In their paper, $\Pi(\by)$, the expected cost function when the actually implemented inventory level is $\by$, is defined as
\[
\Pi(\by) = \E[\bm{c}\cdot \by + (\bh-\bm{c})\cdot (\by-\bD)^+ + \bp\cdot(\bD-\by)^+],
\]
where $\bx \cdot \by$ is the inner product of $\bx$ and $\by$, $\bm{c} = (c^1,\dots,c^n)$, $\bm{h} = (h^1,\dots,h^n)$, and $\bm{p} = (p^1,\dots,p^n)$ is cost parameters and $\bD = (D^1,\dots,D^n)$ is the random demand vector.}

{
The Data-Driven Method (DDM) algorithm in \cite{shi2016nonparametric} generates two sequences: $(\hat{\by})_{t\geq 0}$ is the target inventory level updated by SGD and $({\by})_{t\geq 0}$ is the actually implemented inventory level due to inventory carry-over. Then the total regret incurred by DDM is
\[
\mathcal{R}(T) := \E\left[\sum_{t=1}^T \Pi(\by_t) \right] - \sum_{t=1}^T \Pi(\by^*),
\]
where $\by^*$ is the minima of $\Pi(\by)$ under the single warehouse constraint.}

{
First, \cite{shi2016nonparametric} decomposed the regret into two parts as
\[
\mathcal{R}(T) = \E\left[\sum_{t=1}^T \Pi(\by_t) -\sum_{t=1}^T \Pi(\hat{\by}_t)\right] + \E\left[\sum_{t=1}^T \Pi(\hat{\by}_t)- \sum_{t=1}^T \Pi(\by^*). \right]
\]
Recall $(\hat{\by})_{t\geq 0}$ are the target inventory levels, and the second term
\[
\E\left[\sum_{t=1}^T \Pi(\hat{\by}_t)- \sum_{t=1}^T \Pi(\by^*). \right]
\]
can be bound based on the ideas and techniques used in standard SGD.}

{
By the Lipschitz continuity of $\Pi(\by)$,
\[
\E\left[\sum_{t=1}^T \Pi(\by_t) -\sum_{t=1}^T \Pi(\hat{\by}_t)\right] \leq C_1 \cdot \E\left[\sum_{t=1}^T \sum_{i=1}^n |y_t^i - \hat{y}_t^i|\right],
\]
where $C_1$ is a universal constant, $n$ is the dimension of $\by$, and $y_t^i$, $\hat{y}_t^i$ is the $i$-th component of $\by_t$ and $\hat{\by}_t$, respectively. Next, \cite{shi2016nonparametric} proposed a queueing-theory-based method to bound the term on the above of the right-hand side. }

{
When there is inventory carry-over, \cite{shi2016nonparametric} analyzed the behavior of the inventory system under the warehouse constraint and proved some intermediate lemmas.  For brevity, we list one of them as follows.  \cite{shi2016nonparametric} frequently used the tight warehouse constraint $\sum_{i=1}^n y_t^i = M$ to prove this lemma. However, such proof techniques can not work under general multiple constraints.}
{
\begin{lemma}[Lemma 5 of \cite{shi2016nonparametric}]
Let $x_{t+1}^i$ be the on-hand inventory of product i at the beginning of period $t+1$.
In each period $t+1$, we bound the distance function for all $i \in \bar{J} := \left\{i: x_{t+1}^i \leq \hat{y}_{t+1}^i\right\}$ as follows. If $J=\emptyset$, we have $\sum_{i \in \bar{J}}\left|\hat{y}_{t+1}^i-y_{t+1}^i\right|=0$. Otherwise, if $J \neq \emptyset$, we have
$$
\sum_{i \in \bar{J}}\left|\hat{y}_{t+1}^i-y_{t+1}^i\right| \leq \sum_{i \in \bar{J}}\left|\hat{y}_t^i-y_t^i\right|+\eta_t \sum_{i \in \bar{J}}\left(p^i-c^i\right)-\sum_{j \in J} d_t^j.
$$
\end{lemma}

Based on the above lemma, \cite{shi2016nonparametric} derived the following inequality.}

{

\begin{lemma}[Lemma 6 of \cite{shi2016nonparametric}]
In each period $t+1$, we bound the sum of distance functions as follows.
$$
\sum_{i=1}^n\left|y_{t+1}^i-\hat{y}_{t+1}^i\right| \leq\left(\sum_{i=1}^n\left|y_t^i-\hat{y}_t^i\right|+\eta_t\left(\sum_{i=1}^n\left(h^i+2\left(p^i-c^i\right)\right)\right)-\min _{j=1, \ldots, n} d_t^j\right)^{+},
$$
where $\eta_t$ is the stepsize of DDM and $d_t^j$ is a realization of $D_t^j$.
\end{lemma}
Motivated by the above inequality, \cite{shi2016nonparametric} defined a stochastic process $(Z_t \mid t \geq 0)$ as
\[
Z_{t+1} = \left[Z_t + \frac{S_t}{\sqrt{t}} - \tilde{D}_t \right]^+,\quad Z_0 = 0,
\]
where $S_t = \sum_{i=1}^n (h^i + 2(p^i-c^i))$ and $\tilde{D}_t$ is a random variable satisfying $\tilde{D}_t \leq D_t^j,\forall j$ in stochastic order. Based on stochastic comparison, they proved the following lemma.
\begin{lemma}[Lemma 7 of \cite{shi2016nonparametric}]
The total expected distance function
$$
\mathbb{E}\left[\sum_{t=1}^T \sum_{i=1}^n\left|y_t^i-\hat{y}_t^i\right|\right] \leq \mathbb{E}\left[\sum_{t=1}^T Z_t\right]
$$
where $Z_{t+1}^i$ is a stochastic process defined above.
\end{lemma}}

{
With the above Lemma, it suffices to bound the expectations of the process $(Z_t \mid t \geq 0)$. Note that $(Z_t \mid t \geq 0)$ is very similar to a GI/GI/1 queue and they considered a GI/G/1 process as
\[
W_{t+1} = [W_t + S_t-\tilde{D}_t]^+,\quad,W_0=0.
\]
Let $B$ be the busy period of the above queue.
They proved the following lemma.
\begin{lemma}[Lemma 8 of \cite{shi2016nonparametric}]
The total expected distance function
$$
\mathbb{E}\left[\sum_{t=1}^T \sum_{i=1}^n\left|y_t^i-\hat{y}_t^i\right|\right] \leq C'\E[B] \sqrt{T},
$$
where $C'$ is a constant depending on the queue.
\end{lemma}}{
Therefore, they establish an $\O(\sqrt{T})$ bound for regret provided that $\E[B] < \infty$.
Unfortunately, if the demand is dependent among products, the definition of $\tilde{D}_t$ such that $\tilde{D}_t \leq D_t^j,\forall j$ may cause $\tilde{D}_t = 0$. Therefore, the inter-arrival time $\tilde{D}_t$ of the queue may be zero and the busy period can be infinite. The bound of Lemma 8 in \cite{shi2016nonparametric} no longer holds and their proof fails.

{
In our numerical experiments, the SGD-based algorithms also have a good regret performance, but require longer running time. It seems quite difficult to establish the theoretical guarantee for SGD-based algorithms for complex inventory management problems.
The state-of-the-art queueing-theory-based analytic techniques in~\cite{huh2009nonparametric,shi2016nonparametric} are complex in nature and the analyses in the simple settings~\citep{huh2009nonparametric,shi2016nonparametric} already require quite a lot technical effort.
This complexity hinders the further application of the techniques to even more complicated inventory systems. 
Given the limitation of the state-of-the-art analytic techniques for the SGD-based inventory algorithms, we believe that new analysis techniques are needed for learning complex inventory systems. }

{\subsection{Upper bound of $R$}
\begin{lemma}\label{le:Upper bound of R}
   Note that $R = \max_{\by,\by' \in \Gamma} \Vert \by-\by' \Vert_2$ and $\Gamma = \{\by \in \mathbb{R}^n_+ \mid A \by \leq \bm{\rho} \}$,  we have \[R\leq  2\sqrt{\sum_{j=1}^n(\min_{i\in[m]}\rho_j/A_{i,j})^2}.\]
\end{lemma}
\proof{Proof.}
By the nonnegativity of elements of $A \in \mathbb{R}^{m\times n}_+$ and the definition of  $\Gamma = \{\by \in \mathbb{R}^n_+ \mid A \by \leq \bm{\rho} \}$, for each$j\in[n]$ we have $$\bm A_jy_j\leq \bm\rho,$$ where $\bm A_j$ is the $j$-th column of matrix. Thus  we have \[y_j\leq \min_{i\in[m]}\rho_i/A_{i,j}.\] Consequently, \[R\leq 2\max_{\by\in\Gamma}\|\by\|_2\leq 2\sqrt{\sum_{j=1}^n(\min_{i\in[m]}\rho_j/A_{i,j})^2},\]
where we prove this lemma.\Halmos
\endproof
}}

\subsection{Proof of Lemma~\ref{le:app1properties}}
We only validate the first and second properties, since the other properties are easy to validate.
\proof{Proof.}
Since the expected one-period  cost when the after-delivery inventory level is $\by$ as $$Q(\by) = \bh^\top \E[(\by - \bD)^+] + \bb ^\top \E[(\bD-\by)^+],$$
and the  gradient of $Q(\by)$ at $\by$ is given by $$\nabla Q(\by) = \E \left[ \bh\odot\mathbb{I}[\by>\bD]-\bb\odot\mathbb{I}[\by\leq\bD] \right],$$
we could get the Hessian Matrix at $\by$, which is a diagonal matrix given as follows
$$\nabla^2 Q(\by) = \mathrm{Diag}\{(h_1+b_1)f_{D_1}(y_1),\dots,(h_n+b_n)f_{D_n}(y_n)\}.$$
Under Assumption \ref{ass1:generalass}, we have
$$0\preceq\nabla^2 Q(\by) \preceq (\max_{i\in[n]} (h_i+b_i) \beta_0)\mathbf{I_n}.$$
Therefore, $Q(\by)$ is convex and $(\max_{i\in[n]}(h_i+b_i)\beta_0)$-smooth under Assumption \ref{ass1:generalass}.

Under Assumption \ref{ass1:generalass} and \ref{ass2:strongly}, we have
$$(\min_{i\in[n]}(h_i+b_i)\alpha_0)\mathbf{I_n}\preceq\nabla^2 Q(\by) \preceq (\max_{i\in[n]} (h_i+b_i) \beta_0)\mathbf{I_n}.$$
Therefore, $Q(\by)$ is $(\min_{i\in[n]}(h_i+b_i)\alpha_0)$-strongly convex and $(\max_{i\in[n]}(h_i+b_i)\beta_0)$-smooth under Assumption \ref{ass1:generalass} and \ref{ass2:strongly}.
\Halmos
\endproof
\subsection{Proof of Lemma~\ref{app1:le:transition solver}}
\begin{lemma}[Lipschitz Hitting Time of Ascending Random Walk \citep{yuan2021marrying}]\label{le:yuanhao}
 Suppose $\{D_t\}_{t=1}^{\infty}$ is a sequence of i.i.d.~positive random variables with densities bounded by a constant $\rho$. Define an ascending random walk $W_0=0$ and $W_t=\sum_{i=1}^t D_i$ for $t=1,2, \ldots$. For any $\delta>0$, let $L(\delta)$ be the hitting time to the interval $[\delta, \infty)$, i.e., $L(\delta)=\min \{t: W_t \geq \delta\}$. Then, $\mathbb{E}[L(\delta)]$ is Lipschitz in $\delta$, with the Lipschitz constant being $6 \rho$.
\end{lemma}
\begin{lemma}[Bound on Hitting Time of Multi-dimensional Ascending Random Walk]
\label{lem:hittingtime}
Suppose $\{\bD_t = (D_{1,t},D_{2,t},\dots,D_{n,t})^{\top}\}_{t=1}^\infty$ is a sequence of i.i.d.~positive  random vectors with a density of each component bounded by a constant $\beta > 0$. Define an ascending random walk $\bm{W}_0= \bm{0}$, and $\bm{W}_t = \sum_{i=1}^t \bD_{i}$. For any $\bm{\delta}=(\delta_1,\delta_2,\dots,\delta_n)^{\top}$ with $\delta_i>0$, let $L(\bm{\delta})$ be the hitting time to the set $[\bm{\delta},+\infty )$, i.e., $L(\bm{\delta}) = \min \{t \mid \bm{W}_t \geq \bm{\delta} \}$.
Then $\E[L(\bm{\delta})]$ is bounded by $n + 6\beta \sum_{i=1}^n \delta_i$.
\end{lemma}
\proof{Proof.}
From the definition of $L(\bm{\delta})$, we know that
\begin{align*}
L(\bm{\delta}) &= \min \{t \mid \bm{W}_t \geq \bm{\delta} \} \\
&= \min \{t \mid {W}_{1,t} \geq \delta_1,{W}_{2,t} \geq \delta_2  ,\dots,{W}_{n,t} \geq \delta_n \}
\\
&\leq \min \{t \mid {W}_{1,t} \geq \delta_1\} + \dots + \min \{t \mid {W}_{n,t} \geq \delta_n\}.
\end{align*}
The second equation is by expressing $\bm{W}_t \leq \bm{\delta}$ as components. Then we explain why the last inequality holds. Suppose that $t'$ is the first time $t$ such that $\bm{W}_t \geq \bm{\delta}$. We must have $W_{i,t'} \geq \delta_i$ for all $i \in [n]$ and $W_{i,t'} < \delta_i$ for some $i \in [n]$. So the summation of all one-dimensional hitting time contains an item $t'$, which gives that
\[
t' = L(\bm{\delta}) \leq \min \{t \mid {W}_{1,t} \geq \delta_1,{W}_{2,t} \geq \delta_2  ,\dots,{W}_{n,t} \geq \delta_n \}.
\]
On the other hand, by  Lemma~\ref{le:yuanhao}, we know that $\min \{t \mid {W}_{i,t} \geq \delta_i \} \leq 1+6\beta \delta_i$. Summation from $i\in[n]$ gives the conclusion
\begin{align*}
\E[L(\bm{\delta})] 
&\leq \E[\min \{t \mid {W}_{1,t} \geq \delta_1\} + \dots + \min \{t \mid {W}_{n,t} \geq \delta_n\}]\\
&\leq n + 6\beta \sum_{i=1}^n \delta_i.
\end{align*}
\Halmos
\endproof
\proof{Proof of Lemma~\ref{app1:le:transition solver}.}
Consider the worst case that target inventory levels of all products are infeasible. For each product $i$, the capacity for the product is upper bounded by the diameter $R$ of $\Gamma$. 
Therefore, by Lemma \ref{lem:hittingtime}, we can bound $\E[s(\bx_0,\bw)]$ as follows
$$\E [s(\bx_0,\bw)] \leq n + 6n\beta_0 R,$$
where we complete the proof of this lemma.
\Halmos
\endproof

\section{Proofs and Explanations Omitted in Section~\ref{app: owms}} \label{appendix:app3}

\subsection{Transform the Initial Problem into Linear Program}\label{subsec: proofoflemmaapp3}
In this section, we transform Problem~\eqref{eq:delivery} into a linear program. The convexity of $Q(\by)$ is proved as a by-product. Recall that $Q(\by) = \E[C(\by,D)]$ and $C(\by,\bd)$ is defined by the following second-stage optimization problem
\begin{equation}\label{eq:twostage}
\begin{aligned} 
C(\by,\bm{d}) := &\min \bm{c}^\top \bm{z} + h_0(y_{i}-\sum_{i=1}^n z_{i})^+ + \sum_{i=1}^n [h_i(z_i+y_i-d_i)^+ + b_i(d_i-z_i-y_i)^+]\\
&\begin{array}{r@{\quad}r@{}l@{\quad}l}
s.t. 
     &\sum_{i=1}^n z_i &\leq y_0,\\
     &z_i &\geq 0, i \in [n].\\
\end{array}
\end{aligned}
\end{equation}
Introducing variables $z_0^{(1)},z_i^{(1)},z_i^{(2)},i \in [n]$ satisfying
\begin{equation}\nonumber
\begin{aligned}
\left\{
\begin{array}{cl}
    y_0-z_1-\dots-z_n \ &\ \leq  z_0^{(1)},\\
    z_i + y_i-d_i \ &\ \leq z_i^{(1)}, i \in [n] \\
      d_i-z_i-y_i \ &\ \leq z_i^{(2)},i \in [n],
\end{array}
    \right.
\end{aligned}
\end{equation}
we can transform Problem (\ref{eq:twostage}) into the following linear program.
\begin{equation}\nonumber
\begin{aligned} 
&\min \bm{c}^\top \bm{z} + \sum_{i=0}^n h_i z_i^{(1)} + \sum_{i=1}^n b_iz_i^{(2)}\\
&\begin{array}{r@{\quad}r@{}l@{\quad}l}
s.t. 
     &\sum_{i=1}^n z_i &\leq y_0,\\
     &y_0-\sum_{i=1}^n z_i &\leq z_0^{(1)}, \\
    &z_i + y_i-d_i &\leq z_i^{(1)},i \in [n], \\
     &d_i - z_i-y_i &\leq z_i^{(2)},i \in [n],\\
     &z_i,z_0^{(1)},z_i^{(1)},z_i^{(2)}  &\geq 0,i \in [n].\\
\end{array}
\end{aligned}
\end{equation}
Next, we add some slack variables $z_0^{(3)},z_0^{(4)},z_i^{(4)},z_i^{(5)},i \in [n]$ to the inequality constraints and change them into equality constraints.
\begin{equation}\nonumber
\begin{aligned} 
&\min \bm{c}^\top \bm{z} + \sum_{i=0}^n h_i z_i^{(1)} + \sum_{i=1}^n b_iz_i^{(2)}\\
&\begin{array}{r@{\quad}r@{}l@{\quad}l}
s.t. 
     &y_0 &= \sum_{i=1}^n z_i +z_0^{(3)},\\
     &z_0^{(1)} &= y_0-\sum_{i=1}^n z_i + z_0^{(4)}, \\
    &z_i^{(1)} &= z_i + y_i-d_i+z_i^{(4)},i \in [n], \\
     &z_i^{(2)} &= d_i - z_i-y_i+z_i^{(5)},i \in [n], \\
     &z_i,z_0^{(1)},z_i^{(1)},z_i^{(2)}  &\geq 0,i \in [n],\\
&z_0^{(3)},z_0^{(4)},z_i^{(4)},z_i^{(5)}  &\geq 0, i \in [n].\\
\end{array}
\end{aligned}
\end{equation}
or equivalently
\begin{equation}
\label{eq:twostageLP}
\begin{aligned} 
&\min \bm{c}^\top \bm{z} + \sum_{i=0}^n h_i z_i^{(1)} + \sum_{i=1}^n b_iz_i^{(2)}\\
&\begin{array}{r@{\quad}r@{}l@{\quad}l}
s.t. 
     &y_0-\sum_{i=1}^n z_i -z_0^{(3)} &= 0,\\
     &y_0-\sum_{i=1}^n z_i - z_0^{(1)}+z_0^{(4)} &= 0, \\
    &y_i+z_i-z_i^{(1)}+z_i^{(4)} &=d_i,i \in [n], \\
     &y_i+z_i+z_i^{(2)}-z_i^{(5)} &=d_i,i \in [n], \\
     &z_i,z_0^{(1)},z_i^{(1)},z_i^{(2)}  &\geq 0,i \in [n],\\
&z_0^{(3)},z_0^{(4)},z_i^{(4)},z_i^{(5)}  &\geq 0, i \in [n].\\
\end{array}
\end{aligned}
\end{equation}
Let $$\bm{z}' =(z_1,\dots,z_n,z_0^{(1)},\dots,z_n^{(1)},z_1^{(2)},\dots,z_n^{(2)},z_0^{(3)},z_0^{(4)},\dots,z_n^{(4)},z_0^{(5)},\dots,z_n^{(5)}), $$ and $$\bd'= (0,0,d_1,\dots,d_n,d_1,\dots,d_n)^{\top}.$$
Let $\bm{c}'$ be the coefficient of $\bm z'$ in the objective function of Problem (\ref{eq:twostageLP}). And we denote by $T, W$ the coefficient matrix of $\by',\bm{z}'$ in Problem (\ref{eq:twostageLP}), respectively. With the notations above we can write $C(\by,\bd)$ as the following problem
\begin{equation}\label{eq:twostageVec}
\begin{aligned} 
C(\by,\bm{d}) := &\min (\bm{c}')^\top \bm{z}'\\
&\begin{array}{r@{\quad}r@{}l@{\quad}l}
s.t. 
     & H\by  + W\bm{z}' &= \bd', \\
     &\bm{z}' &\geq \bm{0}.\\
\end{array}
\end{aligned}
\end{equation}

\subsection{Proof of Lemma~\ref{le:app3-properties-estimator}}\label{sec:app-app3-properties-estimator}
The proof of Lemma~\ref{le:app3-properties-estimator} is divided into the following two lemmas.
\begin{lemma}
    Given $\by \geq \bm{0}$ and $\bd$ a realization of demand $\bD$, the gradient estimator \eqref{estimator:app3} is unbiased.
\end{lemma}
\proof{Proof.}
Recall the definition of gradient estimator \eqref{estimator:app3}.
We write dual problem of (\ref{eq:twostageLP}) in the following form
\begin{equation}
    \max_{\bm{\pi}} \bm{\pi}^\top (\bd'- H\by ) 
    \quad\text{s.t. }  W^\top \bm{\pi} \leq \bm{c}'.
\end{equation}
For any $\by \geq 0$ and sample $\bd$, we know $C(\by,\bd)$ is finite. By proposition 2.2 of \cite{shapiro2021lectures}, we have
\begin{equation}
    \label{eq:app3estimator}
    \nabla_{\bm y} C(\by,\bd) = - H^\top \bm{\pi}^*(\by,\bd),
\end{equation}
where $$ \bm{\pi}^*(\by,\bd) = \argmax_{\bm{\pi} \in \{W^\top \bm{\pi} \leq \bm{c}'\}} \bm{\pi}^{\top} (\bd'- H\by ).$$
Note such a conclusion could also be obtained by Danskin's theorem.
We can choose any solution in $\bm{\pi}^*(\by,\bd)$ to get a gradient estimator $\hat{\nabla} Q(\by)$. By the conclusion of \cite{ermoliev1983stochastic}, we know $\E [\hat{\nabla} Q(\by)] = \nabla Q(\by)$. In other word, $\hat{\nabla} Q(\by)$ is unbiased gradient estimator of $Q(\by)$.
\Halmos
\endproof

Next, we prove the that gradient estimator $\hat{\nabla} Q(\by)$ is bounded.
\begin{lemma}
    There exist constant $\sigma_0$ such that, $$\Vert \hat{\nabla} Q(\by) \Vert _2 \leq \sigma_0.$$
    Moreover, $\Vert\nabla Q(\by) \Vert \leq \sigma_0$ and $\E \left[\Vert \hat{\nabla} Q(\by) \Vert _2^2 \right] \leq \sigma_0^2.$
\end{lemma}
\proof{Proof.}
It is easy to note that the primal problem (\ref{eq:twostageLP}) is feasible and has an optimal solution. Therefore the dual problem has a bounded optimal solution. To find a specific bound on the dual variable, we can analyze the structure of the dual program. Instead, we give an intuitive proof, based on the practical background of the second-stage problem.

Suppose that the inventory of the warehouse is $y_0$, and $\epsilon$ is a small positive number. If we change $y_0$ to $y_0+ \epsilon$, the extra $\epsilon$ inventory may increase the holding cost of the warehouse by $b_0 \epsilon$ or the extra inventory is delivered to the other store and decrease the lost-sales cost of some store. Noting the delivery requires some transportation cost, so the cost decreases less than $b_i \epsilon$. With the discussion above, we have that
$$\left\vert \frac{\partial C(\by,\bd)}{\partial y_0}\right\vert \leq h_0 + \max_{i \in [n]} b_i.$$
Similarly, for any store $i \in [n]$, adding $\epsilon$-inventory may at most increase holding cost $h_i \epsilon$ or decrease the lost-sales at most $b_i \epsilon$. Note that extra inventory at store $i$ may cause some inventory of the warehouse to be delivered to other stores, but we still have that the cost decreases less than $b_i \epsilon$. So we have $$\left\vert \frac{\partial C(\by,\bd)}{\partial y_i}\right\vert \leq h_i +b_i .$$

Define $\sigma_0 = \left[(h_0 + \sum_{i=1}^n b_i)^2+ \sum_{i=1}^n(h_i+b_i)^2 \right]^\frac 12$, then we know $$\left\Vert \hat{\nabla}Q(\by) \right\Vert_2 \leq \sigma_0 . $$
Note that $\left\Vert \hat{\nabla}Q(\by) \right\Vert_2 \leq \sigma_0$ implies that $\left\Vert\nabla Q(\by) \right\Vert \leq \sigma_0$ and $\E \left[\Vert \hat{\nabla} Q(\by) \Vert _2^2 \right] \leq \sigma_0^2.$ We complete the proof.
\Halmos
\endproof

\subsection{Proof of Lemma~\ref{app3:le:transition solver}}\label{sec:app-lemma-transition-solver}
\proof{Proof of Lemma~\ref{app3:le:transition solver}.}

If the inventory of all stores and the warehouse are consumed, the target will become feasible, which provides an upper bound for the hitting time. 

It is easy to note that each store will consume its own inventory and at least the $i^*$-th store consumes the inventory at the warehouse. 
Recalling that the inventory capacity for installation $i$ is $\rho_i$, by Lemma \ref{lem:hittingtime} we know the expected hitting time for each store is $1+6\beta_0 \rho_i$ and the expected hitting time for store $i^*$ is $1+6\beta_0 \rho_0+6\beta_0 \rho_{i^*}$. 

Therefore, we can bound $\E[s(\bx_0,\bw)]$ as follows
$$\E [s(\bx_0,\bw)] \leq n +6 \beta_0 \times (\rho_0+\rho_{i^*})+ 6\beta_0 \sum_{i=1,i \neq i^*}^n \rho_i = n+6\beta_0 \sum_{i=0}^n \rho_i,$$
where we complete the proof of this lemma.
\Halmos
\endproof

\subsection{Detailed Explanation about the Techniques in the Smoothness Proof of Lemma~\ref{le:app3-properties-qy}}\label{sec:app3smoothnesshighlevel}
\noindent\underline{\bf The sample-based analysis framework.}
Since $C(\by, \bd)$ is defined by a linear program (Problem \eqref{eq:delivery}), every value of $\by$, the gradient $\nabla_{\bm y} C(\by,\bd)$ is a piece-wise constant function of the sample $\bd$ -- indeed, the samples of $\bd$ in the same piece correspond to the same optimal solution of the LP (more precisely, the same extreme point of the feasible polyhedron), and give the same value of $\nabla_{\bm y} C(\by,\bd)$. These pieces naturally divide the probability space of $\bD$ into regions.
Each region of the probability space will shift slightly if we perturb $\by$ to $\by + \Delta \by$. This shift may cause a sample $\bd$ to give a different gradient value, and thus $\nabla_{\by} C(\by+\Delta \by,\bd)- \nabla_{\by} C(\by,\bd)$ becomes non-zero. Taking expectation over $\bd$, and assuming that the function $f_{\bD}(\cdot)$ is bounded, we can upper bound $\E[\|\nabla_{\by}C(\by+\Delta \by,\bd)-\nabla_{\by}C(\by,\bd)\|]$ by the volume of the set of $\bd$ whose gradient value is affected by the shift (multiplied by the upper bound of $f_{\bD}(\cdot)$), and hence upper bound $\|\nabla Q(\by+\Delta \by)-\nabla Q(\by)\|$, which relates to the smoothness of $Q(\by)$.

To bound the volume of $\bd$ affected by the shift due to $\Delta\by$, \cite{wang1985distribution} naturally related this quantity to the measure of the boundaries separating the above-mentioned regions which further relates to the number of the regions. Since the number of the regions can be as many as the extreme points of the LP polyhedron (which may be exponential in $n$), the smoothness bound derived following the method by \cite{wang1985distribution} is at the order of $2^{\O(n)}$.

\noindent\underline{\bf Exponential improvement via cost decomposition and coordinate-wise analysis.} We observe that the optimization problem $C(\by, \bd)$ can be solved by greedily delivering the inventory to the most profitable stores first. This helps us to decompose the total cost $C(\by, \bd)$ to the summation of $\O(n)$ terms, namely $\{C_k(\by, \bd)\}$ (the detailed definition of $C_k$ is deferred to Section \ref{appendix:app3}). These $C_k$'s correspond to the costs of different types of installations. We would like to apply the above-mentioned sample-based analysis to upper bound the smoothness parameter of $\mathbb{E}[C_k(\by, \bd)]$ for each $k$, and the smoothness bound of $Q(\by)$ is at most $\O(n)$ times the maximum bound among these terms. 

For each term $C_k$, however, if we apply the same method of \cite{wang1985distribution} for the smoothness bound, we would find that the number of pieces (regions) associated with $\nabla_{\by}C_k(\by,\bd)$ might still be exponential in $n$. Fortunately, for the one-warehouse multi-store system, we are able to leverage the structural property of each term $C_k$ and apply a coordinate-wise analysis --- for each coordinate $i \in [n]$, the partial derivative $\partial C_k(\by, \bd) / \partial y_i$ turns out to be a piece-wise constant function with only $\O(1)$ pieces. In this way, we may bound $\E[|\partial C_k(\by + \Delta \by, \bd) / \partial y_i - \partial C_k(\by, \bd) / \partial y_i|]$ by the perturbation $\|\Delta\by\|$ up to only an extra factor of $n^{\O(1)}$ (instead of $2^{\O(d)}$ as in the general case). Finally, we may polynomially upper bound the gradient difference $\E[\|\nabla_{\by}C(\by+\Delta \by,\bd)-\nabla_{\by}C(\by,\bd)\|]$ since the difference along each coordinate direction is bounded.

\subsection{Proof of Lemma~\ref{le:app3-properties-qy}}\label{proofofsmoothnessofapp3}

We divide the proof of Lemma~\ref{le:app3-properties-qy} into two parts: convexity and smoothness.
\proof{Proof of convexity.}
In the dual problem \eqref{problem:dual},  $C(\by,\bd)$ is defined as the maximization of a set of linear functions $\bm{\pi}^\top(\bd'- H\by )$ indexed by $\{\bm{\pi}\mid W^\top \bm{\pi} \leq \bm{c}'\}$. Since convexity is preserved in maximization operation, we know $C(\by,\bd)$ is convex. Consequently, $Q(\by) = \E[C(\by,\bD)]$ is also convex. 
 \Halmos \endproof
\begin{algorithm}[t]
\caption{Optimal Solution of Problem (\ref{eq:delivery}) }
\label{alg:greedy}
\begin{algorithmic}[1]
\State \textbf{Input: $h_0,\dots,h_n$, $b_i,i\in[n]$, $c_i,i\in[n]$, $y_0,\dots,y_n$ and 
 $d_i,i\in[n].$ } Without loss of generality we assume $\{b_i-c_i\}_{i\in[n]}$ is in a descending order (i.e., $b_1-c_1\geq b_2-c_2\geq  \dots\geq b_n-c_n$) and initialize $z_i=0$ for all $i\in[n]$.
\While{the inventory in the warehouse is not depleted ($y_0>0$)}
\For{$i = 1,2,\dots,n$}
\If{$d_i>z_i$} 
\State $z_i= \max\{y_0, d_i-y_i\}$.
\Else{ $z_i=0$}
\EndIf
\State $y_0= y_0-z_i$.
\EndFor
\EndWhile
\State \textbf{Output:} $z_i,i\in[n].$
\end{algorithmic}
\end{algorithm} 
\proof{Proof of smoothness.}
We adopted an innovative sample-based analysis to prove the smoothness of $Q(\by)$. Recall that $Q(\by) = \E[C(\by,\bD)]$ and $C(\by,\bd)$ is defined by the second-stage problem \eqref{eq:twostage} given sample $\bd$. By definition, $C(\by,\bd)$ is the cost of all installations when the problem \eqref{eq:twostage} is solved optimally. 

First, we make an important observation that problem \eqref{eq:twostage} can be solved greedily (Algorithm~\ref{alg:greedy}). 
Specifically, if store $i$ needs to replenish its inventory from the warehouse, the transportation of inventory per unit from the warehouse to the store will result in a reduction of $b_i-c_i$ cost.
Thus it holds that the optimal solution is to deliver inventory to all the stores requiring delivery in the descending order of $\{b_i-c_i\}_{i \in [n]}$. With such an observation, we assume that all stores are relabelled in the descending order of $\{b_i-c_i\}_{i \in [n]}$ without loss of generality. Therefore,  the optimal strategy is to deliver inventory to all stores in need, starting from store $1$ up to store $n$, until the warehouse's inventory is exhausted.

Suppose that demand is realized as $\bd$ and order-up-to level is $\by$. Under the optimal delivery decision, we introduce some notations to describe the delivery process as follows.  
\begin{enumerate}[nolistsep]
    \item Let $s_i,i \in [n]$ be the remaining inventory after delivery to the first $i$ stores. Note that $s_n$ is the left-over inventory of the warehouse, after delivery to all stores. For convenience, we define $s_0=y_0$.
    \item Let $z_i,i \in [n]$ be the inventory delivered to the $i$-th store. 
    \item Let $o_i,i \in [n]$ be the left-over inventory of the $i$-th store, after the delivery from the warehouse.
    \item Let $l_i,i \in [n]$ be the unsatisfied demand of the $i$-th store, after the delivery from the warehouse.
\end{enumerate}
With the notations, it is easy to verify that $s_i,z_i,o_i,l_i,i \in [n]$ satisfy the following system dynamics
\begin{align*}
    s_i &= (y_0-(d_1-y_1)^+ -\dots - (d_i-y_i)^+)^+, i\in [n], \\
    z_i &= \min ((d_i-y_i)^+,s_{i-1}),i \in [n], \\
    o_i &= (y_i-d_i)^+,i \in [n],\\
    l_i &= ((d_i-y_i)^+-s_{i-1})^+, i \in [n].
\end{align*}
and we can decompose $C(\by,\bd)$ as
\[
C(\by,\bd) = h_0s_n+\sum_{i=1}^n (c_i z_i +b_i l_i + h_i o_i). 
\]
Note that all $s_i,z_i,o_i,l_i,i \in [n]$ are functions of $\by,\bd$. We aim to analyze each term separately and to prove that each term is smooth under expectation, which further implies the smoothness of $Q(\by)$. 

We explain the idea of the sample-based analysis. First, we fix a sample $\bd$. Then, we analyze the solution under inventory $\by$. Next, we perturb $\by$ to $\by'$ and examine how this perturbation affects the terms $s_i, z_i, o_i,$ and $l_i$ for all $i\in [n]$. Let $\by' = \by + \Delta \by$, where $\Delta \by = (\Delta y_0,\dots,\Delta y_n)^{\top}$. We assume, without loss of generality, that $\Delta y_i \geq 0$ for all $i=0,\dots,n$. For ease of notation, we define $\Delta \by_k = (0,\dots,0,\Delta y_k,0,\dots,0)$, where $\Delta y_k$ is the $k$-th component. Next, we consider each $\Delta \by_k$, $k=0,\dots,n$, and use the triangle inequality to draw conclusions about $\Delta \by$. 

\noindent \underline{\bf Analysis of terms $o_i,i \in [n]$ as an illustrative example.} Obviously, by the definition of $o_i = (y_i-d_i)^+,i \in [n]$, we know that if $i \neq j$,  $\nabla_{y_i} o_j(\by,\bd) = 0$. Then we consider $\nabla_{y_i} o_i(\by,\bd)$. Since
\begin{equation}\nonumber
    \nabla_{y_i} o_i(\by,\bd) = 
    \begin{cases}
         1, & \text{if } y_i \geq d_i;   \\
         0, & \text{if } y_i < d_i,
    \end{cases}
\end{equation}
it holds that the difference $|\nabla_{y_i} o_i(\by+\Delta \by_i,\bd) - \nabla_{y_i} o_i(\by,\bd)|=1$ if and only if $y_i < d_i$ and $y_i + \Delta y_i \geq d_i$. Direct computing gives
\begin{align*}
    \left\vert\E[\nabla_{y_i} o_i(\by+\Delta \by_i,\bD)] - \E[\nabla_{y_i} o_i(\by,\bD)] \right\vert
    &= |\E_{\bd \sim \bD}[\nabla_{y_i} o_i(\by+\Delta \by_i,\bd)-\nabla_{y_i} o_i(\by,\bd)]|\\
    &\leq \E_{\bd \sim \bD}[|\nabla_{y_i} o_i(\by+\Delta \by_i,\bd)-\nabla_{y_i} o_i(\by,\bd)|] \\
    &\leq \P_{\bd \sim \bD} \{y_i \leq d_i \leq y_i + \Delta y_i \}\\
    &\leq \beta_0 \Delta y_i.
\end{align*}
where the last inequality is because the density function of $\bD$ is bounded by $\beta_0$ and for $j \neq i$, $|\E[\nabla_{y_i} o_i(\by+\Delta \by_j,\bd)]- \E[\nabla_{y_i} o_i(\by,\bd)]|=0$. With all the discussions above, we conclude that
\[\nabla_{y_i} \E[\nabla_{y_i} o_i(\by,\bd)]\leq \beta_0 \text{  ~~and~~   } \nabla_{y_j} \E[\nabla_{y_i} o_i(\by,\bd)]= 0, \]
for all $i \in [n],j = 0,\dots, n \text{ and } j\neq i $. 
Therefore, we have $\|\nabla \E[\nabla o_i(\by,\bD)]\|_F\leq \beta_0$
By the quasi-mean-value theorem, we have $\E[\nabla o_i(\by,\bD)]$ is $\beta_0$-Lipschitz, which implies that $\E[o_i(\by,\bD)]$ is $\beta_0$-smooth.

\noindent \underline{\bf Analysis of term $s_n$.}
We consider $\nabla_{y_i} s_n(\by,\bd)$, for each $i = 0,\dots,n$, separately.

If $i=0$, by the definition of $s_n$, we know
\begin{equation}\nonumber
    \nabla_{y_i} s_n(\by,\bd) = 
    \begin{cases}
         1, & \text{if } y_0 \geq \sum_{i'=1}^n (d_{i'}-y_{i'})^+;   \\
         0, & \text{if } y_0 < \sum_{i'=1}^n (d_{i'}-y_{i'})^+.
    \end{cases}
\end{equation}
Consider adding inventory $\Delta y_k$ to installation $k$ for any fixed $k=0,\dots,n$. There are two cases in which the value of $\nabla_{y_i} s_n(\by,\bd)$ changes.
\begin{enumerate}[label=\textbf{Case \arabic*:}]
    \item
    Before adding $\Delta \by_k$, we have $y_0 < \sum_{i'=1}^n (d_{i'}-y_{i'})^+$. After adding $\Delta \by_k$, this inequality is violated due to that $y_0$ changes to $y_0 + \Delta y_0$. This case only happens when $k=0$.
    We can infer that before adding $\Delta \by_k$, there exists a store $k'$ for which the delivery from the warehouse is insufficient to satisfy its demand. Specifically, for store $k'$, we have $z_{k'}+y_{k'} < d_{k'}$. After adding $\Delta y_k$ inventory, the warehouse delivery can now satisfy the demand at store $k'$. Thus, we have $d_{k'} \leq z_{k'} + y_{k'}+ \Delta y_k$.
Given the assumption that $f_{\bD}(\bx) \leq \beta_0$, we have $\P \{z_{k'}+y_{k'} \leq d_{k'} \leq z_{k'}+y_{k'}+\Delta y_k \} \leq \beta_0 \Delta y_k$. Furthermore, there are at most $n$ possible stores that can cause such an event. Using the union bound, we conclude that under case 1, the probability of the event $|\nabla_{y_i} s_n(\by+\Delta \by_k,\bd) - \nabla_{y_i} s_n(\by,\bd)|=1$ is at most $n\beta_0 \Delta y_k$.
Similar to the proof of term $o_i$, we can obtain $\nabla_{y_k} \E[\nabla_{y_i} s_n(\by,\bd)]\leq n\beta_0$, for $i=0,k=0.$
    \item
    Before adding $\Delta \by_k$, we have $y_0 < \sum_{i'=1}^n (d_{i'}-y_{i'})^+$. After adding $\Delta \by_k$, this inequality is violated due to that $y_0$ changes to $y_0 + \Delta y_0$. This case only happens when $k\geq 1$.
    We can infer that before adding $\Delta \by_k$, there exists a store $k'$ such that delivery from the warehouse is not enough to satisfy its demand. After adding extra $\Delta y_k$ inventory, the delivery from the warehouse can satisfy the demand at the store $k'$. Note that $\Delta y_k$ is added to the store $k$, which causes at most $\Delta y_j$ inventory delivered to the store $k'$. It follows that $\P\{z_{k'}+y_{k'} \leq d_{k'} \leq z_{k'}+y_{k'}+\Delta y_k \} \leq \beta_0 \Delta y_k$. The event that $|\nabla_{y_i} s_n(\by+\Delta \by_k,\bd) - \nabla_{y_i} s_n(\by,\bd)|=1$ has probability at most $n\beta_0 \Delta y_k.$ We obtain  $\nabla_{y_k} \E[\nabla_{y_i} s_n(\by,\bd)]\leq n\beta_0$, for $i=0,k\geq1.$
\end{enumerate}
If $i \geq 1$, by the definition of $s_n$, we know
\begin{equation}\nonumber
    \nabla_{y_i} s_n(\by,\bd) = 
    \begin{cases}
         1, & \text{if } y_0 \geq \sum_{i'=1}^n (d_{i'}-y_{i'})^+ \text{ and } d_i \geq y_i;   \\
         0, & \text{otherwise}.
    \end{cases}
\end{equation}
For any $k = 0,\dots,n$, consider adding $\Delta y_k$ inventory to the installation $k$. There are three cases that may change the value of $\nabla_{y_i} s_n(\by,\bd)$.
\begin{enumerate}[label=\textbf{Case \arabic*:}]
    \item
    When $k=0$, $y_0 < \sum_{i'=1}^n (d_{i'}-y_{i'})^+$ is violated after adding $\Delta \by_k$ due to that $y_0$ changes to $y_0 + \Delta y_k$. Similar to Case 1 when $i=0$, we obtain the event has probability at most $n\beta_0 \Delta y_k$.
    \item
    When $k \geq 1$, $y_0 < \sum_{i'=1}^n (d_{i'}-y_{i'})^+$ is violated after adding $\Delta \by_k$ due to that $\sum_{i'=1}^n (d_{i'}-y_{i'})^+$ decreases after adding $\Delta \by_k$. Similar to Case 2 when $i=0$, we obtain the event has probability at most $n\beta_0 \Delta y_k$.
    \item When $k = i$, $d_i \geq y_i$ is violated after adding $\Delta \by_k$, due to that $y_i$ increases. It follows that $\P\{ y_i + \Delta y_k \geq d_i \geq y_i \} \leq \beta_0 \Delta y_k$.
\end{enumerate}
Combining the three cases above, we obtain $\nabla_{y_k} \E[\nabla_{y_i} s_n(\by,\bd)]\leq (n+1)\beta_0$. Similar to the analysis of $o_i,i \in [n]$, we can prove analogously that $\E[s_n(\by,\bD)]$ is $(n+1)^2\beta_0$-smooth.

\noindent \underline{\bf Analysis of terms $l_i,i \in [n]$.} We fix a $i \in [n]$ and analyze $\nabla_{y_j} l_i(\by,\bd)$ for all $j$.

If $j=0$, then direct computing gives
\begin{equation}\nonumber
    \nabla_{y_j} l_i(\by,\bd) = 
    \begin{cases}
         -1, & \text{if } (d_i-y_i)^+ > s_{i-1} \text{ and } y_0 \geq \sum_{i'=1}^{i-1} (d_{i'}-y_{i'})^+;   \\
         0, & \text{otherwise}.
    \end{cases}
\end{equation}
For any $k =0,\dots,n$, consider adding $\Delta \by_k$ inventory. There are four cases that may change the value of $\nabla_{y_j} l_i(\by,\bd)$. 
\begin{enumerate}[label=\textbf{Case \arabic*:}]
    \item When $k=i$, $(d_i-y_i)^+ > s_{i-1}$ is violated after adding $\Delta \by_k$ due to that $(d_i-y_i)$ decreases, We obtain that $s_{i-1}<d_i-y_i$ and $d_i-y_i-\Delta y_k \leq s_{i-1}$. It is easy to see that $\P \{s_{i-1}< d_i-y_i\leq \Delta y_k+ s_{i-1}\} \leq \beta_0 \Delta y_k.$
    \item For any $k$, $(d_i-y_i)^+ > s_{i-1}$ is violated after adding $\Delta \by_k$ due to that $s_{i-1}$ decreases. Adding $\Delta y_k$ inventory to store $k$ will save at most $\Delta y_k$ inventory to the warehouse, which implies that $ s_{i-1} + \Delta y_k\leq d_i-y_i > s_{i-1}$.
    \item  When $k=0$, $y_0 < \sum_{i'=1}^{i-1} (d_{i'}-y_{i'})^+$ is violated after adding $\Delta \by_k$ due to that $y_0$ increases to $y_0 + \Delta y_k$. We obtain that $s_{i-1}<d_i-y_i$ and $d_i-y_i-\Delta y_k \leq s_{i-1}$. It is easy to see that $\P \{s_{i-1}< d_i-y_i\leq \Delta y_k+ s_{i-1}\} \leq \beta_0 \Delta y_k.$
    \item For any $k$, $y_0 < \sum_{i'=1}^{i-1} (d_{i'}-y_{i'})^+$ is violated after adding $\Delta \by_k$ due to that $\sum_{i'=1}^{i-1} (d_{i'}-y_{i'})^+$ decreases. The analysis of the case is similar to Case 2 of $s_n$ when $i=0$.
\end{enumerate}
Combining the four cases above, we obtain $\nabla_{y_k} \E[\nabla_{y_j} l_i(\by,\bd)]\leq (n+3)\beta_0$ for $j=0.$

If $0 < j < i$, then direct computing gives
\begin{equation}\nonumber
    \nabla_{y_j} l_i(\by,\bd) = 
    \begin{cases}
         -1, & \text{if } (d_i-y_i)^+ > s_{i-1} , d_j > y_j,\text{ and } y_0 \geq \sum_{i'=1}^{i-1} (d_{i'}-y_{i'})^+;   \\
         0, & \text{otherwise}.
    \end{cases}
\end{equation}
For any $k =0,\dots,n$, consider adding $\Delta \by_k$ inventory. There are two cases that may change the value of $\nabla_{y_j} l_i(\by,\bd)$.
\begin{enumerate}[label=\textbf{Case \arabic*:}]
    \item  When $k=i$, $(d_i-y_i)^+ > s_{i-1}$ is violated after adding $\Delta \by_k$. The analysis is similar to Case 1 when $j=0.$
    \item  For any $k$, $(d_i-y_i)^+ > s_{i-1}$ is violated, since $s_{i-1}$ decreases after adding $\Delta \by_k$. The analysis is similar to Case 2 when $j=0.$
    \item  When $k=j$, $d_j > y_j$ is violated after adding $\Delta \by_k$, which implies that $y_j + \Delta y_k \geq d_j \geq y_j.$
    \item  When $k=0$, $y_0 < \sum_{i'=1}^{i-1} (d_{i'}-y_{i'})^+$ is violated after adding $\Delta \by_k$. The analysis is similar to Case 3 when $j=0.$
    \item  For any $k$, $y_0 < \sum_{i'=1}^{i-1} (d_{i'}-y_{i'})^+$ is violated, since $\sum_{i'=1}^{i-1} (d_{i'}-y_{i'})^+$ decreases after adding $\Delta \by_k$. The analysis is similar to Case 4 when $j=0.$
\end{enumerate}
Combining the five cases above, we obtain $\nabla_{y_k} \E[\nabla_{y_j} l_i(\by,\bd)]\leq (n+4)\beta_0$ for $0 < j < i.$

If $j = i$, then direct computing gives
\begin{equation}\nonumber
    \nabla_{y_j} l_i(\by,\bd) = 
    \begin{cases}
         -1, & \text{if } (d_i-y_i)^+ > s_{i-1} \text{ and } d_i \geq y_i;   \\
         0, & \text{otherwise}.
    \end{cases}
\end{equation}
For any $k =0,\dots,n$, consider adding $\Delta \by_k$ inventory. There are three cases that may change the value of $\nabla_{y_j} l_i(\by,\bd)$.
\begin{enumerate}[label=\textbf{Case \arabic*:}]
    \item  When $k=i$, $(d_i-y_i)^+ > s_{i-1}$ is violated after adding $\Delta \by_k$, since $(d_i-y_i -\Delta y_i) \leq s_{i-1}.$ The analysis is similar to Case 1 when $j=0.$
    \item  For any $k$, $(d_i-y_i)^+ > s_{i-1}$ is violated, since $s_{i-1}$ decreases after adding $\Delta \by_k$. The analysis is similar to Case 2 when $j=0.$
    \item  When $k=i$, $d_i \geq y_i$ is violated after adding $\Delta \by_k$. The analysis is similar to Case 3 when $0 <j<i.$
\end{enumerate}
Combining the three cases above, we obtain $\nabla_{y_k} \E[\nabla_{y_j} l_i(\by,\bd)]\leq 3\beta_0$ for $j=i.$

If $j > i$, then it is easy to find that $\nabla_{y_j} l_i(\by,\bd) =0$, since $y_j$ will not affect the delivery of the stores before the store $j$ and not appears in $l_i(\by,\bd).$

With all the discussion above, we can prove that $\E[l_i],i \in [n]$ is $(n+4)(n+1)\beta_0$-smooth.

\noindent \underline{\bf Analysis of terms $z_i,i \in [n]$.} It is remarkable to observe the following equation
\[
l_i + z_i = (d_i-y_i)^+,i \in [n],
\]
since the quantity $(d_i-y_i)^+$ required by store $i$ is partially satisfied by actually delivery quantity $z_i$ from the warehouse and the other is unsatisfied $l_i$.
The analysis of term $(d_i-y_i)^+$ is analogous to the term $o_i = (y_i-d_i)^+$, therefore it is easy to verify that $\E[(d_i-y_i)^+]$ is $\beta_0$-smooth.
Combining with the conclusion of $\E[l_i]$, we know that $\E[z_i],i \in [n]$ is $[(n+4)(n+1)+1]\beta_0$-smooth.

\noindent \underline{\bf Final step of proof.} Recall the definition of $C(\by,\bd)$ is
\[
C(\by,\bd) = h_0s_n+\sum_{i=1}^n c_i z_i +b_i l_i + h_i o_i. 
\]
All the terms on the right-hand side are proved to be smooth under expectation. Combining their smoothness coefficients, we conclude that $Q(\by) = \E[C(\by,\bD)]$ is $\beta$-smooth, and $\beta$ is given by
\begin{equation}
    \beta = (n+1)^2\beta_0h_0+\beta_0 \sum_{i=1}^n h_i + (n+4)(n+1)\beta_0 \sum_{i=1}^n b_i+[(n+4)(n+1)+1]\beta_0 \sum_{i=1}^n c_i .
\end{equation}
It is easy to obtain a compact upper bound as follows
\begin{equation}
\label{eq:app3beta}
    \beta \leq (n+4)^2\beta_0 \left(h_0 + \sum_{i=1}^n (h_i + b_i + c_i)\right) = \O (n^3).
\end{equation}


{
\section{Details of Numerical Experiments}\label{detailedexp}
In this section, we present the detailed setup of our numerical experiments and the numerical results. In the following, we first introduce some general settings used in this section. Then in each subsection, we present the details and numerical results as mentioned in Section~\ref{sec:numerical experiments}.
All simulations were implemented by MATLAB R2023a on an AMD Ryzen 3.00GHz PC.

\medskip \noindent  {\bf Demand distributions.} We conduct our numerical experiments under several different demand distributions. The distributions include the uniform distribution $U(a,b)$ supported on $[a,b]$, the Poisson distribution $\mathrm{Poiss}(\lambda)$ with mean $\lambda$, the geometric distribution $\mathrm{Geom}(p)$ with success probability $p$, the gamma distribution $\mathrm{Gamma}(r,\lambda)$ with shape parameter $r$ and inverse scale parameter $\lambda$, and normal distribution $N(\mu,\sigma^2)$ with mean $\mu$ and variance $\sigma^2$.\footnote{{We set the negative realization of gamma and normal distribution to be zero, and we still use $\mathrm{Gamma}(r,\lambda)$ and $N(\mu,\sigma^2)$ to represent the clipped distribution for simplicity.}}
For different experiments, the distribution parameters are either written in the title of corresponding figures or specified in the description of the experiments.

\medskip \noindent  {\bf Evaluation metrics. }The performance of an algorithm $\pi$ is measured by the \emph{relative average regret} defined as
\begin{equation}
\nonumber
    \frac{C^{\pi}_P(T)-T\cdot C^*_P}{T\cdot C_P^{*}} \times 100 \%,
\end{equation}
where $P$ denotes the problem instance, $C_P^{\pi}(T)$ is the expected cost incurred by the algorithm $\pi$ and $C_P^*$ is the optimal cost. Except for the experiments in Section \ref{subsec:num-mpmc},\footnote{{Since solving the quadratic optimization problem for the experiments in Section \ref{subsec:num-mpmc} is quite time-consuming, we repeat them 10 times to evaluate the relative average regret.}} we repeat all experiments 1000 times to evaluate the relative average regret. The random seed is set as $1$.

}{

\medskip \noindent {\bf Labels of the tested algorithms. } If there is no additional description in an experiment, we generally label and set the hyperparameters of SGD and our minibatch-SGD-based meta-policy as follows.
\begin{enumerate}
    \item {\bf Algorithm derived from minibatch-SGD-based meta-policy (MS-X-Y):} The notation X represents the problem and Y is the batch size scheme used by the algorithm (Y$=$L for the linear batch size scheme and Y$=$E for the exponential batch size scheme). The stepsize is set to be a constant $\eta$. The linear batch size scheme is fixed as $n_{\tau} = \tau$ and the base of the exponential batch size scheme is chosen by cross-validation and grid search in the set $\{1.05,1.15,1.25,1.5,2\}$.
    \item {\bf Algorithm based on Stochastic Gradient Descent (SGD-Z):} The notation Z represents the stepsize pattern of the algorithm (Z$=0.5$ for the stepsize pattern $\eta/t^{0.5}$ and Z$=1$ for the stepsize pattern $\eta/t^{1}$). When the target order-up-to level updated by SGD is feasible, the algorithm sets this level as the actually implemented order-up-to level. Otherwise, the algorithm sets the actually implemented order-up-to level following the suggestion of the transition solver used by the meta-policy.
\end{enumerate}

 The parameter $\eta$ for the stepsize in both algorithms is chosen by cross-validation and grid search in the set $\{0.01,0.03,0.1,0.3,1,3,10,30,100\}$. The initial inventory of all experiments and the initial order-up-to level of all algorithms is set to be $0$.

\subsection{Single Product Newsvendor Problem}
\label{subsec:num-nvp}
In this section, we conduct numerical experiments for the single product NewsVendor Problem (NVP) under four different distributions ($N(5,1)$, $U(0,10)$, $\mathrm{Poiss}(5)$, and $\mathrm{Geom}(0.2)$) and set per-unit holding cost $h=1$ and per-unit lost-sale cost $b=50$. We test the performance of the following classic algorithms in the literature.
\begin{enumerate}[nolistsep]
    \item \textbf{Minibatch SGD (MS-NV):} This algorithm is our Minibatch-SGD-based (MS) meta-policy applied to the NewsVendor (NV) problem.
    \item \textbf{Stochastic Gradient Descent (SGD)} in \cite{huh2009nonparametric}: This algorithm uses stochastic gradient descent to update the order-up-to level.
    \item \textbf{Sample Average Approximation (SAA)} in \cite{levi2015data}: This algorithm uses all history demand data to construct an empirical quantile of demand and applies the solution as order-up-to level in the next period. This algorithm requires uncensored demand data.
    \item \textbf{Kaplan-Meier (KM)} in \cite{huh2011adaptive}: This algorithm uses the Kaplan-Meier estimator from survival analysis to estimate the quantile of demand from the censored data.
    \item \textbf{Binary Search (BS)} in \cite{chen2020optimal}: This algorithm constructs confidence intervals of the gradient of cost function combined with binary search to find optimal order-up-to level.
\end{enumerate}}

{
For each algorithm, we fine-tune their parameters to achieve their peak performance for meaningful comparisons. For example, we choose the stepsize pattern of SGD to be $\eta/t$ and the batch size of MS-NV to be exponential. The relative average regret curves are presented in Figure~\ref{fig:nvp}.

Through the numerical results presented in Figure~\ref{fig:nvp}, we note that MS-NV derived from our meta-policy performs very well and is competitive among many classical algorithms.
Specifically, the performance of MS-NV is only slightly worse than that of SAA. This is reasonable since SAA requires full observation of the demands, while our algorithm only requires the censored demand.

\begin{figure}[th]
\begin{center}
\includegraphics[width =0.48\textwidth]{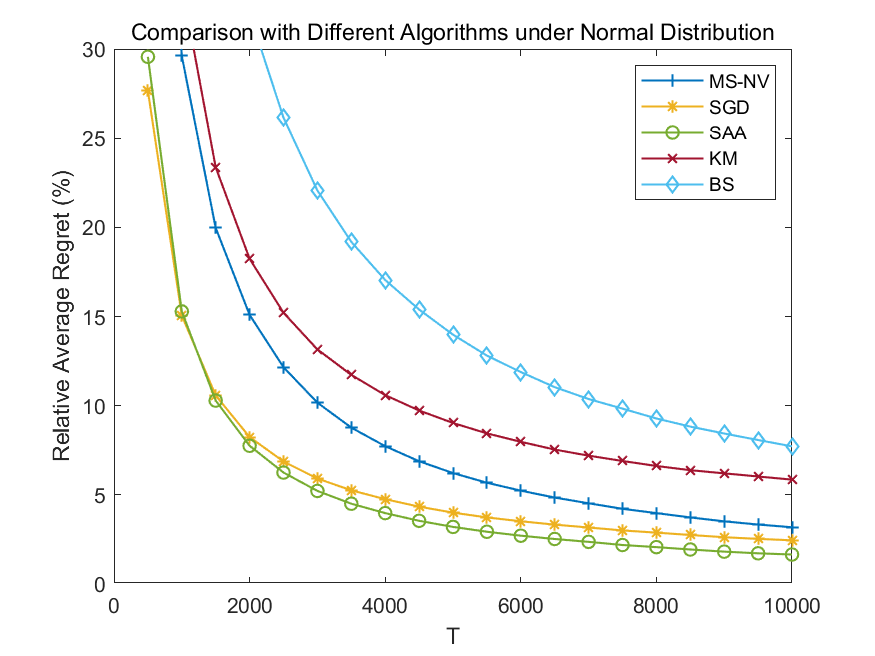}
\includegraphics[width =0.48\textwidth]{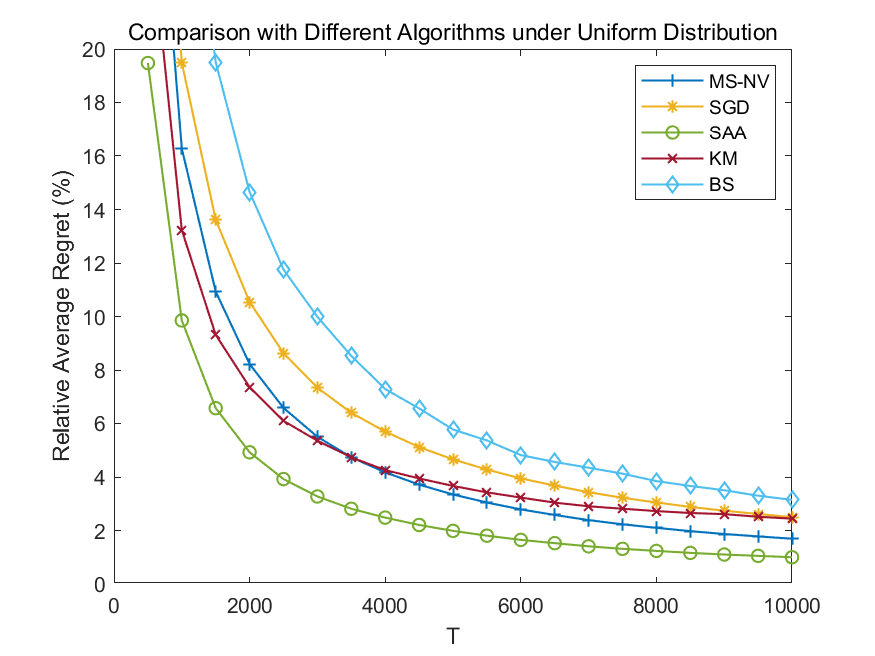}
\includegraphics[width =0.48\textwidth]{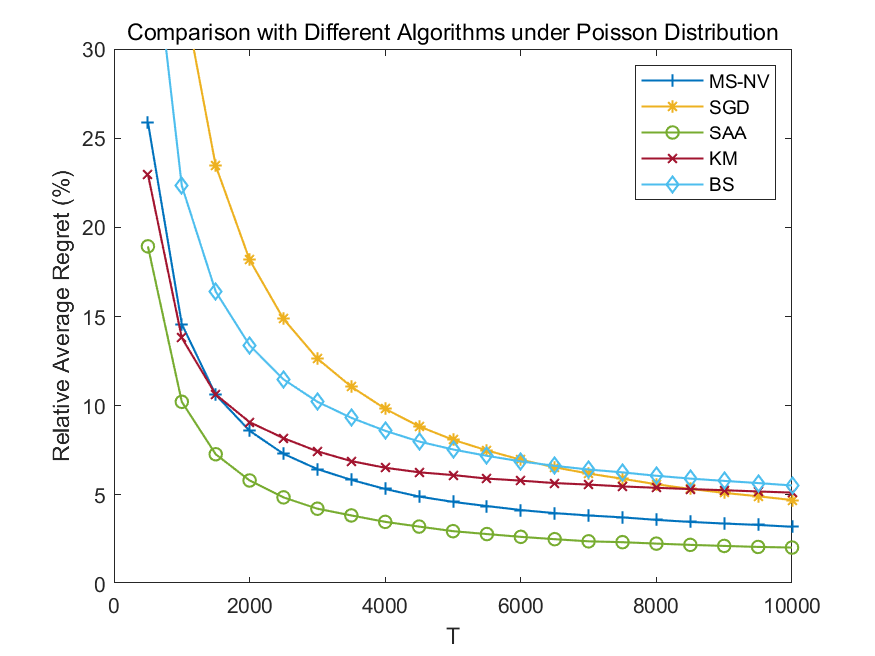}
\includegraphics[width =0.48\textwidth]{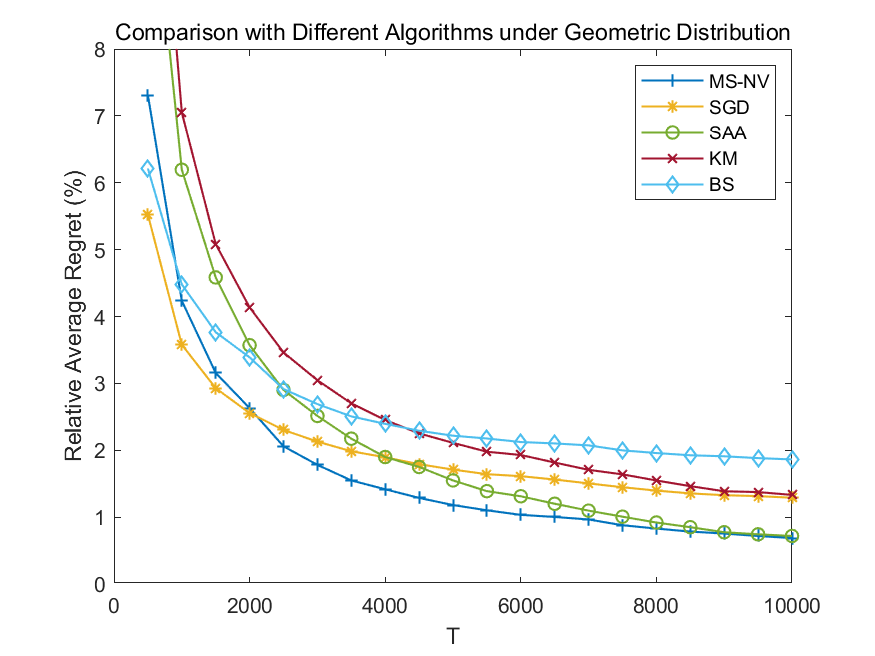}
\end{center}
\caption{{NVP: Comparision with Different Algorithms under Different Distributions}}\label{fig:nvp}
\end{figure}

{
\medskip \noindent {\bf Robustness in parameter $b$.} In most experiments, we choose $b \gg h$ because in many practical scenarios the lost-sales penalty $b$ is generally much larger than the holding cost $h$. For example, in retail and service parts maintenance \citep{huh2009asymptotic}, the storage of products is cheap, and $h$ will be small. The lost sales will cause a loss in revenue and customer goodwill, which is significantly larger than the holding cost $h$. To test the robustness of the minibatch SGD algorithm, we conduct numerical experiments and compare it with other algorithms under normal distribution and smaller parameters $b=5$ and $b=25$. The numerical results are presented in Figure.~\ref{fig:nvp-robust}. In \ref{subsec:num-leadtime}, we also conduct numerical experiments for inventory systems with lead times under different parameters $b$. In these figures, we can see minibatch SGD performs well and is competitive compared with other algorithms. These experiments show the robustness of minibatch SGD in parameter $b$.

\begin{figure}[th]
\begin{center}
\includegraphics[width =0.48\textwidth]{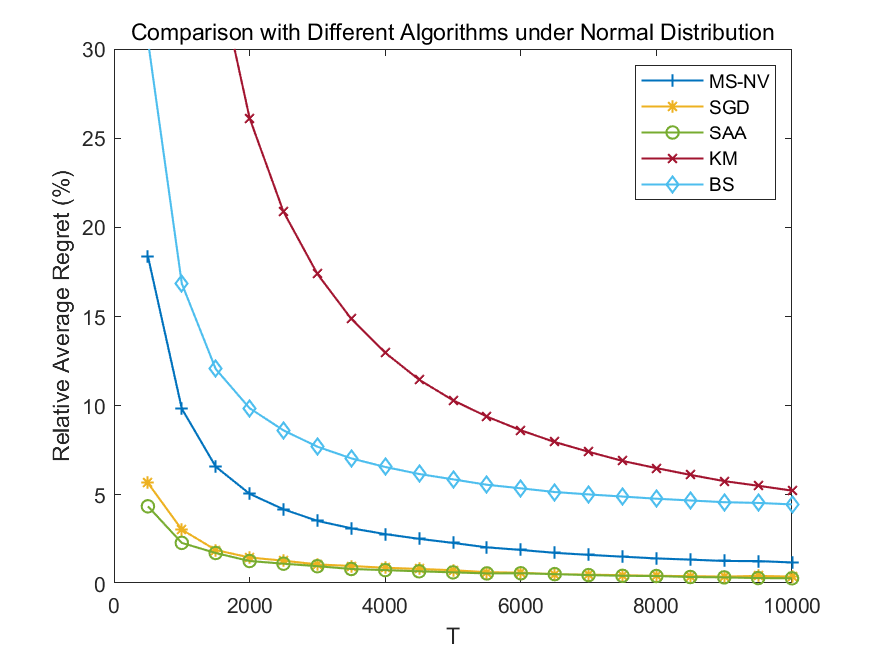}
\includegraphics[width =0.48\textwidth]{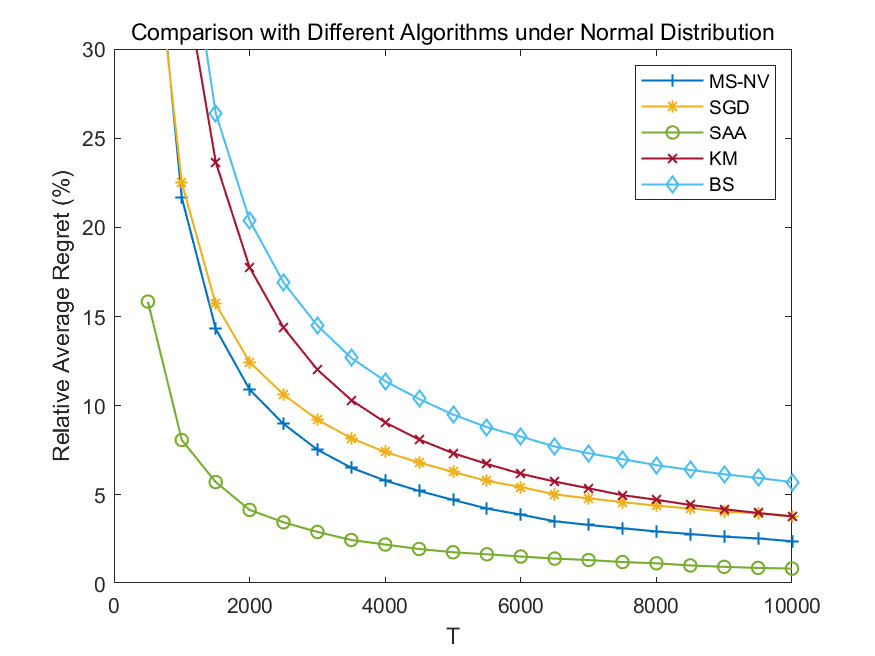}
\end{center}
\caption{{NVP: Comparision with Different Algorithms under Different $b$ (Left Panel: $b=5$, Right Panel: $b = 25$)}}\label{fig:nvp-robust}
\end{figure}
}

\medskip \noindent {\bf Variance of the learned solution.} 
We compare the variances of the learned solution of the minibatch-SGD based methods and the SGD-based methods. we consider the experiment shown on the top-right of Figure~\ref{fig:nvp} and draw the relative average regret curve and the curve of the distance from the optimal solution to learned one, and include the standard deviations. The numerical results are presented in Figure~\ref{fig:nvp-var}. The shaded regions represent the standard deviation from the mean curve. We can see that regret and distance from the optimal solution of MS-NV have smaller variances in these instances.

\begin{figure}[th]
\begin{center}
\includegraphics[width =0.48\textwidth]{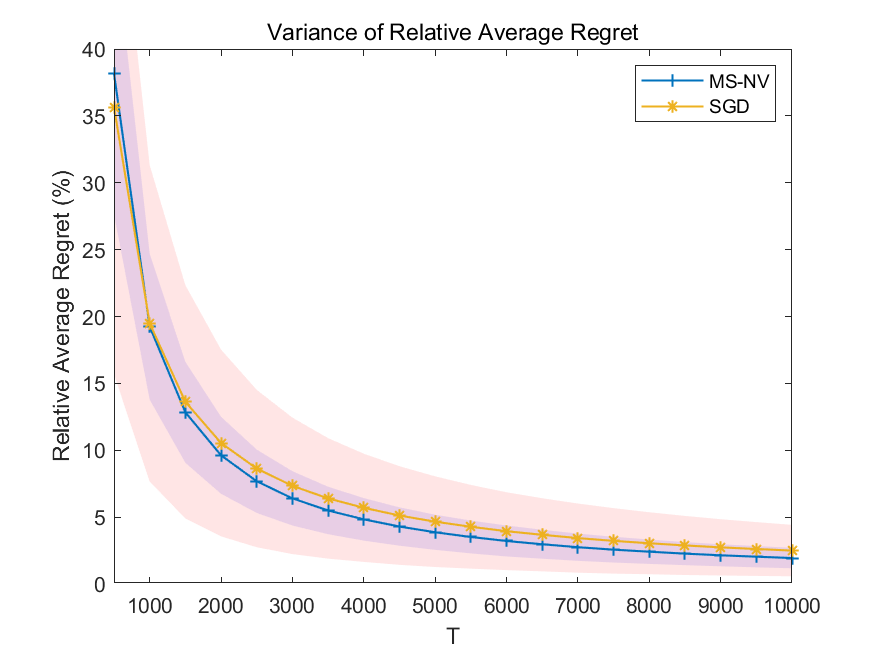}
\includegraphics[width =0.48\textwidth]{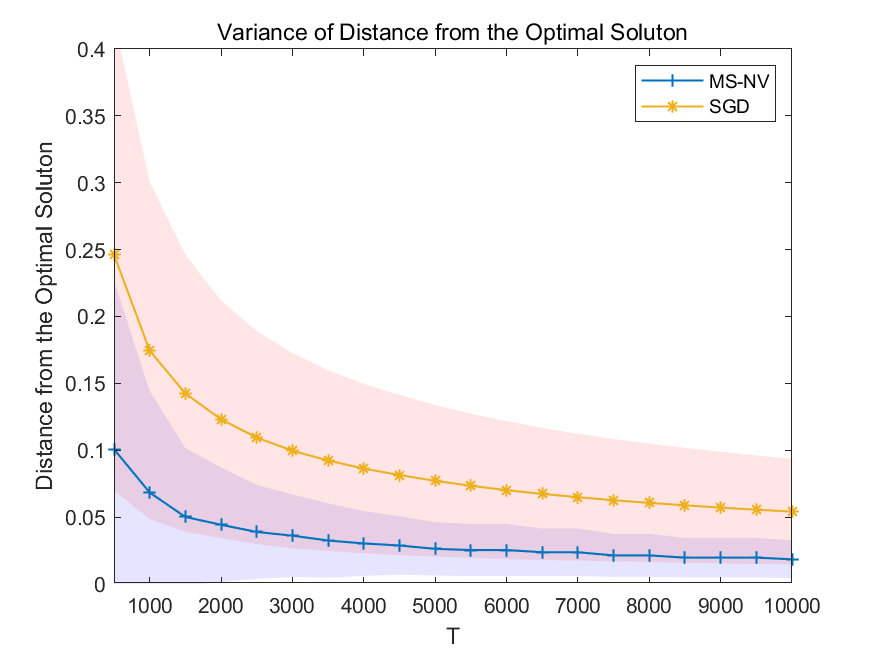}
\end{center}
\caption{{NVP: Variance of Regret and Distance}}\label{fig:nvp-var}
\end{figure}
}

{
\subsection{Application I: Multi-Product Multi-Constraint Inventory System}
\label{subsec:num-mpmc}
In this section, we conduct numerical experiments for both Multi-Product Single Constraint (MPSC) system management problem and Multi-Product Multi-Constraint (MPMC) system management problem. We will compare both the regret and running time of SGD-based algorithms and meta-policy.
We consider two stepsize patterns ($\eta/\sqrt{t}$ and $\eta/{t}$) and two types of batch size (linear batch size and exponential batch size).

The holding costs of all products are selected as $1$. Each product's lost-sales penalty cost $b_i$ is drawn uniformly from interval $[5,30]$. The constraint set is represented by $\Gamma = \{\by \mid A\by \leq c\cdot \bm{1}\}$. All elements of $A$ are drawn uniformly from $[0,1]$. We select the scale factor $c$ as the sum of all products' mean to get a reasonable constraint set.

\medskip \noindent {\bf The Multi-Product Single Constraint (MPSC) problem.} We consider an MPSC instance with three products.  We compare the SGD-based algorithm (DDM in~\cite{shi2016nonparametric})\footnote{{In each period DDM first conducts SGD algorithm to update the targeted order-up-to level. When the target order-up-to level updated by SGD is feasible, DDM sets this level as the actually implemented order-up-to level. Otherwise,  DDM uses a greedy projection algorithm to obtain the actually implemented order-up-to level.}} and our meta-policy with the greedy projection transition solver.

\begin{figure}[th]
\begin{center}
\includegraphics[width =0.48\textwidth]{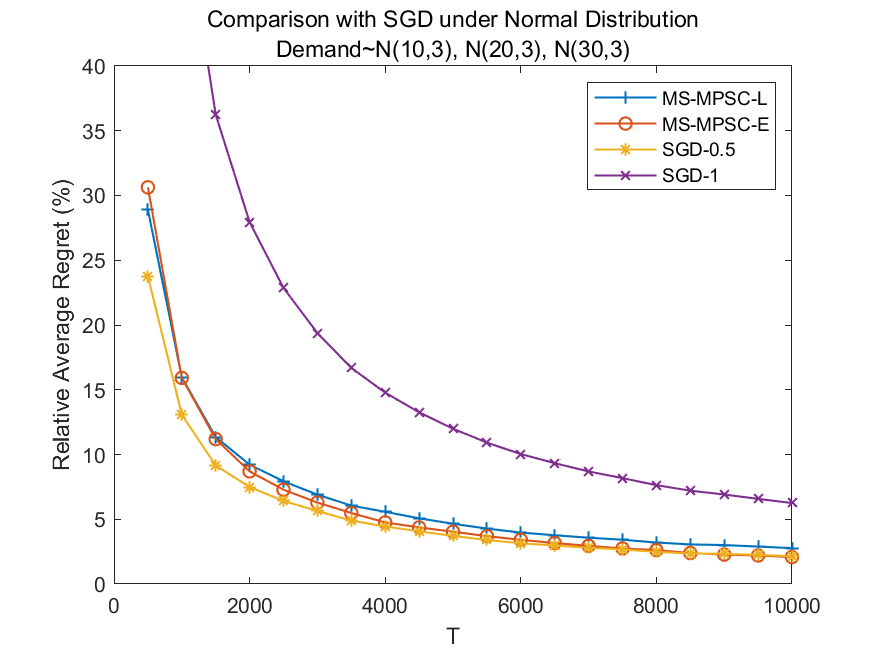}
\includegraphics[width =0.48\textwidth]{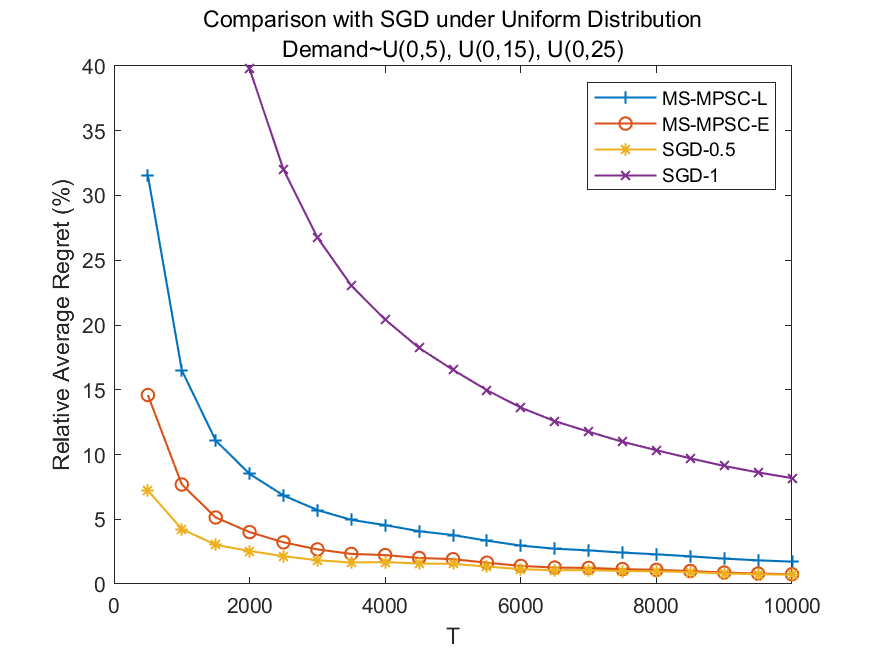}
\includegraphics[width =0.48\textwidth]{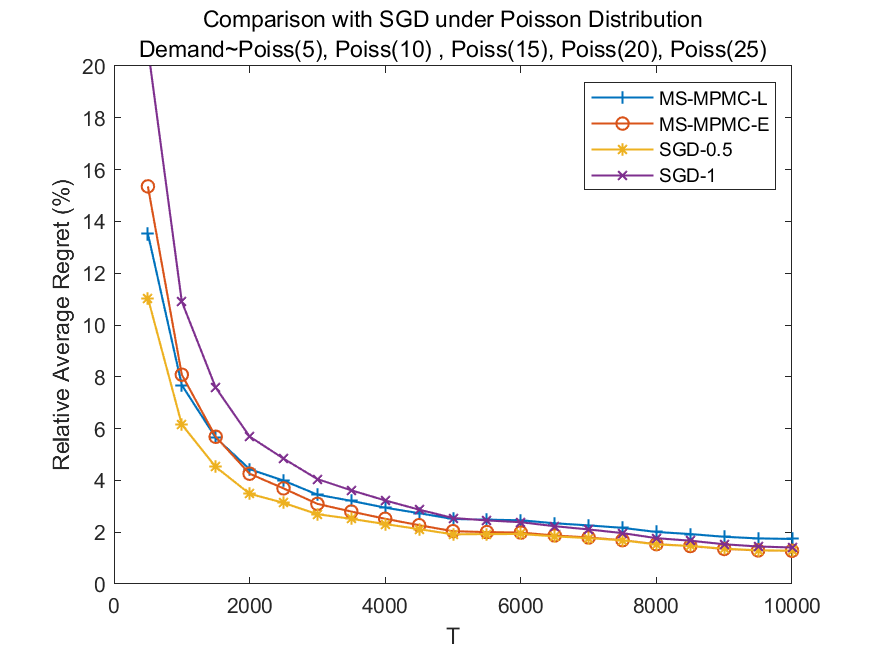}
\includegraphics[width =0.48\textwidth]{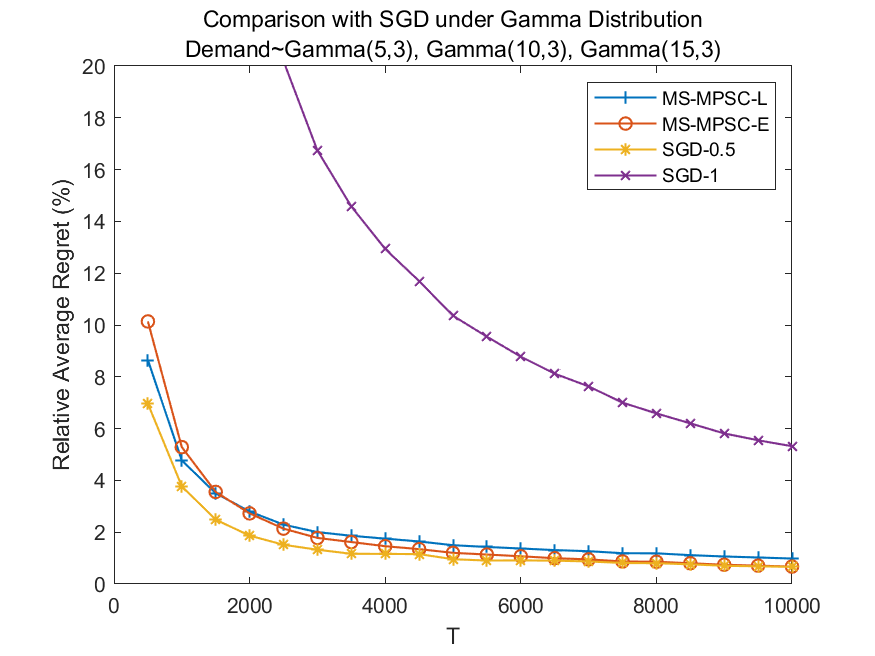}
\end{center}
\caption{{MPSC: Comparison with SGD under Different Distributions}}\label{fig:MPSC}
\end{figure}

\medskip \noindent {\bf The Multi-Product Multi-Constraint (MPMC) problem.} We consider an MPMC problem with three constraints and five products (i.e., $m=3$ and $n=5$). The SGD-based algorithm in this problem generalizes the projection step of \cite{shi2016nonparametric} to the multi-constraint version. 

\medskip 
The relative average regret curves of the MPSC and MPMC problems are presented in Figure~\ref{fig:MPSC} and Figure~\ref{fig:MPMC}, respectively.

\medskip \noindent \textbf{Comparison of regret performance.} Through the numerical results in both Figure~\ref{fig:MPSC} and Figure \ref{fig:MPMC}, we could note that the regret performances of our meta-policy and the SGD algorithm are about the same with each other.
For our meta-policy, the exponential batch size has a better performance than the linear batch size performance, which is in consistent with our theoretical results.
Regarding to the SGD algorithm, the stepsize pattern $\eta/\sqrt{t}$ has a better performance than the stepsize pattern $\eta/{t}$, and in some instances the stepsize pattern $\eta/t$ performs very poorly (e.g., the three instances in Figure~\ref{fig:MPSC} except the bottom-left one). 
However, as shown in \cite{shi2016nonparametric} theoretically, SGD with the stepsize pattern $\eta/t$ will achieve an $\mathcal O(\log T)$ regret, which is better than the $\mathcal O(\sqrt{T})$ regret of the stepsize pattern $\eta/\sqrt{t}$. 
Similar numerical phenomena have also been observed in other numerical experiments, for example, the numerical experiments in Sections~\ref{subsec:num-me} and~\ref{subsec:num-owms}. 
In our opinion, this is because the stepsize pattern $\eta/t$ converges to zero too fast.

\begin{figure}[th]
\begin{center}
\includegraphics[width =0.48\textwidth]{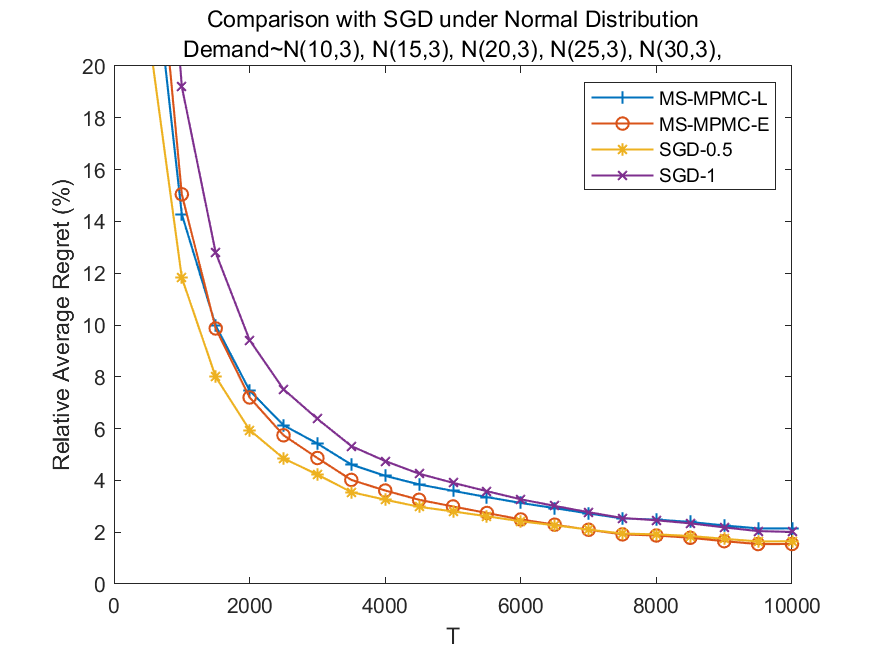}
\includegraphics[width =0.48\textwidth]{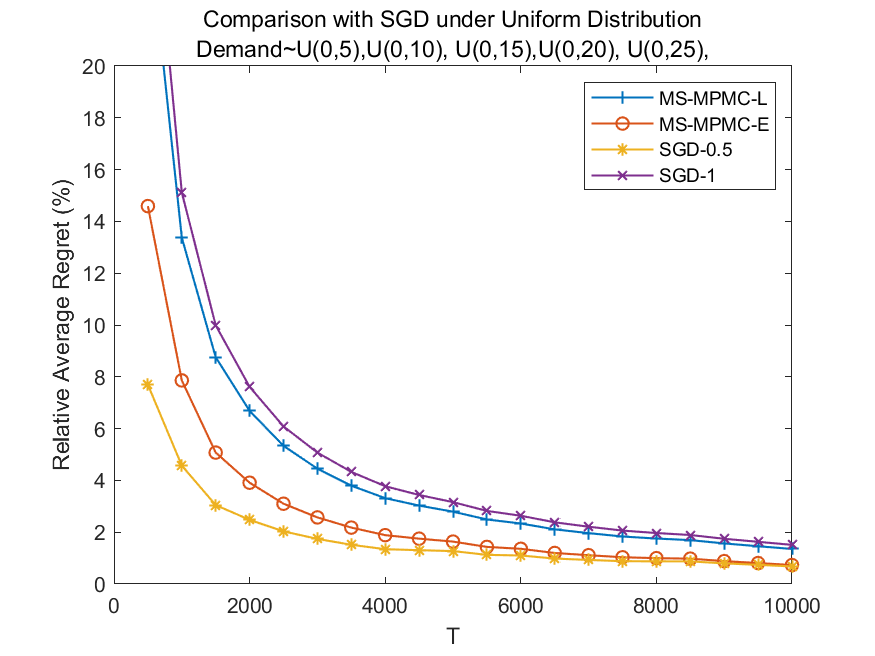}
\includegraphics[width =0.48\textwidth]{MPMC_poisson.png}
\includegraphics[width =0.48\textwidth]{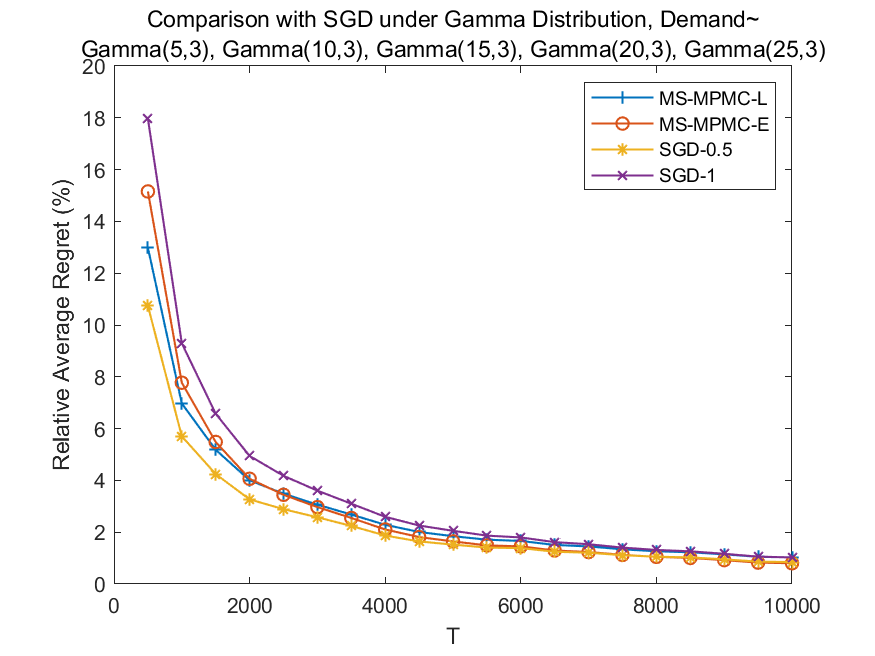}
\end{center}
\caption{{MPMC: Comparison with SGD under Different Distributions}}\label{fig:MPMC}
\end{figure}

\medskip \noindent {\bf Comparison of running time.} 
Each projection on the constraint set considered in the MPSC problem and MPMC problem requires solving a quadratic programming problem, which is quite time-consuming. 
SGD-based algorithms have to update their decision policies at every time step, while each update has to perform a projection operation due to the constrained gradient descent.
In contrast, our meta-policy greatly reduces the number of projection operations due to the infrequent decision updates. Therefore, theoretically,  an important strength of our meta-policy is its remarkable computational efficiency due to the infrequent decision updates. In the following, we verify this strength through experimental studies.
}

{
We consider the MPSC problem and MPMC problem as described before. The quadratic programming problems of both algorithms are solved by the function \texttt{fmincon()} of MATLAB. We test the running time of both algorithms for many different period numbers.
The results under the MPSC problem are presented in Table~\ref{tab:run1} and the results under the MPMC problem are presented in Table~\ref{tab:run2}. 
We can observe that our method achieves about $15\sim 150$ times speed-up for different values of $T$. 
These results indicate that our minibatch-SGD-based algorithm is much faster than the SGD-based algorithm and very computationally efficient.

\begin{table}[]
\centering
\caption{{MPSC: Running time (in seconds) compared with SGD}}\label{tab:run1}{
\begin{tabular}{l|cccccccccc}
\hline
Time periods & 500 & 1000 & 2000 & 3000 & 5000 & 7500 & 10000 & 15000 & 20000 & 50000 \\ \hline
SGD-0.5 & 12.12 & 21.81 & 41.25 & 60.32 & 98.24 & 143.32 & 188.77 & 283.3 & 375.38 & 914.35 \\ \hline
SGD-1 & 10.29 & 19.84 & 39.05 & 57.94 & 95.16 & 141.87 & 189.29 & 282.26 & 377.02 & 896.5 \\ \hline
MS-MPSC-L & 0.65 & 0.94 & 1.34 & 1.67 & 2.16 & 2.63 & 3.05 & 3.8 & 4.45 & 7.48 \\ \hline
MS-MPSC-E & 0.63 & 0.74 & 0.89 & 0.99 & 1.14 & 1.26 & 1.38 & 1.57 & 1.76 & 2.74 \\ \hline
\end{tabular}}
\end{table}

\begin{table}[]
\centering
\caption{{MPMC: Running time (in seconds) compared with SGD}}\label{tab:run2}{
\begin{tabular}{l|cccccccccc}
\hline
Time periods & 500 & 1000 & 2000 & 3000 & 5000 & 7500 & 10000 & 15000 & 20000 & 50000 \\ \hline
SGD-0.5 & 12.99 & 25.31 & 52.19 & 76.57 & 129.1 & 198.14 & 266.72 & 417.47 & 561.92 & 1387.47 \\ \hline
SGD-1 & 10.09 & 22.31 & 48.86 & 72.98 & 125.36 & 200.16 & 278.36 & 434.25 & 575.31 & 1412.35 \\ \hline
MS-MPMC-L & 0.5 & 0.65 & 0.88 & 1.04 & 1.33 & 1.62 & 1.86 & 2.32 & 2.73 & 4.81 \\ \hline
MS-MPMC-E & 0.39 & 0.56 & 0.78 & 0.94 & 1.24 & 1.51 & 1.74 & 2.21 & 2.65 & 4.66 \\ \hline
\end{tabular}}
\end{table}

{In Figure~\ref{fig:MPSC} and Figure~\ref{fig:MPMC}, it turns out that the minibatch SGD with exponential and linear batch sizes performs quite similarly across all distributions except for the uniform distribution.  We think the possible reasons are as follows.
\begin{enumerate}
    \item {\bf Lower bound of density.} By the regret analysis under strong convexity (Theorem \ref{thm:inv,str}), we know that regret of order $\mathcal{O}(\log T /\alpha)$ and would decrease if parameter $\alpha$ becomes larger, where $\alpha$ is a positive lower bound on the density of demand $f_D(x)$.
    Note the density of the uniform distribution is constant over its domain. Therefore, compared with other distributions, the uniform distribution has the largest $\alpha$ and the minibatch SGD with exponential batch size has the best theoretical guarantee. Compared with the minibatch SGD with linear batch size, the difference caused by $\alpha$ is maximized by the uniform distribution. Therefore, under the uniform distribution, the theoretical performance gap between the minibatch SGD with two different batch sizes is the largest.
    \item {\bf Tail behavior of distribution.} We take the single product case as an example. Note the optimal solution of the single product case is the newsvendor solution $x^* = F^{-1}(b/(h+b))$. As $b$ is large compared with $h$, $x^*$ will be close to the upper bound  $u$ of demand under uniform distribution. Therefore, if a solution $x_{\tau}$ exceeds $x^*$ by a small distance greater than $u-x^*$, due to the boundedness of $D$, the gradient estimator will be positive, and in the next iteration $x_{\tau}$ must decrease. That is the overestimation of $x^*$ can be corrected easily. Therefore, the frequent updates in the early periods will find the $x^*$ more quickly and not need to pay much cost. As the minibatch SGD with exponential batch size in the above experiments has a base of $1.05$ or $1.15$, it updates more frequently than the minibatch SGD with linear batch size. Thus, the minibatch SGD with the exponential batch size is more likely to perform better under the uniform distribution.
    \end{enumerate}}

\subsection{Application II: Multi-Echelon Serial Inventory System}
In this section, we consider the Multi-Echelon (ME) serial inventory system management problem described in Section~\ref{app: multi-echelon} with three echelons. We consider four different demand distributions ($N(5,1)$, $U(0,10)$, $\mathrm{Poiss}(5)$, and $\mathrm{Geom}(0.2)$) and set $\bh=(1,1,1)$ and $\bb=(50,50,50)$. We compare the performance of SGD under different stepsizes and meta-policy under different batch sizes. 
The numerical results are presented in Figure~\ref{fig:multi-echelon}.
\label{subsec:num-me}
\begin{figure}[th]
\begin{center}
\includegraphics[width =0.48\textwidth]{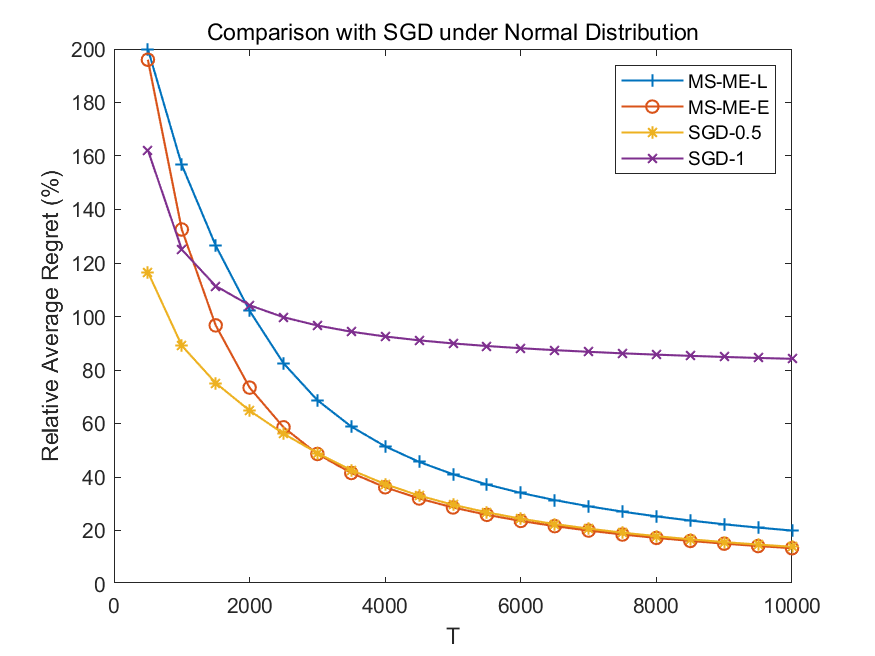}
\includegraphics[width =0.48\textwidth]{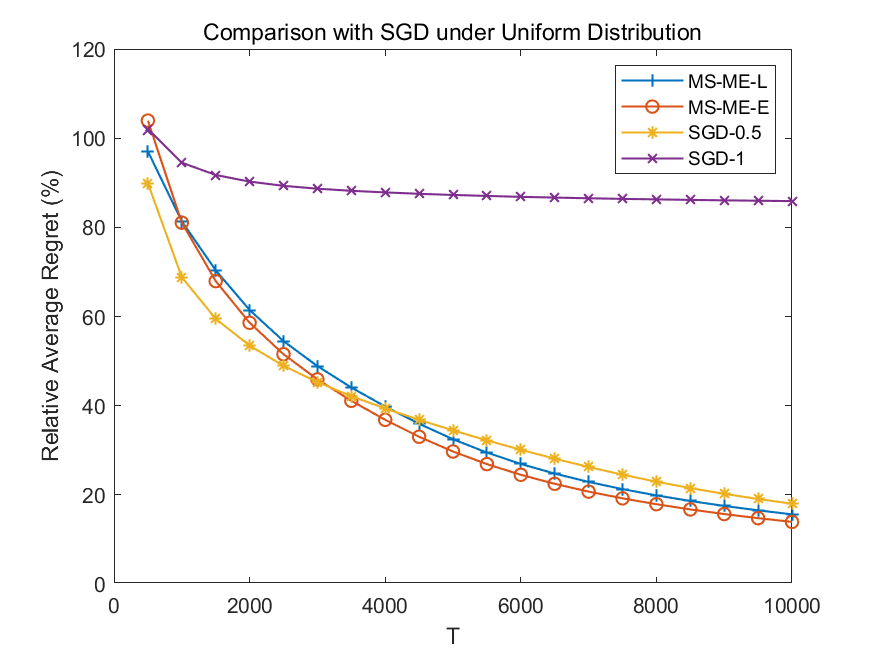}
\includegraphics[width =0.48\textwidth]{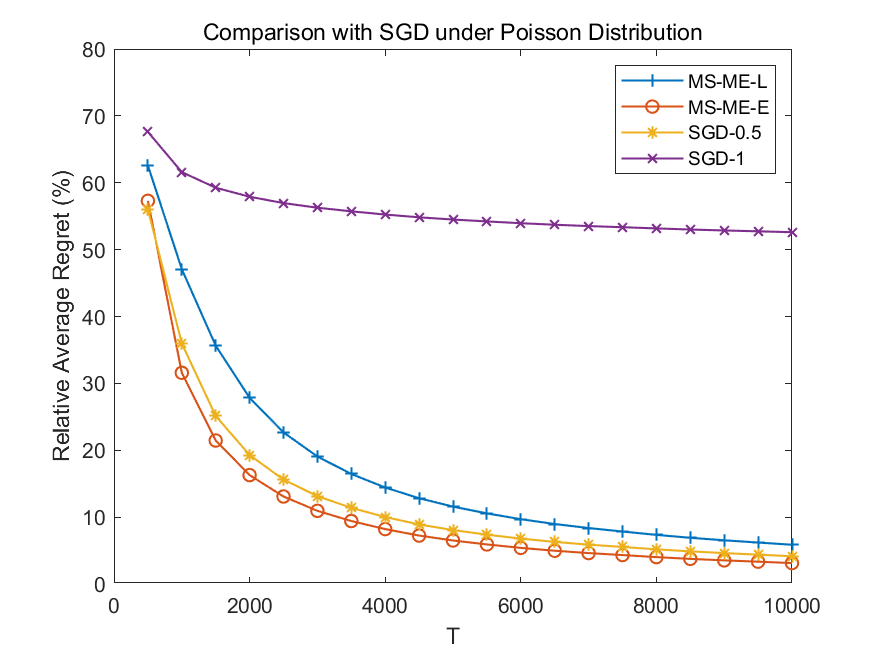}
\includegraphics[width =0.48\textwidth]{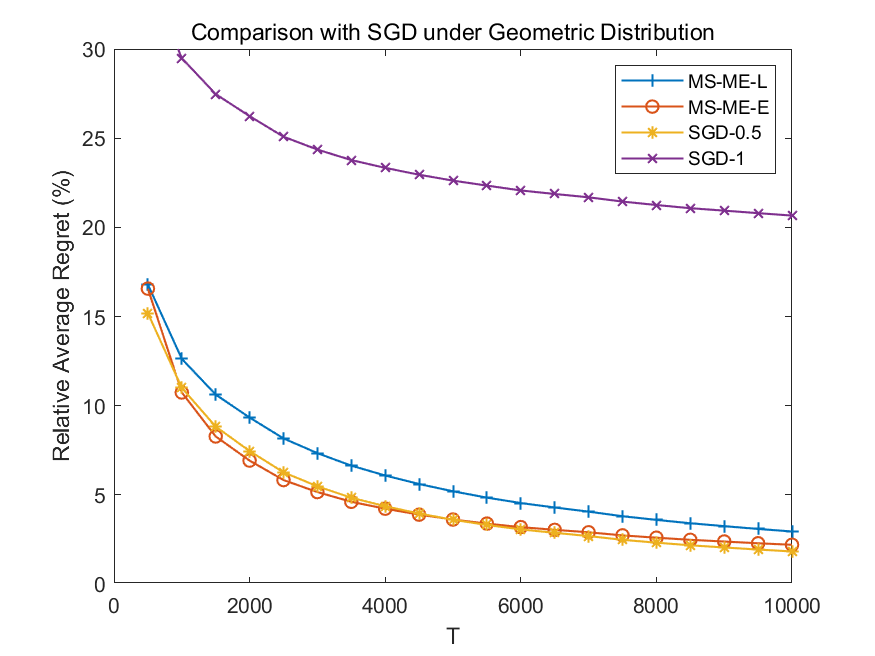}
\end{center}
\caption{{Multi-Echelon: Comparison with Different Algorithms under Different Distributions}}\label{fig:multi-echelon}
\end{figure}
\subsection{Application III: One-Warehouse Multi-Store Inventory System}
In this section, we consider the One-Warehouse Multi-Store (OWMS) system management problem described in Section~\ref{app: owms} with three stores. We consider four different demand distributions ($N(5,1)$, $U(0,10)$, $\mathrm{Poiss}(5)$, and $\mathrm{Geom}(0.2)$) and set $\bh=(0.5,1,1,1)$, $\bb=(70,50,30)$ and $\bm{c} = (10,20,30)$. We compare the performance of SGD under different stepsizes and meta-policy under different batch sizes. 
The numerical results are presented in Figure~\ref{fig:owms-exp}. 
\label{subsec:num-owms}
\begin{figure}[th]
\begin{center}
\includegraphics[width =0.48\textwidth]{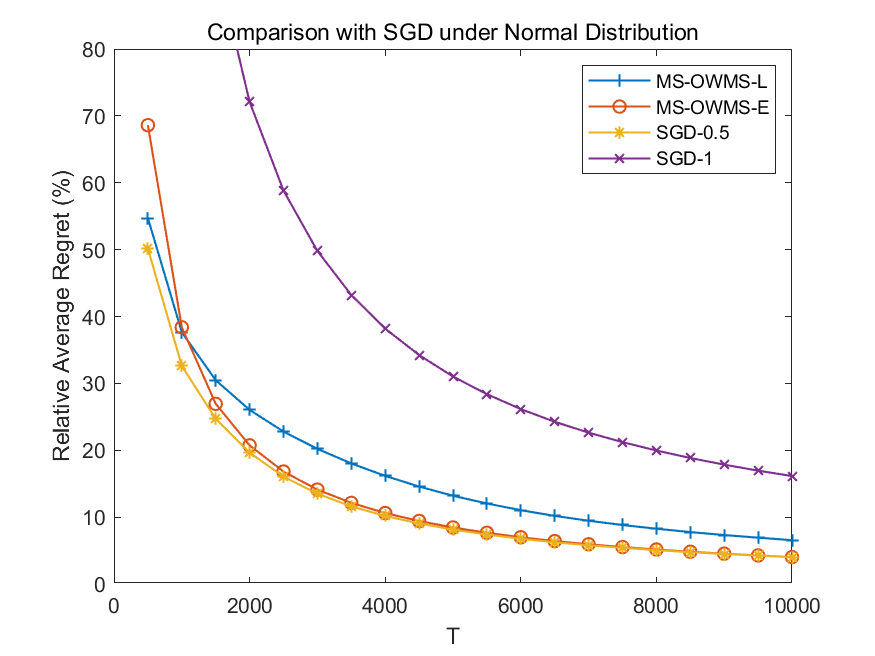}
\includegraphics[width =0.48\textwidth]{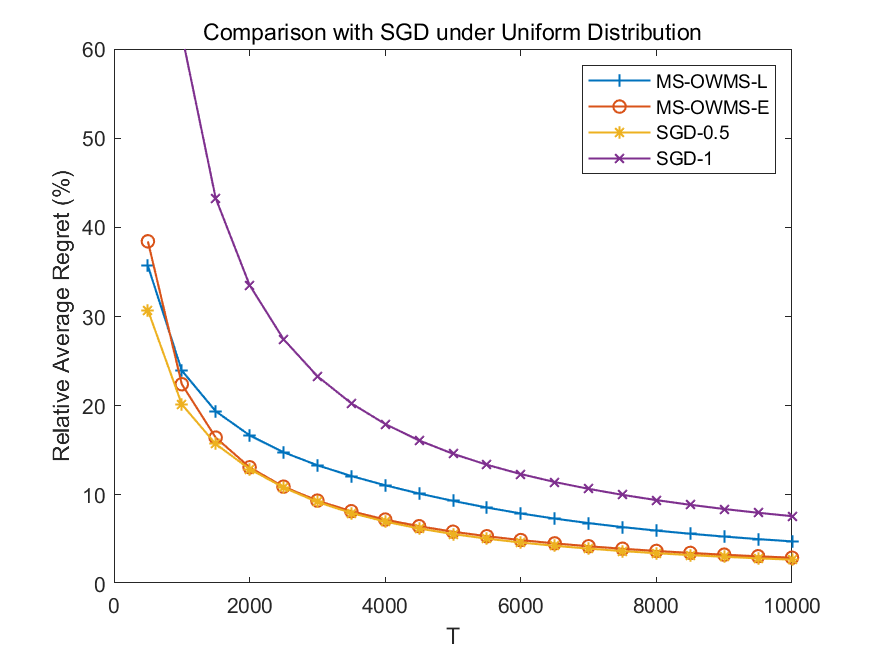}
\includegraphics[width =0.48\textwidth]{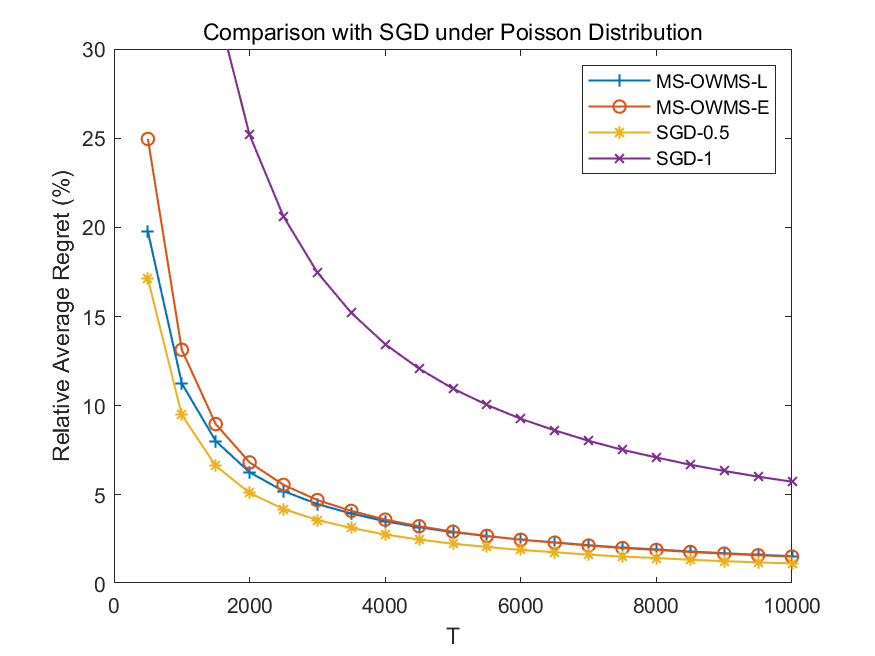}
\includegraphics[width =0.48\textwidth]{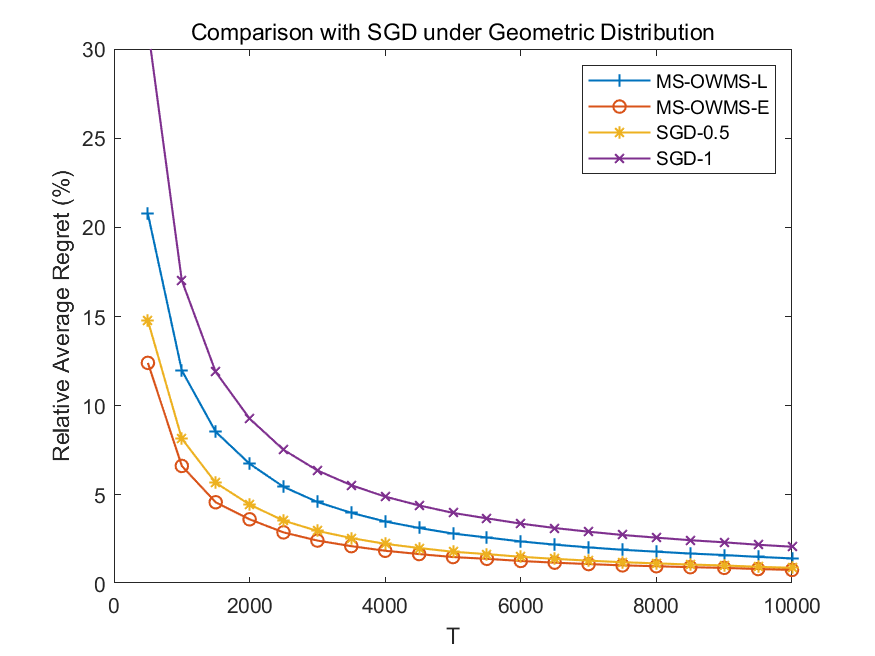}
\end{center}
\caption{{OWMS: Comparison with SGD under Different Distributions}}\label{fig:owms-exp}
\end{figure}
\subsection{Inventory System with Positive Lead Times}
\label{subsec:num-leadtime}
In this section, we test the performance of our minibatch-SGD-based meta-policy in learning the optimal base-stock policy of lost-sales inventory systems with positive lead times. The lost-sales system has holding cost $h$, lost-sales cost $b$, and deterministic lead time $L$.

We compute the sample path derivative of the cost function (see \cite{huh2009adaptive} for detailed expression) as a gradient estimator of the long-run cost function. The batch size is linear and the transition solver is replaced by a fixed warm-up stage of $\log T$ periods. The algorithm is named MS-LT (Minibatch-SGD-based (MS) meta-policy for inventory systems with Lead Times (LT)).
}

{
To test whether our algorithm tailored to this problem is efficient and invites further theoretical research, we consider various problem parameters as follows.
\begin{enumerate}[nolistsep]
    \item \textbf{Different lost-sales costs:} four lost-sales costs $b = 10,25,50,100$ are tested.
    \item \textbf{Different demand distributions:} four different distributions are tested ($N(5,1)$, $U(0,10)$, $\mathrm{Poiss}(5)$, and $\mathrm{Geom}(0.2)$).
    \item \textbf{Different lead times:} four different lead times $L = 5,10,15,20$ are tested.
    \item \textbf{Different variances of demand:} clipped normal distributions with a mean of $20$ and four different standard deviations of $\sigma = 1,2,3,4$ are tested.
\end{enumerate}

We use the default problem parameter setting: $b=25$, normal distribution $N(5,1)$, and $L = 10$. For each of the problem parameters, we test the parameter by varying its value according to the above list, while keeping the other parameters as default. The relative average regret curves are reported in Figure~\ref{fig:leadtime}.

As shown in the figures, the relative average regret of MS-LT converges to zero consistently under different lost-sales costs, distributions, lead times, and variances. We can also observe the influence of problem parameters on regret. For example, larger lost-sales costs and lead times cause larger regret. One can also observe in the bottom-right figure that a larger variance causes a lower relative regret. We think it is mainly because the larger variance will cause a much larger optimal cost compared to the regret. 
\begin{figure}[th]
\begin{center}
\includegraphics[width =0.48\textwidth]{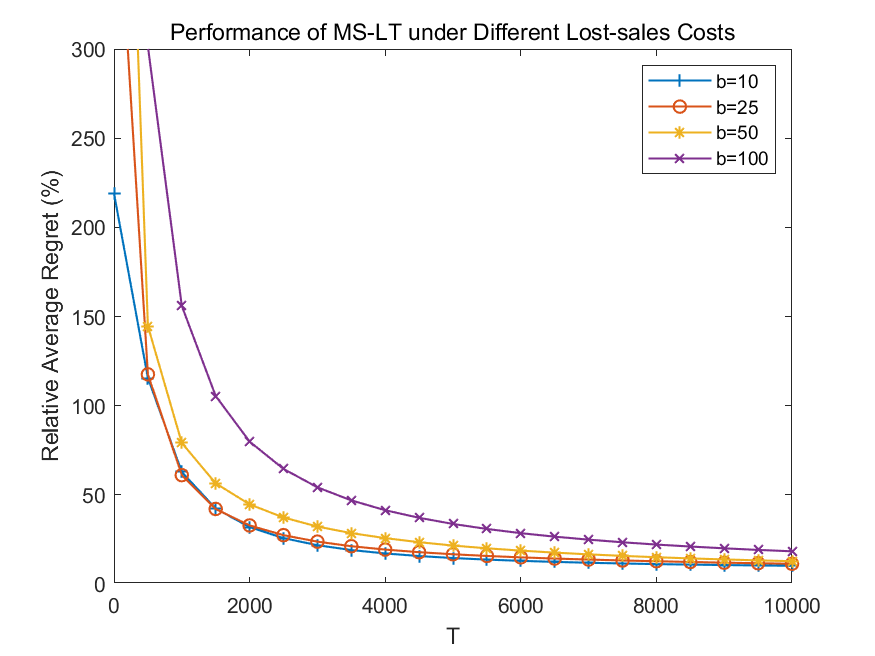}
\includegraphics[width =0.48\textwidth]{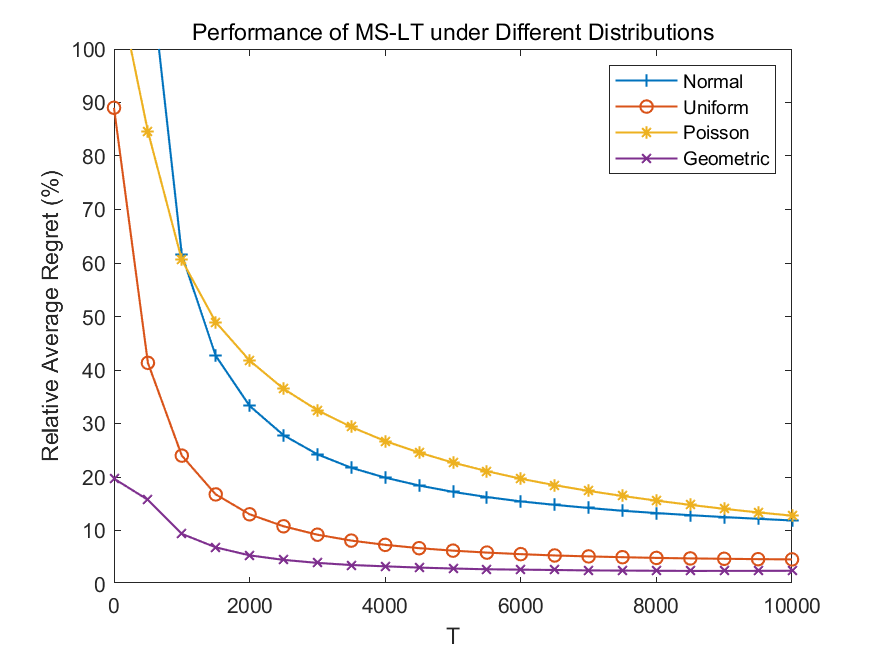}
\includegraphics[width =0.48\textwidth]{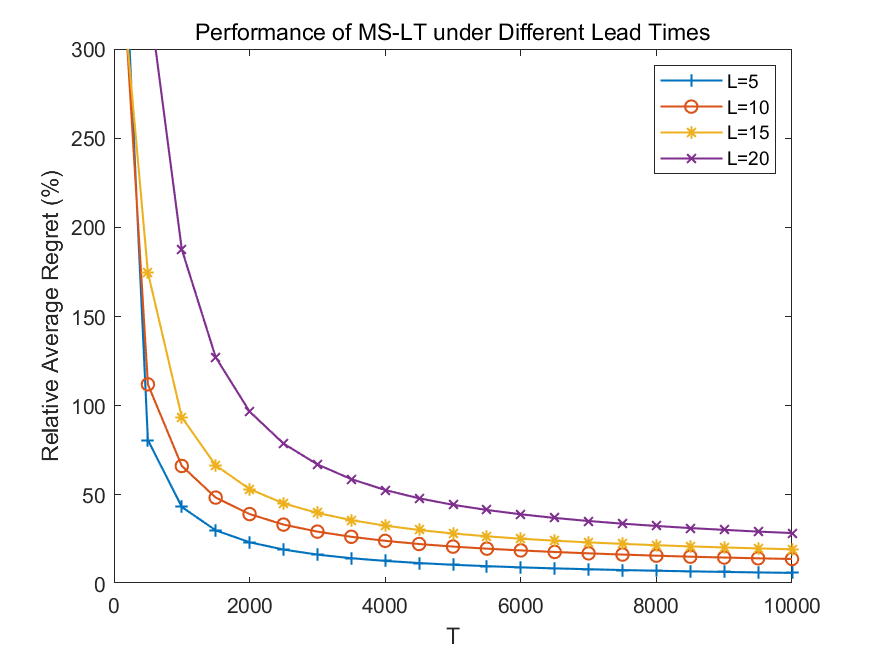}
\includegraphics[width =0.48\textwidth]{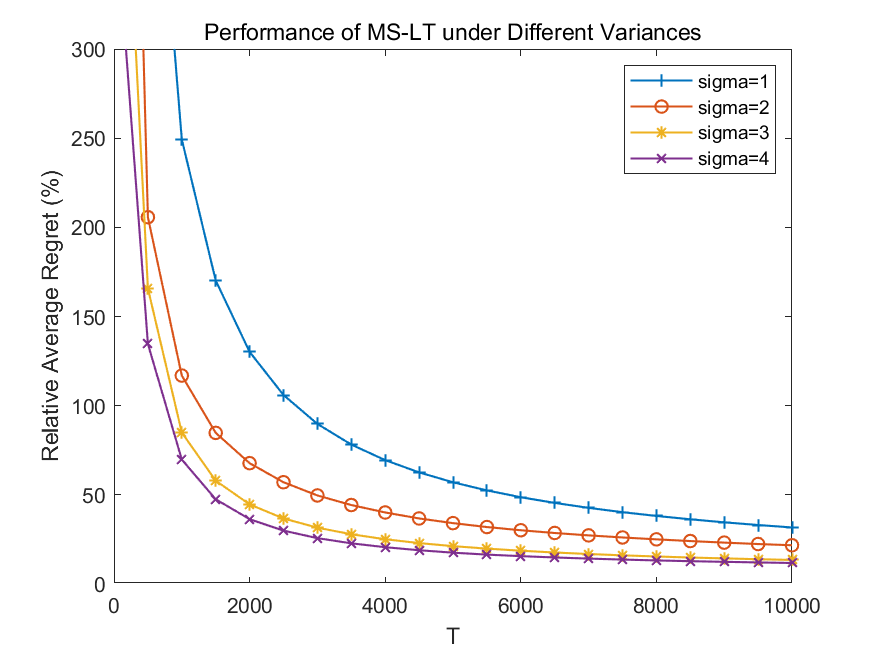}
\end{center}
\caption{{Performance of MS-LT under Different Problem Settings}}\label{fig:leadtime}
\end{figure}

\subsection{Two-Echelon Inventory System in~\cite{zhang2022no}}
\label{subsec:num-te}
In this section, we consider the centralized Two-Echelon (TE) problem studied in~\cite{zhang2022no}.  In their Algorithm 1, the authors used a batched SAA algorithm to estimate the optimal decision of the retailer, and under the estimated decision of the retailer they used a constant-stepsize SGD algorithm (Algorithm 2 therein) to compute the optimal solution of the supplier. In the following, we inherit the formulations and notations in ~\cite{zhang2022no}, please refer to their paper for detailed definitions.
}

{
To study the efficacy of our minibatch-SGD-based meta-policy, we consider adapting our meta-policy to this problem in the following manner. For the decision of the retailer, we use a similar SAA method to estimate the decision of the retailer, then we apply the minibatch-SGD to optimize the solution of the supplier. The central planner is presented in Algorithm~\ref{alg:plan-agent}.

{Algorithm \ref{alg:plan-agent} precedes in epochs. In epoch $m$ of Algorithm \ref{alg:plan-agent}, the central planner collects data set $\mathcal{D}_m$ and calls Algorithm \ref{alg:agent2} with input $\mathcal{D}_m$ to obtain decision of Agents 1 and 2 as $s_{m+1,1}$ and $s_{m+1,2}$. In Algorithm \ref{alg:agent2}, $s_{m+1,1}$ is solved by SAA and $s_{m+1,2}$ is given by minibatch SGD. Specifically, we use data set $\mathcal{D}_m$ to run a minibatch SGD. For each time $L/2+1 \leq i \leq L$, the minibatch SGD gives an iteration point $s_{m+1,2,i}$ and we take the average to obtain $s_{m+1,2}$.}

\begin{algorithm}[t!]
\caption{Central Planner for Agent 1 and Agent 2}
\label{alg:plan-agent}
{\begin{algorithmic}[1]
\State \textbf{Input: }Algorithm $\mathcal{A}$ (Algorithm \ref{alg:agent2}).
\State \textbf{Initialize: } Set $s_{1,1}$ and $s_{1,2}$ arbitrarily. Epoch length $L_1 = 2$, $t = 1$.
\For{$m = 1,2,\dots$}
\State Define $I_m = \{t,t+1,\dots,t+L_m-1\}$.

\While{$t \in I_m$}
\State Decide the desired inventory level: $s_{m,1}$ for Agent 1 and $s_{m,2}$ for Agent 2.
\State Receive the realized demand $d_t$, $t \leftarrow t+1$.
\EndWhile
\State \textbf{end}
\State Collect $\mathcal{D}_m =\{d_i\}_{i \in I_m}$; send $\mathcal{D}_m$ to $\mathcal{A}$ and get $s_{m+1,1}$ and $s_{m+1,2}$; and set $L_{m+1}=2L_m$.
\EndFor
\State \textbf{end} 
\end{algorithmic}}
\end{algorithm} 

\begin{algorithm}[t!]
\caption{Optimization Algorithm for Agent 1 and Agent 2}
\label{alg:agent2}
{\begin{algorithmic}[1]
\State \textbf{Input: }A set of demand value $\mathcal{D}_m = \{ d_1,\dots,d_L\}$.
\State \textbf{Initialize: } Set a minibatch size sequence $\{n_1, n_2,\dots,n_{\tau_{\mathrm{max}}}\}$ where $\tau_{\mathrm{max}} = \min\{k:\sum_{\tau = 1}^k n_{\tau} \geq T\}$  and stepsize $\eta$. Set a working period counter $l=0$, epoch step $\tau=1$, estimated gradient sequence $\mathcal G = \emptyset$ for the updating of minibatch SGD. Set $w_1 \leq D - \frac{h_2}{\Gamma(h_2+p_1)} = s_{\mathrm{max}}$ arbitrarily.

\State Compute the empirical cumulative density function $\hat{\Phi}_{L}(x) = \frac{2}{L}\sum_{i=1}^{L/2} \mathbb{I}\{ d_i \leq x\}$.

\State Set $s_1 = \hat{\Phi}^{-1}_{L}(\frac{h_2+p_1}{h_1+p_1})$.

\For{$i = L/2+1,L/2+2,\dots,L$}

\State Set $s_{m+1,2,i} = w_{\tau}$ and let $l = l+1$.
\State Compute $m_i = \mathbb{I}[s_{m+1,2,i} \leq d_{i-1}] \cdot [(h_1+p_1)\mathbb{I}[\hat{s}_{i,1} \geq d_i]-p_1] + h_2 \mathbb{I}[s_{m+1,2,i} \geq d_i] + C_1 (h_1+p_1)\sqrt{\frac{2\log(T^3D)}{L}} \cdot \hat{\Phi}_{L}(s_{m+1,2,i})$ and $\hat{s}_{i,1} = s_1$ if $d_{i-1} \leq s_{m+1,2,i}$ and $\hat{s}_{i,1} = s_1 + s_{m+1,2,i} - d_{i-1}$ otherwise.

\State Add $m_i$ to the estimated gradient sequence $\mathcal G$.
    \If{ $l = n_\tau$} 
    \State Update $\bw_{\tau+1}$ by minibatch SGD: $\displaystyle{
    w_{\tau+1}=\min \left\{ s_{\mathrm{max}}, \max\left\{0,w_{\tau}- \frac{\eta}{n_\tau} \sum_{\hat g\in\mathcal G}\hat g ]\right\}\right\}}$. 
    \State Set $l=0$, $\tau=\tau+1$ and $\mathcal G = \emptyset$. 
    \EndIf

\EndFor
\State \textbf{end for}

\State \textbf{return} $s_{m+1,1} =s_1$ and ${s}_{m+1,2} = \frac{2}{L}\sum_{i=L/2+1}^L s_{m+1,2,i}$.
\end{algorithmic}}
\end{algorithm} 

\medskip \noindent {\bf Numerical results.} We compare the SGD-based method of \cite{zhang2022no} with the minibatch-SGD-based method under two demand distributions ($N(5,1)$ and $U(0,10)$) and set $(h_1,h_2)=(2,1)$ and $(b_1,b_2)=(50,50)$. 
The numerical results are presented in Figure~\ref{fig:twoechelon}. Through the numerical results, we observe that under both distributions, both versions of MS-TE perform better than the SGD-based method of \cite{zhang2022no}.

\begin{figure}[th]
\begin{center}
\includegraphics[width =0.48\textwidth]{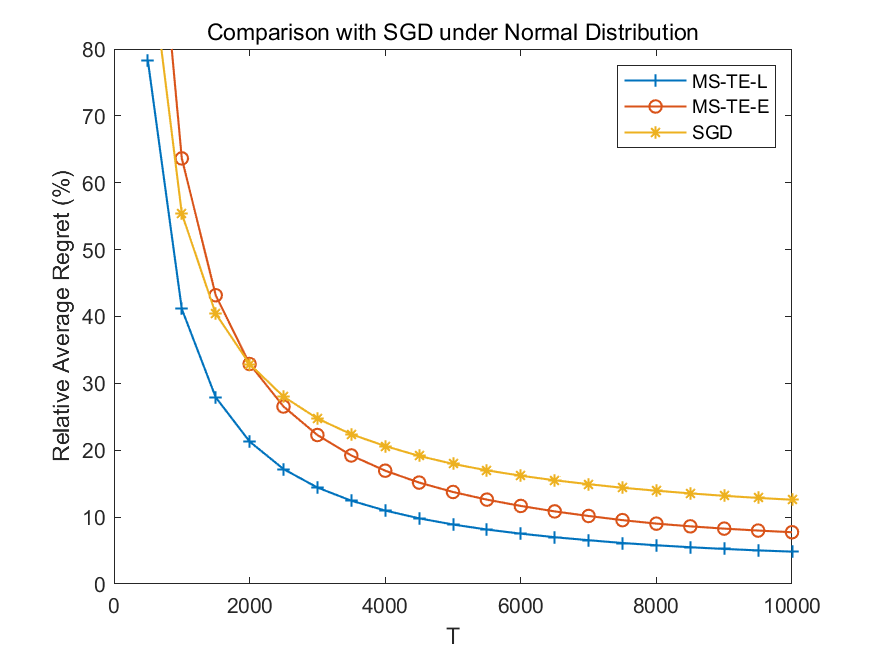}
\includegraphics[width =0.48\textwidth]{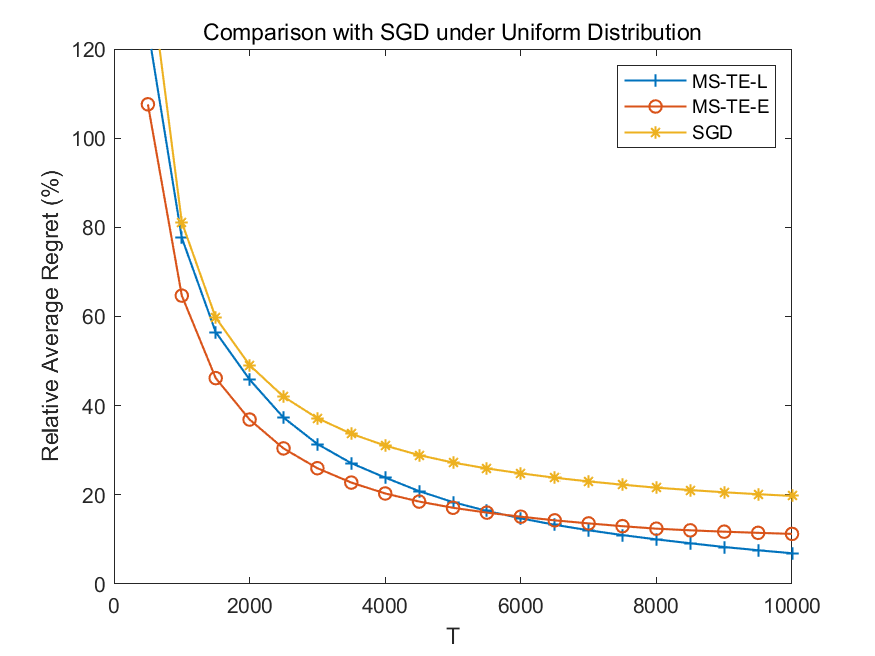}
\end{center}
\caption{{Two-Echelon: Comparison with SGD under Different Distributions}}\label{fig:twoechelon}
\end{figure}

\medskip \noindent {\bf Theoretical results.} Motivated by the above numerical results, we provide a regret analysis for our Algorithm~\ref{alg:plan-agent}. The main difference between Algorithm~\ref{alg:plan-agent} and Algorithm 1 of \cite{zhang2022no} is that we use the minibatch-SGD to optimize the decision of Agent 2. Therefore, we will re-use many results in~\cite{zhang2022no} to simplify the presentation.

\begin{theorem}
Run Algorithm \ref{alg:plan-agent} with appropriate stepsize and batch size $n_{\tau} = K\tau$, where $K$ is a positive integer, we have $\mathcal{R}(T) \leq \tilde{O}(\sqrt{T})$.
\end{theorem}
\proof{Proof.}
Recall the definition of regret is
\[
\mathcal{R}(T) = \E \left[\sum_{t=1}^T \tilde{H}_t -\sum_{t=1}^T {H}(s_1^*,s_2^*)  \right].
\]
By the low-switching property of Algorithm~\ref{alg:plan-agent}, it suffices to provide an upper bound on
\[
\E \left[\sum_{m=1}^M \left( \sum_{t \in I_m} H(s_{m,1},s_{m,2}) -\sum_{t \in I_m} H(s_1^*,s_2^*) \right)  \right],
\]
where detailed proof can be found in Eq. (19) of~\cite{zhang2022no}.

Now we focus on the analysis of an epoch $I_m$ with length $L$ for a fixed $m \in [M]$. To establish an upper bound on the above term, we define
\[
H'_L(s_{1},s_2) = H(s_1,s_2) + C_1 (h_1+p_1)\sqrt{\frac{2\log(T^3D)}{L}} \int_0^{s_2} \Phi(x) \mathrm{d}x.
\]

From the assumption on the demand distribution, it is easy to verify that $H'_L(s_{m,1},\cdot)$ is $\beta$-smooth (we omit the dependence on constants). Besides, $m_i$ is bounded and an unbiased gradient estimator of the function $H'_L(s_{m,1},\cdot)$. By an argument similar to the proof of Lemma A.5 of \cite{zhang2022no}, we know $H'_L(s_{m,1},\cdot)$ is convex with probability\footnote{{Note that we divide the data set $\mathcal{D}$ into two parts, the first part is used to solve $s_{m,1}$ and the second part is used to optimize $s_{m,2}$. The high probability statement holds for the first part and the following expectation is taken with respect to the second part.}} at least $1-1/T^2$. Choose stepsize $\eta < 1/\beta$, from the regret analysis of minibatch-SGD (Lemma \ref{le:cvx}), we have
\[
\E \left[\sum_{i=L/2+1}^L H'_L(s_{m,1},s_{m,2,i}) - \sum_{i=L/2+1}^L H'_L(s_{m,1},s_{2}) \right] \leq \tilde{O}(\sqrt{L}),
\]
for any $s_2$. Therefore, we have
\begin{align*}
    &\E \left[\sum_{i=L/2+1}^L H'_L(s_1^*,s_{m,2,i}) - \sum_{i=L/2+1}^L H'_L(s_1^*,s_{2}) \right] \\
    &\qquad\qquad \leq \E \left[\sum_{i=L/2+1}^L H'_L(s_{m,1},s_{m,2,i}) - \sum_{i=L/2+1}^L H'_L(s_{m,1},s_{2}) \right] + \O(\sqrt{L \log(T^3 D)}),
\end{align*}
where we use the Lipschitz property of $H'$ and convergence result on $s_{m,1}$ (detailed proof based on SAA can be found in Equation (11) of \cite{zhang2022no}).

Choose $s_2 = s_2^*$ and using the fact that $H'_L(s_1,s_2)$ is a small perturbation of $H(s_1,s_2)$, we have
\begin{align*}
    &\E\left[\sum_{i=L/2+1}^L H(s^*_{1},s_{m,2,i}) -\sum_{i=L/2+1}^L H(s_1^*,s_2^*) \right]\leq \E\left[\sum_{i=L/2+1}^L H_L'(s^*_{1},s_{m,2,i}) -\sum_{i=L/2+1}^L H'(s_1^*,s_2^*) \right] + \tilde{\O}(\sqrt{L}).
\end{align*}
From the convexity of $H(s_1^*,\cdot)$ and $s_{m,2} = \frac{2}{L} \sum_{i=L/2+1}^L s_{m,2,i}$, we have
\begin{align*}
    &L_m \cdot \E[H(s_1^*,s_{m,2}) -H(s_1^*,s_2^*)] \leq 2\E\left[\sum_{i=L/2+1}^L H(s^*_{1},s_{m,2,i}) -\sum_{i=L/2+1}^L H(s_1^*,s_2^*) \right] \leq \tilde{\O}(\sqrt{L_m}).
\end{align*}
Therefore, we have
\[
\E \left[\sum_{m=1}^M \left( \sum_{t \in I_m} H(s_{m,1},s_{m,2}) -\sum_{t \in I_m} H(s_1^*,s_2^*) \right)  \right] \leq \sum_{m=1}^M \tilde{\O} (\sqrt{L_m}) = \tilde{\O} (\sqrt{T}),
\]
where we complete the proof.
\Halmos
\endproof}

\end{document}